\documentclass[11pt]{article}
\usepackage{curves}
\usepackage{amsmath, amssymb, lamsarrow}
\usepackage{amssymb}
\def\mineappendix{
        \setcounter{section}{1}
        \setcounter{subsection}{0}
        \def\thesection{\Alph{section}}
        \def\sectionap{\@startsection  {section}{1}{\z@}
                        {-3.5ex plus-1ex minus-.2ex} {0ex plus.2ex}
                        {\reset@font\Large\bf  Appendix:  \, }
                        }
        }
\makeatother
\def\Proclaim #1. #2\par{\bigbreak\noindent{\sc#1.\enspace}{\it#2}\par}

\newcommand{\gwii}[1]{\left< \hspace{-2pt} \left< \, #1 \,
        \right>  \hspace{-2pt} \right>_{0}}

\newcommand{\gwiione}[1]{\left< \hspace{-2pt} \left< \, #1 \,
        \right> \hspace{-2pt} \right>_{1}}

\newcommand{\gwiitwo}[1]{\left< \hspace{-2pt} \left< \, #1 \,
        \right> \hspace{-2pt} \right>_{2}}

\newcommand{\gwih}[2]{ \left< \, #2 \, \right>_{#1}}
\newcommand{\gwiig}[1]{\left< \hspace{-2pt} \left< \, #1 \,
    \right> \hspace{-2pt} \right>_{g}}
\newcommand{\gwiih}[2]{\left< \hspace{-2pt} \left< \, #2 \,
    \right> \hspace{-2pt} \right>_{#1}}

\newcommand{\grav}[2]{\tau_{#1}(\gamma_{#2})}
\newcommand{\grava}[1]{\tau_{#1}(\gamma_{\alpha})}

\newcommand{\ga}{\gamma_{\alpha}}
\newcommand{\gua}{\gamma^{\alpha}}
\newcommand{\gb}{\gamma_{\beta}}
\newcommand{\gub}{\gamma^{\beta}}
\newcommand{\gm}{\gamma_{\mu}}
\newcommand{\gum}{\gamma^{\mu}}

\newcommand{\vs}{{\cal S}}

\newcommand{\vd}{{\cal D}}
\newcommand{\vw}{{\cal W}}

\newcommand{\qp}{\circ}

\newtheorem{lem}{Lemma}[section]

\newtheorem{thm}[lem]{Theorem}

\newtheorem{rem}[lem]{Remark}

\title{Topological Recursion Relations on $\overline{\cal M}_{3,2}$}

\author{Takashi Kimura
\thanks{Research of the first author was partially supported by NSA grant H98230-10-1-0179.}
, \,\,\, \,\,\,
Xiaobo Liu
\thanks{Research of the second author was partially supported by NSF grant
DMS-0905227, NSFC Tianyuan special fund 11326023, SRFDP grant
20120001110051, and Peking University 985 fund.} }

\date{}

\begin{document}

\maketitle
\allowdisplaybreaks

\begin{abstract}
In this paper, we give some new genus-3 universal equations for
Gromov-Witten invariants of compact symplectic manifolds. These
equations were obtained by studying relations in the tautological ring
of the moduli space of 2-pointed genus-3 stable curves. A byproduct of our search
for genus-3  equations is a new genus-2 universal equation
for Gromov-Witten invariants.
\end{abstract}

It is well known that relations in the tautological ring of moduli space of stable
curves $\overline{\cal M}_{g, n}$ produce universal equations for Gromov-Witten invariants
of all compact symplectic manifolds. A typical genus-0 example is the associativity equation for
the quantum cohomology, also known as the WDVV equation. Finding explicit
higher genus universal
equations is a very difficult problem. Some genus-1 and genus-2
universal equations were discovered in \cite{Ge1}, \cite{Ge2}, and
\cite{BP}. Relations among known genus-2 equations were studied
in \cite{L2}.
 For manifolds with semisimple quantum cohomology, these equations determine
the genus-1 and genus-2 Gromov-Witten invariants in terms of genus-0 invariants
(cf. \cite{DZ} for the genus-1 case and \cite{L1} for the genus-2 case). In \cite{KL},
the authors proved a genus-3 topological recursion relation by studying a tautological
relation on $\overline{\cal M}_{3, 1}$. Certain universal equations of all genera
were conjectured in \cite{LX} and proved in \cite{LP}.
Despite all of this progress,
the understanding of universal equations is
still very limited and unsatisfactory. For example, so far the known genus-3 equations still can
not determine the genus-3 generating function even for manifolds with semisimple
quantum cohomology. Therefore it is very interesting to find more genus-3 universal equations.
The main purpose of this paper is to obtain new genus-2 and genus-3
universal equations by studying tautological relations
on $\overline{\cal M}_{3, 2}$.

Universal equations play very important roles in studying Gromov-Witten invariants. Although lots of progress
has been made for Gromov-Witten theory of manifolds with semisimple quantum cohomology (see, for example, Teleman's work \cite{T}),
very little is known for non-semisimple case. However most manifolds do not have semisimple
quantum cohomology. Therefore it is very important to study universal equations
which could be applied to manifolds whose quantum cohomology may not be semisimple.
For each fixed genus $g$, one may expect to obtain many genus-$g$ universal equations by studying
moduli spaces $\overline{\cal M}_{g, k}$ for all possible $k$. However, most of such relations
can be derived from just a few relations which we may call the {\it generators} of genus-$g$ relations.
For genus-0 and genus-1 case, this is more or less understood from the work of \cite{Ke} and \cite{Pe}.
Starting from genus-2, much less is known. It is extremely interesting to find out generators for universal
equations in each genus. Other universal equations can then be easily obtained from the generators by taking derivatives
and other algebraic manipulations. Part of the motivations for this paper is to find universal equations which, together with other equations, could serve as generators for genus-3 and genus-2 universal equations. We hope that these low genus examples
will be useful for finding general patterns for generators of higher genus universal equations.

To describe the universal equations, we need to use some
operators introduced in \cite{L1}.
Let $M$ be a compact symplectic manifold. The {\it big phase space}
for Gromov-Witten invariants of $M$ is a product of infinitely many copies of
$H^{*}(M; {\mathbb C})$. We will choose a basis $\{ \ga \mid \alpha = 1, \ldots, N \}$
of $H^{*}(M; {\mathbb C})$. The quantum product $\vw_{1} \qp \vw_{2}$
of two vector fields $\vw_{1}$ and $\vw_{2}$ on the big phase space was introduced
in \cite{L1}. This is an associative product without an identity element. An
operator $T$ on the space of vector fields on the big phase space was also
introduced in \cite{L1} to measure the failure of the string vector field
to be an identity element with respect to this product.
This operator turns out to be a very useful device to translate relations in the
tautological rings of $\overline{\cal M}_{g,n}$ into universal equations for
Gromov-Witten invariants. We will write universal equations of Gromov-Witten
invariants as equations among tensors $\gwiig{\vw_{1} \, \cdots \, \vw_{k}}$
which are defined to be the $k$-th covariant derivatives of the generating
functions of genus-$g$ Gromov-Witten invariants with respect to the trivial
connection on the big phase space. We will briefly review these definitions
in Section \ref{sec:Pre} for completeness.

For vector fields $\vw_1$ and $\vw_2$ on the big phase space, we define
\allowdisplaybreaks
\begin{eqnarray}
&&   Q(\vw_{1}, \vw_{2})
    \nonumber \\
&:=& \frac{2}{9} \, \gwiitwo{\vw_{2} \, T(\gua)} \gwii{\ga \, \vw_{1} \, \gub \, \gb}
    + \frac{5}{24} \, \gwiitwo{T(\vw_{2}) \, \gua} \gwii{\ga \, \vw_{1} \, \gub \, \gb}
    \nonumber \\
&&
    + \frac{16}{3} \, \gwiitwo{\vw_{2} \, T(\gua)} \gwiione{\left\{ \ga \qp \vw_{1} \right\}}
    + 5  \, \gwiitwo{T(\vw_{2}) \, \gua} \gwiione{\left\{ \ga \qp \vw_{1} \right\}}
     \nonumber \\
&&
    + \frac{40}{3} \, \gwiitwo{T(\gua)} \gwiione{\left\{ \ga \qp \vw_{1} \right\} \, \vw_{2} }
   + \frac{1}{6} \, \gwiitwo{\vw_{2} \, T(\gua) \, \left\{\ga \qp \vw_{1} \right\}}
     \nonumber \\
&&
    + \frac{1}{2} \, \gwiitwo{T(\vw_{2}) \, \gua \, \left\{\ga \qp \vw_{1} \right\}}
    + \frac{2}{9} \, \gwiitwo{\vw_{1} \, \left\{\vw_{2} \qp \gua \qp \ga \right\}}
    \nonumber \\
&&
    + \frac{1}{18} \, \gwiione{\vw_{1} \, \gua } \gwiione{\ga \, \vw_{2} \, \left\{ \gub \qp \gb \right\}}
    + \frac{1}{18} \, \gwiione{\vw_{1} \, \gua } \gwiione{ \ga \, \gub } \gwii{\gb \, \vw_{2} \, \gum \, \gm}
    \nonumber \\
&&
    + \frac{1}{30} \, \gwiione{ \gua } \gwiione{ \ga \, \vw_{2} \, \gub} \gwii{ \gb \, \vw_{1} \, \gum \, \gm}
    + \frac{1}{30} \, \gwiione{  \gua} \gwiione{\ga \, \vw_{2} \, \gub \, \left\{ \gb \qp \vw_{1} \right\}}
    \nonumber \\
&&
    + \frac{9}{10} \, \gwiione{ \vw_{2} \, \gua \, \gub} \gwiione{ \ga \, \left\{ \gb \qp \vw_{1} \right\}}
    + \frac{1}{30} \, \gwiione{\gua \, \ga \, \vw_{2} \, \gub} \gwiione{\left\{ \gb \qp \vw_{1} \right\}}
            \nonumber \\
&&
    + \frac{2}{15} \, \gwiione{\gua \, \ga \,  \gub} \gwiione{\left\{ \gb \qp \vw_{2} \right\} \, \vw_{1} }
    + \frac{3}{5} \, \gwiione{\vw_{2} \, \gua \, \gub} \gwii{\ga \, \gb \, \vw_{1} \, \gum} \gwiione{\gm}
    \nonumber \\
&&
    + \frac{1}{6} \, \gwiione{\vw_{1} \, \gua \, \gub} \gwiione{ \left\{ \ga \qp \gb \right\} \, \vw_{2}}
    + \frac{16}{15} \, \gwiione{ \gua \, \gub} \gwii{\ga \, \gb \, \vw_{1} \, \gum} \gwiione{ \gm \, \vw_{2} }
    \nonumber \\
&&
    + \frac{41}{90} \, \gwiione{\vw_{2} \, \gua} \gwii{\ga \, \vw_{1} \, \gub \, \gb \, \gum} \gwiione{\gm}
    + \frac{4}{5} \, \gwiione{  \gua} \gwiione{\ga \, \vw_{2} \, \gub}
                \gwiione{ \left\{ \gb \qp \vw_{1} \right\}}
            \nonumber \\
&&
    + \frac{16}{5} \, \gwiione{  \gua} \gwiione{\ga  \, \gub}
                \gwiione{ \left\{ \gb \qp \vw_{2} \right\} \, \vw_{1}}
    + \frac{92}{15} \, \gwiione{\vw_{2} \, \gua} \gwiione{ \gub} \gwiione{ \gum}
                \gwii{\ga \, \gb \, \gm \, \vw_{1}}
                \nonumber \\
&&
    + \frac{1}{720} \, \gwiione{\vw_{2} \, \gua \, \ga \, \gub} \gwii{\gb \, \vw_{1} \, \gum \, \gm}
    + \frac{1}{720} \, \gwiione{\gua \, \ga \, \vw_{2} \, \gub \, \left\{ \gb \qp \vw_{1} \right\}}
    \nonumber \\
&&
    + \frac{1}{135} \, \gwiione{\vw_{2} \, \gua \, \gub \, \gum} \gwii{\ga \, \gb \, \gm \, \vw_{1}}
    +\frac{1}{40} \, \gwiione{\vw_{2}  \, \gua \, \gub} \gwii{\ga \, \gb \, \vw_{1} \, \gum \, \gm}
    \nonumber \\
&&
    + \frac{1}{120}  \, \gwiione{  \vw_{2} \, \gua} \gwii{\ga \, \vw_{1} \, \gub \, \gb \, \gum \, \gm}.
    \label{eqn:Qdef}
\end{eqnarray}
Note that $Q(\vw_1, \vw_2)$ only involves data with genus $\leq$ 2.
In this paper we will prove the following genus-3 universal equation
\begin{thm} \label{thm:g3T2Tskew}
 For all compact symplectic manifolds,
\[
\gwiih{3}{T^{2}(\vw_{1}) T(\vw_{2})} - \gwiih{3}{T^{2}(\vw_{2}) T(\vw_{1})}
= \frac{1}{7} \left\{ Q(\vw_1, \vw_2) - Q(\vw_2, \vw_1) \right\}
\]
for all vector fields $\vw_1$ and $\vw_2$ on the big phase space.
\end{thm}
Note that in general
\[ \gwiih{3}{T^{2}(\vw_{1}) T(\vw_{2})} \, \neq \, \frac{1}{7} \, Q(\vw_1, \vw_2).\]
Instead $\gwiih{3}{T^{2}(\vw_{1}) T(\vw_{2})}$ is given by a much more complicated
genus-$3$ universal equation in Theorem~\ref{thm:g3T2Ta2},
which is one of the main results of this paper. Theorem~\ref{thm:g3T2Tskew}
is actually a corollary of Theorem~\ref{thm:g3T2Ta2}.
We would like to point out that the genus-3
topological recursion relation in \cite{KL} also follows from the formula in
Theorem~\ref{thm:g3T2Ta2} (see the last paragraph of Section~\ref{sec:compare}).
Therefore the equation in Theorem~\ref{thm:g3T2Ta2} is a better candidate for
the generators of genus-3 universal equations.

During the proof of Theorem~\ref{thm:g3T2Ta2},
we also discovered a new
 genus-2 universal equation
which will be given in Theorem~\ref{pro:g2uni}.
This equation does not
follow from  known genus-2 equations. Comparison of universal equations given in this paper and
previously known universal equations will be given in Section \ref{sec:compare}.
By a result of \cite{L1}, known genus-2 equations already completely determine
the genus-2 generating function for manifolds with semisimple quantum cohomology.
Because of this fact, the new genus-2 equation in Theorem~\ref{pro:g2uni}
is quite a surprise. Theorem~\ref{pro:g2uni} is another main result of this paper.

The second author would like to thank
R. Pandharipande for helpful discussions on whether the tautological ring
$R^*(\overline{\cal M}_{3, 2})$ is Gorenstein. The authors would also like to thank an anonymous referee for
pointing out the fact that the tautological cohomology of $\overline{\cal M}_{3, 2}$ is Gorenstein,
which follows from
results in \cite{B} and \cite{Y}.

\begin{rem}
This paper is a longer version of a paper which will appear in SCIENCE CHINA Mathematics.
\end{rem}

\section{Preliminaries}
\label{sec:Pre}

Let $M$ be a compact symplectic manifold.
The {\it big phase space} is by definition the infinite product
\[  P := \prod_{n=0}^{\infty} H^{*}(M; \mathbb{C}). \]
Fix a basis $\{ \gamma_{0}, \ldots, \gamma_{N} \}$ of
$H^{*}(M; \mathbb{C})$, where $\gamma_{0}$ is the identity element, of the ordinary
cohomology ring of $M$. Then we denote the corresponding basis for
the $n$-th copy of $H^{*}(M; \mathbb{C})$ in $P$ by
$\{\tau_{n}(\gamma_{0}), \ldots, \tau_{n}(\gamma_{N}) \}$.
We call $\grava{n}$ a {\it descendant of $\gamma_{\alpha}$ with descendant
level $n$}.
We can think of $P$ as an infinite dimensional vector space with a basis
$\{ \grava{n} \mid 0 \leq \alpha \leq N, \, \, \, n \in \mathbb{Z}_{\geq 0} \}$
where $\mathbb{Z}_{\geq 0} = \{ n \in \mathbb{Z} \mid n \geq 0\}$.
Let
$(t_{n}^{\alpha} \mid 0 \leq \alpha \leq N, \, \, \, n \in \mathbb{Z}_{\geq 0})$
be the corresponding coordinate system on $P$.
For convenience, we identify $\grava{n}$ with the coordinate vector field
$\frac{\partial}{\partial t_{n}^{\alpha}}$ on $P$ for $n \geq 0$.
If $n<0$, $\grava{n}$ is understood to be the $0$ vector field.
We also abbreviate $\grava{0}$ by $\gamma_{\alpha}$.
We use $\tau_{+}$ and $\tau_{-}$ to denote the operators which shift the level
of descendants by $1$, i.e.
\[ \tau_{\pm} \left(\sum_{n, \alpha} f_{n, \alpha} \grava{n}\right)
    = \sum_{n, \alpha} f_{n, \alpha} \grava{n \pm 1} \]
where $f_{n, \alpha}$ are functions on the big phase space.

We will adopt the following {\it notational conventions}:
Lower case Greek letters, e.g. $\alpha$, $\beta$, $\mu$, $\nu$,
$\sigma$,..., etc., will be used to index the cohomology classes on $M$.
These indices run from $0$ to $N$. Lower case English
letters, e.g. $i$, $j$, $k$, $m$, $n$, ..., etc., will be used to
index the level of descendants. These indices run over the set of all
non-negative integers, i.e. $\mathbb{Z}_{\geq 0}$. All summations are
over the entire ranges of the corresponding indices unless otherwise
indicated.
Let
\[ \eta_{\alpha \beta} = \int_{M} \gamma_{\alpha} \cup
    \gamma_{\beta}
\]
 be the intersection form on $H^{*}(M, \mathbb{C})$.
We will use $\eta = (\eta_{\alpha \beta})$ and $\eta^{-1} =
(\eta^{\alpha \beta})$ to lower and raise indices.
For example,
\[ \gua :=  \eta^{\alpha \beta} \gb.\]
Here we are using the summation convention that repeated
indices (in this formula, $\beta$) should be summed
over their entire ranges.

Let
\[ \gwih{g, d}{\grav{n_{1}}{\alpha_{1}} \, \grav{n_{2}}{\alpha_{2}} \,
    \ldots \, \grav{n_{k}}{\alpha_{k}}} :=
    \int_{[\overline{\cal M}_{g,n}(M;d)]^\mathrm{virt}} \bigcup_{i=1}^k
    (\Psi^{n_{i}} \cup \mathrm{ev_{i}}^*\gamma_{\alpha_{i}})
\]
be the genus-$g$, degree $d$,
 descendant Gromov-Witten invariant associated
to $\gamma_{\alpha_{1}}, \ldots, \gamma_{\alpha_{k}}$ and nonnegative
integers $n_{1}, \ldots, n_{k}$
(cf. \cite{W}, \cite{RT}, \cite{LiT}).
Here, $\overline{\cal M}_{g, k}(M; d)$ is the moduli space of stable maps
from genus-$g$, $k$-pointed curves to $M$ of degree $d \in H_{2}(M; \mathbb{Z})$.
$\Psi_{i}$ is the first Chern class of the tautological line bundle
over $\overline{\cal M}_{g, k}(M; d)$ whose geometric fiber over a stable map
is the cotangent space of the domain curve at the $i$-th marked point while
$\mathrm{ev_i}:\overline{\cal M}_{g,n}(M; d)\to M$ is the $i$-th evaluation map for all
 $i=1,\ldots,k$. Finally, $[\overline{\cal M}_{g,n}(M;d)]^\mathrm{virt}$ is the virtual
 fundamental class.
The genus-$g$
generating function is defined to be
\[ F_{g} =  \sum_{k \geq 0} \frac{1}{k!}
         \sum_{ \begin{array}{c}
        {\scriptstyle \alpha_{1}, \ldots, \alpha_{k}} \\
                {\scriptstyle  n_{1}, \ldots, n_{k}}
                \end{array}}
                t^{\alpha_{1}}_{n_{1}} \cdots t^{\alpha_{k}}_{n_{k}}
    \sum_{d \in H_{2}(V, \mathbb{Z})} q^{d}
    \gwih{g, d}{\grav{n_{1}}{\alpha_{1}} \, \grav{n_{2}}{\alpha_{2}} \,
        \ldots \, \grav{n_{k}}{\alpha_{k}}} \]
where $q^{d}$ belongs to the Novikov ring.
This function is understood as a formal power series n the variables
    $\{\,t_{n}^{\alpha}\,\}$ with coefficients in the Novikov ring.

Introduce a $k$-tensor
 $\left< \left< \right. \right. \underbrace{\cdot \cdots \cdot}_{k}
        \left. \left. \right> \right> $
defined by
\[ \gwiig{{\cal W}_{1} {\cal W}_{2} \cdots {\cal W}_{k}} \, \,
         := \sum_{m_{1}, \alpha_{1}, \ldots, m_{k}, \alpha_{k}}
                f^{1}_{m_{1}, \alpha_{1}} \cdots f^{k}_{m_{k}, \alpha_{k}}
        \, \, \, \frac{\partial^{k}}{\partial t^{\alpha_{1}}_{m_{1}}
            \partial t^{\alpha_{2}}_{m_{k}} \cdots
            \partial t^{\alpha_{k}}_{m_{k}}} F_{g},
 \]
for vector fields ${\cal W}_{i} = \sum_{m, \alpha}
        f^{i}_{m, \alpha} \, \frac{\partial}{\partial t_{m}^{\alpha}}$ where
$f^{i}_{m, \alpha}$ are functions on the big phase space.
This tensor is called the {\it $k$-point (correlation) function}.

For any vector fields $\vw_{1}$ and $\vw_{2}$ on the big phase space,
the quantum product of $\vw_{1}$ and $\vw_{2}$ is defined by
\[ \vw_{1} \qp \vw_{2} := \gwii{\vw_{1} \, \vw_{2} \, \gua} \ga. \]
Define the vector field
\[ T(\vw) := \tau_{+}(\vw) - \gwii{\vw \, \gua} \ga \]
for any vector field $\vw$. The operator $T$
was introduced in \cite{L1} as a convenient tool in the study of universal
equations for Gromov-Witten invariants. Let $\psi_{i}$ be the first
Chern class of the tautological line bundle over $\overline{\cal M}_{g, k}$
whose geometric fiber over a stable curve is the cotangent space of the curve
at the $i$-th marked point.  When we translate a relation in the tautological
ring of $\overline{\cal M}_{g, k}$ into differential equations for
generating functions of Gromov-Witten invariants,
each $\psi$ class corresponds to the insertion of the operator $T$. Let
$\nabla$ be the trivial flat connection on the big phase space with respect
to which $\grava{n}$ are parallel vector fields for all $\alpha$ and
$n$. Then the covariant derivative of the quantum product satisfies
\[ 
\nabla_{\vw_{3}} (\vw_{1} \qp \vw_{2})
    = (\nabla_{\vw_{3}} \vw_{1}) \qp \vw_{2}
        +  \vw_{1} \qp (\nabla_{\vw_{3}} \vw_{2})
        + \gwii{\vw_{1} \, \vw_{2} \, \vw_{3} \, \gua} \ga
\]
and the covariant derivative of the operator $T$ is given by
\[ 
 \nabla_{\vw_{2}} \, \, T(\vw_{1}) = T(\nabla_{\vw_{2}} \vw_{1})
    - \vw_{2} \qp \vw_{1}
\]
for any vector fields $\vw_{1}, \vw_{2}$ and  $\vw_{3}$
(cf. \cite[Equation (8) and Lemma 1.5]{L1}).
We need to use these formulas in order to compute derivatives of universal equations.

\section{Topological Recursion Relations on $\overline{\cal M}_{3, 2}$}
\label{sec:TRR}

In this section, we will give a simultaneous proof of  the main theorems in this paper,
i.e. Theorem \ref{pro:g2uni} and  Theorem \ref{thm:g3T2Ta2} below. Theorem~\ref{thm:g3T2Tskew}
will be proved using Theorem \ref{thm:g3T2Ta2}.

It is well known that the cohomology class $\psi_{1}^{2}\, \psi_2$ vanishes on ${\cal M}_{3,2}$ due
to a result of Ionel (cf. \cite{Io}).
 Furthermore, by a result of Faber and Pandharipande
\cite{FP2}, $\psi_{1}^{2}\, \psi_2$ is equal to a  class from the boundary
strata which is tautological, and therefore is a linear combination
of products of $\psi$ and $\kappa$ classes and fundamental classes
of some boundary strata. For any curve in
$\partial {\cal M}_{3,2}:= \overline{\cal M}_{3,2} - {\cal M}_{3,2}$
with a genus-3 component, it has a dual graph
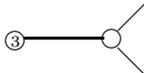
\[ \begin{picture}(120, 30)
\put(5.5, 5){$\scriptscriptstyle 3$}
\put(7, 7){\circle{7}}
\put(10.5, 8){\line(1, 0){30}}
\put(40, 8){ \circle{7}}
\put(46, 11){\line(1, 1){10}}
\put(46, 4.5){\line(1, -1){10}}
\end{picture}
\]
where a hollow circle
$\begin{picture}(10,8)(0,0)\put(5,3){\circle{7}}\end{picture}$
represents a vertex of genus $0$, and
 $\begin{picture}(10,8)(0,0)\put(5,3){\circle{7}} \put(3,2){$\scriptscriptstyle g$}\end{picture}$
represents a vertex of genus $g \geq 1$.
Since this graph is a tree, all such curves are of compact type.
For this stratum to occur in the expression of $\psi_{1}^{2}\, \psi_2$,
it must be multiplied by combinations of $\psi$-classes and $\kappa$-classes
of degree 2 on the  genus-3 component. By a result of Yang \cite{Y},
 $\kappa$-classes in this expression
can be replaced with combinations of $\psi$-classes and fundamental classes
of boundary strata.
On all other boundary strata, all components of curves must have genus
at most $2$, therefore $\kappa$ classes again can be represented as
linear combinations of $\psi$ classes and fundamental classes of
boundary strata (cf. \cite{AC}). Therefore, it follows that
$\psi_{1}^{2}\, \psi_2$ on $\overline{\cal M}_{3,2}$ can be written as a
linear combination of products of the $\psi$ classes and the
fundamental classes of some boundary strata. By taking into
consideration the genus-0 and genus-1 topological recursion
relations as well as Mumford's genus-2  relation, we can translate
these results into the following universal equations for Gromov-Witten
invariants with unknown
constants $a_{1}, \ldots, a_{105}$:

\allowdisplaybreaks
\begin{eqnarray}\label{eq:PhiDef}
0&=&  \Phi(\vw_1, \vw_2) \nonumber \\
&:=& -\gwiih{3}{T^{2}(\vw_{1}) T(\vw_{2})} \nonumber \\
&& + a_{1} \gwiih{3}{T^{2}(\vw_{1} \qp \vw_{2})}
    + a_{2} \gwiitwo{ \vw_{1} \, \vw_{2} \,  T(\ga \qp \gua)}
    \nonumber \\
&&
    + a_{3} \gwiitwo{T(\vw_{1}) \, \gua} \gwii{\ga \, \vw_{2} \, \gub \, \gb}
    +a_{4} \gwiitwo{\vw_{1} \, T(\gua)} \gwii{\ga \, \vw_{2} \, \gub \, \gb}
    \nonumber \\
&&
    + a_{5} \gwiitwo{T(\vw_{2}) \, \gua} \gwii{\ga \, \vw_{1} \, \gub \, \gb}
    +a_{6} \gwiitwo{\vw_{2} \, T(\gua)} \gwii{\ga \, \vw_{1} \, \gub \, \gb}
    \nonumber \\
&&
    +a_{7} \gwiitwo{ T(\gua)} \gwii{\ga \, \vw_{1} \, \vw_{2} \, \gub \, \gb}
    + a_{8} \gwiitwo{\vw_{1} \, T(\gua)} \gwiione{\left\{ \ga \qp \vw_{2} \right\}} \nonumber \\
&&
    + a_{9} \gwiitwo{T(\vw_{2}) \, \gua} \gwiione{\left\{ \ga \qp \vw_{1} \right\}}
    + a_{10} \gwiitwo{\vw_{2} \, T(\gua)} \gwiione{\left\{ \ga \qp \vw_{1} \right\}} \nonumber \\
&&
    + a_{11} \gwiitwo{T(\gua)} \gwiione{\left\{ \ga \qp \vw_{1} \right\} \, \vw_{2} }
    + a_{12} \gwiitwo{T(\gua)} \gwiione{\left\{ \ga \qp \vw_{2} \right\} \, \vw_{1} } \nonumber \\
&&
    + a_{13} \gwiitwo{T(\gua)} \gwii{\ga \, \vw_{1} \, \vw_{2} \, \gub } \gwiione{\gb}
    + a_{14} \gwiitwo{\vw_{1} \, T(\gua) \, \left\{\ga \qp \vw_{2} \right\}}
     \nonumber \\
&&
    + a_{15} \gwiitwo{T(\vw_{2}) \, \gua \, \left\{\ga \qp \vw_{1} \right\}}
    + a_{16} \gwiitwo{\vw_{2} \, T(\gua) \, \left\{\ga \qp \vw_{1} \right\}}
    \nonumber \\
&&
    + a_{17} \gwiitwo{T(\gua) \, \gub} \gwii{\ga \, \gb \, \vw_{1} \, \vw_{2}}
    + a_{18} \gwiitwo{T(\gua) \, \ga \, \left\{\vw_{1} \qp \vw_{2} \right\}}
    \nonumber \\
&&
    + a_{19} \gwiitwo{\gua \, \ga \, T(\vw_{1} \qp \vw_{2})}
    + a_{20} \gwiitwo{T(\gua)} \gwiione{\ga \, \left\{ \vw_{1} \qp \vw_{2} \right\}}
    \nonumber \\
&&
    + a_{21} \gwiitwo{T(\vw_{1} \qp \vw_{2}) \, \gua } \gwiione{\ga}
    + a_{22} \gwiitwo{\left\{ \vw_{1} \qp \vw_{2} \right\} \, T(\gua) } \gwiione{\ga}
    \nonumber \\
&&
    + a_{23} \gwiitwo{\vw_{1} \, \left\{\vw_{2} \qp \gua \qp \ga \right\}}
    + a_{24} \gwiitwo{\vw_{2} \, \left\{\vw_{1} \qp \gua \qp \ga \right\}}
    \nonumber \\
&&
    + a_{25} \gwiitwo{\left\{\vw_{1} \qp \vw_{2} \right\} \, \left\{ \gua \qp \ga \right\}}
    + a_{26} \gwiitwo{\gua} \gwii{\ga \, \left\{ \vw_{1} \qp \vw_{2} \right\} \, \gub \, \gb}
    \nonumber \\
&&
    + a_{27} \gwiitwo{\gua} \gwii{\ga \, \vw_{1} \, \gub \, \left\{ \gb \qp \vw_{2} \right\}}
    + a_{28} \gwiitwo{\gua} \gwii{\left\{ \ga \qp \gub \right\} \, \gb \, \vw_{1} \, \vw_{2}}
    \nonumber \\
&&
    + a_{29} \gwiitwo{ \gua \, \left\{\ga \qp \vw_{1} \, \qp \vw_{2} \right\}}
    + a_{30} \gwiitwo{ \gua} \gwii{\ga \, \vw_{1} \, \vw_{2} \, \left\{ \gub \qp \gb \right\}}
    \nonumber \\
&&
    + a_{31} \gwiitwo{ \gua} \gwiione{\left\{\ga \qp \vw_{1} \qp \vw_{2} \right\}}
    + a_{32} \gwiione{\vw_{1} \, \vw_{2} \, \gua} \gwiione{\ga \, \left\{ \gub \qp \gb \right\}}
    \nonumber \\
&&
    + a_{33} \gwiione{\vw_{1} \, \gua } \gwiione{\ga \, \vw_{2} \, \left\{ \gub \qp \gb \right\}}
    + a_{34} \gwiione{\vw_{2} \, \gua } \gwiione{\ga \, \vw_{1} \, \left\{ \gub \qp \gb \right\}}
    \nonumber \\
&&
    + a_{35} \gwiione{\vw_{1} \, \gua } \gwiione{ \ga \, \gub } \gwii{\gb \, \vw_{2} \, \gum \, \gm}
    + a_{36} \gwiione{\vw_{2} \, \gua } \gwiione{ \ga \, \gub } \gwii{\gb \, \vw_{1} \, \gum \, \gm}
    \nonumber \\
&&
    + a_{37} \gwiione{ \gua} \gwiione{\ga \, \vw_{1} \, \vw_{2} \, \left\{\gub \qp \gb \right\} }
    + a_{38} \gwiione{ \gua } \gwiione{ \ga \, \vw_{1} \, \gub} \gwii{ \gb \, \vw_{2} \, \gum \, \gm}
    \nonumber \\
&&
    + a_{39} \gwiione{ \gua } \gwiione{ \ga \, \vw_{2} \, \gub} \gwii{ \gb \, \vw_{1} \, \gum \, \gm}
    + a_{40} \gwiione{ \gua } \gwiione{ \ga \,  \gub}
                \gwii{ \gb \, \vw_{1} \, \vw_{2} \, \gum \, \gm}
    \nonumber \\
&&
    + a_{41} \gwiione{ \vw_{1} \, \gua} \gwiione{\ga \, \gub \, \left\{ \gb \qp \vw_{2} \right\}}
    + a_{42} \gwiione{ \vw_{2} \, \gua} \gwiione{\ga \, \gub \, \left\{ \gb \qp \vw_{1} \right\}}
    \nonumber \\
&&
    + a_{43} \gwiione{  \gua} \gwiione{\ga \, \vw_{1} \, \gub \, \left\{ \gb \qp \vw_{2} \right\}}
    + a_{44} \gwiione{  \gua} \gwiione{\ga \, \vw_{2} \, \gub \, \left\{ \gb \qp \vw_{1} \right\}}
    \nonumber \\
&&
    + a_{45} \gwiione{ \gua } \gwiione{ \ga \, \gub \, \gum}
            \gwii{\gb \, \gm \, \vw_{1} \, \vw_{2}}
    + a_{46} \gwiione{ \vw_{1} \, \gua \, \gub} \gwiione{ \ga \, \left\{ \gb \qp \vw_{2} \right\}}
    \nonumber \\
&&
    + a_{47} \gwiione{ \vw_{2} \, \gua \, \gub} \gwiione{ \ga \, \left\{ \gb \qp \vw_{1} \right\}}
    + a_{48} \gwiione{ \gua \, \gub} \gwiione{\ga \, \gum} \gwii{ \gm \, \gb \, \vw_{1} \, \vw_{2}}
    \nonumber \\
&&
    + a_{49} \gwiione{ \gua \, \gub} \gwiione{\ga \, \gb \, \left\{ \vw_{1} \qp \vw_{2} \right\} }
    + a_{50} \gwiione{ \gua } \gwiione{\ga \, \gub \, \gb \, \left\{ \vw_{1} \qp \vw_{2} \right\}}
    \nonumber \\
&&
    + a_{51} \gwiione{ \gua \, \ga \, \gub} \gwiione{ \gb \, \left\{ \vw_{1} \qp \vw_{2} \right\} }
    + a_{52} \gwiione{\gua \, \ga \, \vw_{1} \, \gub} \gwiione{\left\{ \gb \qp \vw_{2} \right\}}
    \nonumber \\
&&
    + a_{53} \gwiione{\gua \, \ga \, \vw_{2} \, \gub} \gwiione{\left\{ \gb \qp \vw_{1} \right\}}
    + a_{54} \gwiione{\gua \, \ga \, \gub} \gwii{ \gb \, \vw_{1} \, \vw_{2} \, \gum}
            \gwiione{ \gm}
            \nonumber \\
&&
    + a_{55} \gwiione{\gua \, \ga \,  \gub} \gwiione{\left\{ \gb \qp \vw_{1} \right\} \, \vw_{2} }
    + a_{56} \gwiione{\gua \, \ga \,  \gub} \gwiione{\left\{ \gb \qp \vw_{2} \right\} \, \vw_{1} }
    \nonumber \\
&&
    + a_{57} \gwiione{\vw_{1} \, \vw_{2} \, \gua \, \gub} \gwiione{ \left\{ \ga \qp \gb \right\}}
    + a_{58} \gwiione{\vw_{1} \, \gua \, \gub} \gwii{\ga \, \gb \, \vw_{2} \, \gum} \gwiione{\gm}
    \nonumber \\
&&
    + a_{59} \gwiione{\vw_{2} \, \gua \, \gub} \gwii{\ga \, \gb \, \vw_{1} \, \gum} \gwiione{\gm}
    + a_{60} \gwiione{ \gua \, \gub} \gwii{\ga \, \gb \, \vw_{1} \, \vw_{2} \, \gum} \gwiione{\gm}
    \nonumber \\
&&
    + a_{61} \gwiione{\vw_{1} \, \gua \, \gub} \gwiione{ \left\{ \ga \qp \gb \right\} \, \vw_{2}}
    + a_{62} \gwiione{\vw_{2} \, \gua \, \gub} \gwiione{ \left\{ \ga \qp \gb \right\} \, \vw_{1}}
    \nonumber \\
&&
    + a_{63} \gwiione{ \gua \, \gub} \gwii{\ga \, \gb \, \vw_{1} \, \gum} \gwiione{ \gm \, \vw_{2} }
    + a_{64} \gwiione{ \gua \, \gub} \gwii{\ga \, \gb \, \vw_{2} \, \gum} \gwiione{ \gm \, \vw_{1} }
    \nonumber \\
&&
    + a_{65} \gwiione{ \gua \, \gub} \gwiione{ \left\{ \ga \qp \gb \right\} \, \vw_{1} \, \vw_{2}}
    + a_{66} \gwiione{\vw_{1} \, \vw_{2} \, \gua} \gwii{\ga \, \gub \, \gb \, \gum} \gwiione{ \gm}
    \nonumber \\
&&
    + a_{67} \gwiione{\vw_{1} \, \gua} \gwii{\ga \, \vw_{2} \, \gub \, \gb \, \gum} \gwiione{\gm}
    + a_{68} \gwiione{\vw_{2} \, \gua} \gwii{\ga \, \vw_{1} \, \gub \, \gb \, \gum} \gwiione{\gm}
    \nonumber \\
&&
    + a_{69} \gwiione{\vw_{1} \, \gua} \gwii{\ga \,  \gub \, \gb \, \gum} \gwiione{\gm \, \vw_{2}}
    + a_{70} \gwiione{ \gua} \gwii{\ga \, \vw_{1} \, \vw_{2} \, \gub \, \gb \, \gum} \gwiione{\gm}
    \nonumber \\
&&
    + a_{71} \gwiione{\gua} \gwiione{ \ga \, \gub} \gwiione{\gb \, \left\{\vw_{1} \qp \vw_{2} \right\}}
    + a_{72} \gwiione{ \vw_{1} \, \gua} \gwiione{\ga \, \gub}
                \gwiione{ \left\{ \gb \qp \vw_{2} \right\}}
    \nonumber \\
&&
    + a_{73} \gwiione{ \vw_{2} \, \gua} \gwiione{\ga \, \gub}
                \gwiione{ \left\{ \gb \qp \vw_{1} \right\}}
    + a_{74} \gwiione{  \gua} \gwiione{\ga \, \vw_{1} \, \gub}
                \gwiione{ \left\{ \gb \qp \vw_{2} \right\}}
    \nonumber \\
&&
    + a_{75} \gwiione{  \gua} \gwiione{\ga \, \vw_{2} \, \gub}
                \gwiione{ \left\{ \gb \qp \vw_{1} \right\}}
    + a_{76} \gwiione{\gua} \gwiione{\ga \, \gub} \gwii{\gb \, \vw_{1} \, \vw_{2} \, \gum}
            \gwiione{\gm}
            \nonumber \\
&&
    + a_{77} \gwiione{  \gua} \gwiione{\ga  \, \gub}
                \gwiione{ \left\{ \gb \qp \vw_{1} \right\} \, \vw_{2}}
    + a_{78} \gwiione{  \gua} \gwiione{\ga  \, \gub}
                \gwiione{ \left\{ \gb \qp \vw_{2} \right\} \, \vw_{1}}
                \nonumber \\
&&
    + a_{79} \gwiione{\vw_{1} \, \vw_{2} \, \gua} \gwiione{\gub} \gwiione{\left\{\ga \qp \gb \right\}}
    + a_{80} \gwiione{\vw_{1} \, \gua} \gwiione{ \gub} \gwiione{ \gum}
                \gwii{\ga \, \gb \, \gm \, \vw_{2}}
                \nonumber \\
&&
    + a_{81} \gwiione{\vw_{2} \, \gua} \gwiione{ \gub} \gwiione{ \gum}
                \gwii{\ga \, \gb \, \gm \, \vw_{1}}
                \nonumber \\
&&
    + a_{82} \gwiione{ \gua} \gwiione{ \gub} \gwiione{ \gum}
                \gwii{\ga \, \gb \, \gm \, \vw_{1} \, \vw_{2}}
                \nonumber \\
&&
    + a_{83} \gwiione{\vw_{1} \,  \gua} \gwiione{\vw_{2} \, \gub} \gwiione{\left\{\ga \qp \gb \right\}}
    + a_{84} \gwiione{\vw_{1} \, \vw_{2} \, \gua \, \ga \, \left\{\gub \, \qp \gb \right\}}
    \nonumber \\
&&
    + a_{85} \gwiione{\vw_{1} \, \gua \, \ga \, \gub} \gwii{\gb \, \vw_{2} \, \gum \, \gm}
    + a_{86} \gwiione{\vw_{2} \, \gua \, \ga \, \gub} \gwii{\gb \, \vw_{1} \, \gum \, \gm}
    \nonumber \\
&&
    + a_{87} \gwiione{\gua \, \ga \, \gub} \gwii{\gb \, \vw_{1} \, \vw_{2} \, \gum \, \gm}
    + a_{88} \gwiione{\gua \, \ga \, \gub \, \gb \, \left\{ \vw_{1} \qp \vw_{2} \right\}}
    \nonumber \\
&&
    + a_{89} \gwiione{\gua \, \ga \, \vw_{1} \, \gub \, \left\{ \gb \qp \vw_{2} \right\}}
    + a_{90} \gwiione{\gua \, \ga \, \vw_{2} \, \gub \, \left\{ \gb \qp \vw_{1} \right\}}
    \nonumber \\
&&
    + a_{91} \gwiione{\gua \, \ga \, \gub \, \gum} \gwii{\gb \, \gm \, \vw_{1} \, \vw_{2}}
    + a_{92} \gwiione{\vw_{1} \, \vw_{2} \, \gua \, \gub \, \left\{\ga \qp \gb\right\}}
    \nonumber \\
&&
    + a_{93} \gwiione{\vw_{1} \, \gua \, \gub \, \gum} \gwii{\ga \, \gb \, \gm \, \vw_{2}}
    + a_{94} \gwiione{\vw_{2} \, \gua \, \gub \, \gum} \gwii{\ga \, \gb \, \gm \, \vw_{1}}
    \nonumber \\
&&
    + a_{95} \gwiione{ \gua \, \gub \, \gum} \gwii{\ga \, \gb \, \gm \, \vw_{1} \, \vw_{2}}
    + a_{96} \gwiione{\vw_{1} \, \vw_{2} \, \gua \, \gub} \gwii{\ga \, \gb \, \gum \, \gm}
    \nonumber \\
&&
    + a_{97} \gwiione{\vw_{1}  \, \gua \, \gub} \gwii{\ga \, \gb \, \vw_{2} \, \gum \, \gm}
    + a_{98} \gwiione{\vw_{2}  \, \gua \, \gub} \gwii{\ga \, \gb \, \vw_{1} \, \gum \, \gm}
    \nonumber \\
&&
    + a_{99} \gwiione{ \gua \, \gub} \gwii{\ga \, \gb \, \vw_{1}  \, \vw_{2} \, \gum \, \gm}
    + a_{100} \gwiione{ \vw_{1} \, \vw_{2} \, \gua} \gwii{\ga \, \gub \, \gb \, \gum \, \gm}
    \nonumber \\
&&
    + a_{101} \gwiione{ \vw_{1}  \, \gua} \gwii{\ga \, \vw_{2} \, \gub \, \gb \, \gum \, \gm}
    + a_{102} \gwiione{  \vw_{2} \, \gua} \gwii{\ga \, \vw_{1} \, \gub \, \gb \, \gum \, \gm}
    \nonumber \\
&&
    + a_{103} \gwiione{  \gua} \gwii{\ga \, \vw_{1} \, \vw_{2} \, \gub \, \gb \, \gum \, \gm}
    + a_{104} \gwii{ \vw_{1} \, \vw_{2} \, \gua \, \ga \, \gub \, \gb \, \gum \, \gm}
    \nonumber \\
&& + a_{105} \gwiione{\gua} \gwiione{\ga \, \left\{ \vw_{1} \qp \vw_{2} \right\} \, \gb}
        \gwiione{\gub}
    \label{eqn:g3T2Ta}
\end{eqnarray}
where $\vw_i$ are arbitrary vector fields on the big phase space.
The following terms are omitted from this formula due to the lower genus equations:
\begin{eqnarray}
&&   b_{1} \gwiitwo{T(\vw_{1}) \vw_{2} \left\{\ga \qp \gua\right\}}
    + b_{2} \gwiitwo{\vw_{1} T(\vw_{2}) \left\{\ga \qp \gua\right\}}
        \nonumber \\
&&
    + b_{3} \gwiitwo{T(\vw_{1}) \, \gua \, \left\{\ga \qp \vw_{2} \right\}}
    + b_{4} \gwiitwo{T(\vw_{1}) \, \{ \vw_2 \qp \gua \}} \gwiione{ \ga }
        \nonumber \\
&&
    + b_{5} \gwiitwo{\gua}\gwii{\{\ga \qp \vw_1\} \, \vw_{2} \, \gub \, \gb}
    + b_{6} \gwiitwo{\gua}\gwii{\{ \ga \qp \vw_2\} \, \vw_{1} \, \gub \, \gb}
        \nonumber \\
&&
    + b_{7} \gwiitwo{\gua} \gwii{\ga \, \vw_{2} \, \gub \, \left\{ \gb \qp \vw_{1} \right\}}.
    \label{eqn:bterms}
\end{eqnarray}
The $b_{5}$, $b_{6}$, $b_{7}$ terms in this equation are eliminated using the first derivatives of
the WDVV equation.
To eliminate other terms, we need the Belorousski-Pandharipande equation (abbreviated as BP equation)
whose top order terms have the form
\begin{eqnarray}
 \sum_{\sigma \in S_3} \gwiitwo{{\cal W}_{\sigma(1)} \, T({\cal W}_{\sigma(2)} \qp {\cal W}_{\sigma(3)})}
- \gwiitwo{ T({\cal W}_{\sigma(1)}) \, \{{\cal W}_{\sigma(2)} \qp {\cal W}_{\sigma(3)}\}}
&=& {\rm L.O.T.} \label{eqn:BP}
\end{eqnarray}
for any vector fields $\vw_1, \vw_2, \vw_3$.
Throughout this paper, "{\bf top order terms}" refer to highest order derivatives of the highest genus
generating function involved in each equation. "{\bf L.O.T.}" is an abbreviation for "Lower Order Terms".
The $b_{4}$-term in equation \eqref{eqn:bterms} is eliminated by
using the BP equation applied to $\vw_{1}, \, \vw_{2}, \, \gua$.

Taking covariant derivative of the BP equation with respect to another vector field $\vw_{4}$ on the big phase space,
 we obtain
\begin{eqnarray}
 \sum_{\sigma \in S_3} \gwiitwo{\vw_{4} {\cal W}_{\sigma(1)} \, T({\cal W}_{\sigma(2)} \qp {\cal W}_{\sigma(3)})}
- \gwiitwo{ \vw_{4} T({\cal W}_{\sigma(1)}) \, \{{\cal W}_{\sigma(2)} \qp {\cal W}_{\sigma(3)}\}}
&=& {\rm L.O.T.} \label{eqn:DBP}
\end{eqnarray}
Setting $\vw_{3}=\ga$ and $\vw_{4}=\gua$ in equation \eqref{eqn:DBP}, we have
\begin{eqnarray}
&& \gwiitwo{ \gua \, \ga \, T(\vw_1 \qp \vw_2)}
    - \gwiitwo{ \gua \, T(\ga) \, \{ \vw_1 \qp \vw_2 \} }
    \nonumber \\
&&  + \gwiitwo{ \gua \, \vw_1 \, T(\ga \qp \vw_2)}
    - \gwiitwo{ \gua \, T(\vw_1) \, \{ \ga \qp \vw_2 \} }
    \nonumber \\
&&  + \gwiitwo{ \gua \, \vw_2 \, T(\ga \qp \vw_1)}
    - \gwiitwo{ \gua \, T(\vw_2) \, \{ \ga \qp \vw_1 \} }
=  {\, \, \rm L.O.T.} \label{eqn:DBP-1}
\end{eqnarray}
This equation can be used to eliminate the $b_{3}$-term in equation~\eqref{eqn:bterms}.

First setting $\vw_{2}=\ga$ and $\vw_{3}=\gua$ in equation \eqref{eqn:DBP}, then
changing $\vw_{4}$ to $\vw_{2}$ and summing over $\alpha$, we have
\begin{eqnarray}
\gwiitwo{ \vw_1 \, \vw_2 \, T(\gua \qp \ga)}
    - \gwiitwo{ \vw_2 \, T(\vw_1) \, \{ \gua \qp \ga \} }
&=&  {\rm L.O.T}. \label{eqn:DBP-2}
\end{eqnarray}
We can use this equation to eliminate the
 $b_{1}$-term in equation \eqref{eqn:bterms}.
Similarly the $b_{2}$-term  can be eliminated  by using the
 equation obtained from equation \eqref{eqn:DBP-2} after switching
$\vw_{1}$ and $\vw_{2}$.

Since equation~\eqref{eqn:g3T2Ta} holds for Gromov-Witten invariants of all
compact symplectic manifolds, it must be satisfied for Gromov-Witten invariants of
a point and ${\mathbb P}^{1}$. It is well known that
the Gromov-Witten invariants of these two manifolds satisfy
the Virasoro constraints (see \cite{W} \cite{K} for the point case, and
\cite{EHX} \cite{Gi} for the ${\mathbb P}^{1}$ case). A computer program for
calculating such invariants based on the Virasoro constraints
was written by Gathmann \cite{Ga}. We use Gathmann's program to compute
such invariants and plug in derivatives of
equation~\eqref{eqn:g3T2Ta} to find relations among constants
$a_1, \ldots, a_{105}$.

When the target manifold is a point, the degree of all stable maps
must be $0$. Hence, we will omit any reference to the degrees of
Gromov-Witten invariants of a point. In fact, such Gromov-Witten invariants
are just intersection numbers
\begin{eqnarray*}
\gwih{g}{\tau_{n_{1}} \, \cdots \, \tau_{n_{k}}}
    &=& \int_{\overline{\cal M}_{g, k}} \psi_1^{n_1} \cdots \psi_k^{n_k}
\end{eqnarray*}
over the moduli spaces $\overline{\cal M}_{g, k}$.
Since the cohomology space of a point is one dimensional,
coordinates on the big phase space are simply denoted by $t_{0},
t_{1}, t_{2}, \cdots$. We also identify vector fields $
\frac{\partial}{\partial t_{m}}$ with $\tau_{m}$ on the big phase
space.
Using Gromov-Witten invariants of a point, we obtained 43 linearly independent
relations among constants $a_1, \ldots, a_{105}$.
For example, from  $\Phi(\tau_{0}, \tau_{5})|_{t=0} = 0$, we obtain
\begin{align}
 \frac{a_2}{288}+\frac{a_{15}}{1152}+\frac{a_{16}}{288}+\frac{a_{24}}{1152}
    +\frac{a_{84}}{24}+\frac{a_{90}}{24}+\frac{a_{92}}{24}+a_{104} &= \frac{77}{414720}.
    \label{eqn:ptstart0}
 \end{align}
Similarly, 42 other linearly independent relations among $a_i$
are obtained from the vanishing of the following values of $\Phi$ and derivatives of
$\Phi$ at $t=0$:
\[ \begin{array}{llll}
\Phi(\tau_{1}, \tau_{4}), & \Phi(\tau_{2}, \tau_{3}), & \Phi(\tau_{3}, \tau_{2}),
            & \Phi(\tau_{4}, \tau_{1}), \\
\Phi(\tau_{5}, \tau_{0}), & \tau_{6}\Phi(\tau_{0}, \tau_{0}), &
        \tau_{5}\Phi(\tau_{0}, \tau_{1}), & \tau_{4}\Phi(\tau_{0}, \tau_{2}), \\
\tau_{3}\Phi(\tau_{0}, \tau_{3}), & \tau_{2}\Phi(\tau_{0}, \tau_{4}), &
    \tau_{5}\Phi(\tau_{1}, \tau_{0}), & \tau_{4}\Phi(\tau_{1}, \tau_{1}), \\
\tau_{3}\Phi(\tau_{1}, \tau_{2}), & \tau_{2}\Phi(\tau_{1}, \tau_{3}), &
    \tau_{4}\Phi(\tau_{2}, \tau_{0}), & \tau_{3}\Phi(\tau_{2}, \tau_{1}), \\
\tau_{2}\Phi(\tau_{2}, \tau_{2}), & \tau_{3}\Phi(\tau_{3}, \tau_{0}), &
    \tau_{2}\Phi(\tau_{3}, \tau_{1}), & \tau_{2}\Phi(\tau_{4}, \tau_{0}), \\
\tau_{3}\tau_{4}\Phi(\tau_{0}, \tau_{0}), & \tau_{2}\tau_{5}\Phi(\tau_{0}, \tau_{0}), &
    \tau_{3}\tau_{3}\Phi(\tau_{0}, \tau_{1}), & \tau_{2}\tau_{4}\Phi(\tau_{0}, \tau_{1}), \\
\tau_{2}\tau_{3}\Phi(\tau_{0}, \tau_{2}), & \tau_{2}\tau_{2}\Phi(\tau_{0}, \tau_{3}), &
    \tau_{3}\tau_{3}\Phi(\tau_{1}, \tau_{0}), & \tau_{2}\tau_{4}\Phi(\tau_{1}, \tau_{0}), \\
\tau_{2}\tau_{3}\Phi(\tau_{1}, \tau_{1}), & \tau_{2}\tau_{2}\Phi(\tau_{1}, \tau_{2}), &
    \tau_{2}\tau_{3}\Phi(\tau_{2}, \tau_{0}), & \tau_{2}\tau_{2}\Phi(\tau_{2}, \tau_{1}), \\
\tau_{2}\tau_{2}\Phi(\tau_{3}, \tau_{0}), & \tau_{2}\tau_{3}\tau_{3}\Phi(\tau_{0}, \tau_{0}), &
    \tau_{2}\tau_{2}\tau_{4}\Phi(\tau_{0}, \tau_{0}), & \tau_{2}\tau_{2}\tau_{3}\Phi(\tau_{0}, \tau_{1}), \\
\tau_{2}\tau_{2}\tau_{2}\Phi(\tau_{0}, \tau_{2}), \hspace{20pt} & \tau_{2}\tau_{2}\tau_{2}\Phi(\tau_{1}, \tau_{1}), \hspace{20pt} &
    \tau_{2}\tau_{2}\tau_{2}\Phi(\tau_{2}, \tau_{0}), \hspace{10pt} & \tau_{2}\tau_{2}\tau_{2}\tau_{3}\Phi(\tau_{0}, \tau_{0}), \\
\tau_{2}\tau_{2}\tau_{2}\tau_{2}\Phi(\tau_{0}, \tau_{1}), & \tau_{2}\tau_{2}\tau_{2}\tau_{2}\tau_{2}\Phi(\tau_{0}, \tau_{0}). & &
\end{array} \]
Explicit relations obtained in this way can be found in Appendix \ref{sec:point}.

When the target manifold is ${\mathbb CP}^{1}$, the degrees of the stable
  maps are indexed by $H_{2}({\mathbb CP}^{1}; {\mathbb Z}) \cong {\mathbb
  Z}$. The degree $d$ part of any equation for generating functions of
  Gromov-Witten invariants is the coefficient of $q^{d}$ in the Novikov ring.
We choose the basis $\{ \gamma_{0}, \gamma_{1} \}$ for $H^{*}({\mathbb CP}^{1}; {\mathbb C})$ with
$\gamma_{0} \in H^{0}({\mathbb CP}^{1}; {\mathbb C})$ being the identity of the ordinary
cohomology ring and $\gamma_{1} \in H^{2}({\mathbb CP}^{1}; {\mathbb C})$ the Poincare dual
to a point. Coordinates on the big phase space are denoted by
$\{ t_{n}^{0}, t_{n}^{1} \mid n \in {\mathbb Z}_{+} \}$. We identify
vector fields $\frac{\partial}{\partial t_{n}^{0}}$ and $\frac{\partial}{\partial t_{n}^{1}}$
with $\tau_{n, 0}$ and $\tau_{n, 1}$ respectively.
Using Gromov-Witten invariants of $\mathbb{P}^1$, we obtained 61 linearly independent
relations among constants $a_1, \ldots, a_{105}$.
For example, from the degree $0$ part of $\Phi(\tau_{0,0}, \tau_{2,1})|_{t=0} = 0$, we obtain
\begin{align}
 \frac{7 a_9}{138240} + \frac{7 a_{10}}{46080}+ \frac{a_{75}}{13824} +
 \frac{a_{79}}{13824}+ \frac{a_{82}}{13824} &= \frac{31}{96768}.
\label{eqn:P1start0}
 \end{align}

Similarly, 15 other linearly independent relations among $a_i$
are obtained from the vanishing of the degree-$0$ part of the following
values of $\Phi$ and derivatives of
$\Phi$ at $t=0$:
\[ \begin{array}{llll}
\Phi(\tau_{0,0}, \tau_{3,0}), & \Phi(\tau_{1,0}, \tau_{1,1}),  &
    \Phi(\tau_{1,0}, \tau_{2,0}), & \Phi(\tau_{1,1}, \tau_{1,0}), \\
\Phi(\tau_{2,0}, \tau_{0,1}), & \Phi(\tau_{3,0},\tau_{0,0}), &
    \tau_{3,1}\Phi(\tau_{0,0}, \tau_{0,0}), & \tau_{4,0}\Phi(\tau_{0,0}, \tau_{0,0}), \\
\tau_{2,1}\Phi(\tau_{0,0}, \tau_{1,0}), \hspace{20pt} & \tau_{3,0}\Phi(\tau_{0,0}, \tau_{1,0}), \hspace{20pt} &
    \tau_{2,0}\Phi(\tau_{0,0}, \tau_{2,0}), \hspace{20pt} & \tau_{2,1}\Phi(\tau_{1,0}, \tau_{0,0}), \\
\tau_{2,0}\Phi(\tau_{1,0}, \tau_{1,0}), & \tau_{2,1}\tau_{2,0}\Phi(\tau_{0,0}, \tau_{0,0}), &
    \tau_{2,0}\tau_{3,0}\Phi(\tau_{0,0}, \tau_{0,0}). &
\end{array} \]
45 more linearly independent relations among $a_i$ are
obtained from the vanishing of the degree-$1$ part of the following derivatives of
$\Phi$ at $t=0$:
\[ \begin{array}{llll}
\tau_{5,1}\Phi(\tau_{0,0}, \tau_{0,0}), & \tau_{4,1}\Phi(\tau_{0,0}, \tau_{0,1}), &
    \tau_{5,0}\Phi(\tau_{0,0}, \tau_{0,1}), & \tau_{4,1}\Phi(\tau_{0,0}, \tau_{1,0}), \\
\tau_{3,1}\Phi(\tau_{0,0}, \tau_{1,1}), & \tau_{4,0}\Phi(\tau_{0,0}, \tau_{1,1}), &
    \tau_{3,1}\Phi(\tau_{0,0}, \tau_{2,0}), & \tau_{2,1}\Phi(\tau_{0,0}, \tau_{2,1}), \\
\tau_{3,0}\Phi(\tau_{0,0}, \tau_{2,1}), & \tau_{2,1}\Phi(\tau_{0,0}, \tau_{3,0}), &
    \tau_{2,0}\Phi(\tau_{0,0}, \tau_{3,1}), & \tau_{4,1}\Phi(\tau_{0,1}, \tau_{0,0}), \\
\tau_{5,0}\Phi(\tau_{0,1}, \tau_{0,0}), & \tau_{3,1}\Phi(\tau_{0,1}, \tau_{0,1}), &
    \tau_{4,0}\Phi(\tau_{0,1}, \tau_{0,1}), & \tau_{3,1}\Phi(\tau_{0,1}, \tau_{1,0}), \\
\tau_{4,0}\Phi(\tau_{0,1}, \tau_{1,0}), & \tau_{2,1}\Phi(\tau_{0,1}, \tau_{1,1}), &
    \tau_{3,0}\Phi(\tau_{0,1}, \tau_{1,1}), & \tau_{2,1}\Phi(\tau_{0,1}, \tau_{2,0}), \\
\tau_{3,0}\Phi(\tau_{0,1}, \tau_{2,0}), & \tau_{2,0}\Phi(\tau_{0,1}, \tau_{2,1}), &
    \tau_{2,0}\Phi(\tau_{0,1}, \tau_{3,0}), & \tau_{4,1}\Phi(\tau_{1,0}, \tau_{0,0}), \\
\tau_{3,1}\Phi(\tau_{1,0}, \tau_{0,1}), & \tau_{4,0}\Phi(\tau_{1,0}, \tau_{0,1}), &
    \tau_{3,1}\Phi(\tau_{1,0}, \tau_{1,0}), & \tau_{2,1}\Phi(\tau_{1,0}, \tau_{1,1}), \\
\tau_{3,0}\Phi(\tau_{1,0}, \tau_{1,1}), & \tau_{2,1}\Phi(\tau_{1,0}, \tau_{2,0}), &
    \tau_{2,0}\Phi(\tau_{1,0}, \tau_{2,1}), & \tau_{3,1}\Phi(\tau_{1,1}, \tau_{0,0}), \\
\tau_{4,0}\Phi(\tau_{1,1}, \tau_{0,0}), & \tau_{2,1}\Phi(\tau_{1,1}, \tau_{0,1}), &
    \tau_{3,0}\Phi(\tau_{1,1}, \tau_{0,1}), & \tau_{2,1}\Phi(\tau_{1,1}, \tau_{1,0}), \\
\tau_{2,0}\Phi(\tau_{1,1}, \tau_{1,1}), & \tau_{2,0}\Phi(\tau_{1,1}, \tau_{2,0}), &
    \tau_{3,1}\Phi(\tau_{2,0}, \tau_{0,0}), & \tau_{3,1}\tau_{3,0}\Phi(\tau_{0,0}, \tau_{0,0}), \\
\tau_{2,1}\tau_{3,1}\Phi(\tau_{0,0}, \tau_{0,0}), \hspace{10pt} & \tau_{2,1}\tau_{4,0}\Phi(\tau_{0,0}, \tau_{0,0}), \hspace{10pt} &
    \tau_{2,0}\tau_{4,1}\Phi(\tau_{0,0}, \tau_{0,0}), \hspace{10pt} & \tau_{2,1}\tau_{2,1}\Phi(\tau_{0,0}, \tau_{0,1}),  \\
\tau_{2,1}\tau_{2,1}\Phi(\tau_{0,0}, \tau_{1,0}).  &
\end{array} \]
Explicit relations obtained in the above way can be found in Appendix \ref{sec:P1}.

Combining results from Gromov-Witten invariants of a point and  $\mathbb{P}^1$, we obtained
104 linearly independent relations
among $a_1, \ldots, a_{105}$. Using these relations we can solve all $a_i$ with $i\neq 2$
in terms of $a_2$, and obtain
the following formulas:

\allowdisplaybreaks
\begin{lem} \label{lem:coeff}
\[
\begin{array}{llll}
a_{1} = 5 ,  \hspace{40pt}  & a_2 \,\,\, {\rm is \,\,\, free}, & a_{3} =  0 ,
        \hspace{40pt} & a_{4} =  - \frac{1}{36} - 4 a_{2}  ,
  \\
 a_{5} =  \frac{5}{168} ,   & a_{6} =  \frac{1}{252} - 4 a_{2} ,
 & a_{7} = \frac{7}{24} + 6 a_{2} ,   & a_{8} =  - \frac{2}{3} - 120 a_{2} ,
  \\
 a_{9} =  \frac{5}{7} ,   & a_{10} =  \frac{2}{21} - 120 a_{2} ,
 & a_{11} = \frac{20}{3} + 240 a_{2} ,    & a_{12} =  \frac{100}{21} + 240 a_{2} ,
 \\
 a_{13} =  7 + 120 a_{2} ,   & a_{14} =  - \frac{1}{21} - 6 a_{2} ,
 & a_{15} = \frac{1}{14} ,    & a_{16} =  - \frac{1}{42} -6 a_{2} ,     \\
 a_{17} =  \frac{11}{14} + 24 a_{2} ,   & a_{18} =  \frac{1}{21} + 6 a_{2} ,
  & a_{19} = 0 ,   & a_{20} =  -\frac{100}{21} - 240 a_{2} ,      \\
 a_{21} =  0 ,   & a_{22} =  \frac{2}{3} + 120 a_{2} ,   & a_{23} =  \frac{1}{36} + 4a_{2} ,
   & a_{24} = -\frac{1}{252} + 4 a_{2} ,      \\
 a_{25} =  - \frac{1}{36} - 5 a_{2} ,   & a_{26} =  -\frac{47}{504} + 4 a_{2} ,
 & a_{27} =0 ,    & a_{28} =  -\frac{7}{12} - 12 a_{2} ,     \\
 a_{29} =  -\frac{11}{42} + 6 a_{2} ,   & a_{30} =  0 ,
& a_{31} = -\frac{47}{21} + 120 a_{2} ,   & a_{32} =  \frac{1}{168} - \frac{3}{5} a_{2} ,     \\
 a_{33} = - \frac{1}{72} - \frac{8}{5} a_{2} ,  & a_{34} =  -\frac{11}{504} - \frac{8}{5} a_{2} ,
 & a_{35} = \frac{97}{2520} + \frac{22}{5} a_{2} ,  &   \\
 a_{36} =  \frac{11}{360} + \frac{22}{5} a_{2} ,  & a_{37} =  -\frac{13}{560} - \frac{11}{10} a_{2} ,
  & a_{38} = \frac{89}{1680} + \frac{49}{10} a_{2},   &   \\
 a_{39} =  \frac{97}{1680} + \frac{49}{10} a_{2} ,  & a_{40} =  -\frac{61}{560} - \frac{71}{10} a_{2} ,
 & a_{41} = \frac{11}{210} + 6 a_{2} , &    \\
 a_{42} =  \frac{11}{210} + 6a_{2} , & a_{43} =  \frac{8}{105} + 6 a_{2} ,
 & a_{44} =  \frac{17}{210} + 6 a_{2} ,
  & a_{45} = -\frac{23}{70} - 24 a_{2} ,     \\
 a_{46} =  \frac{12}{35} + 18 a_{2} ,   & a_{47} =  \frac{33}{70} + 18 a_{2} ,
 &  a_{48} =  - \frac{41}{140} - 36 a_{2} ,  & a_{49} =  - \frac{12}{35} - 18 a_{2} ,    \\
 a_{50} =  - \frac{16}{105} - 12 a_{2} ,  & a_{51} =  \frac{19}{210} + 6 a_{2} ,
 & a_{52} =  \frac{8}{105} + 6 a_{2} ,
 & a_{53} =\frac{17}{210} + 6 a_{2} ,     \\
 a_{54} =  -\frac{3}{35} - 6 a_{2} ,   & a_{55} =  -\frac{17}{105} - 12 a_{2} ,
 & a_{56} = - \frac{1}{7} - 12 a_{2} ,   &  a_{57} =  \frac{9}{140} + \frac{36}{5} a_{2} ,    \\
 a_{58} =  - \frac{1}{140} - \frac{54}{5} a_{2} ,  & a_{59} =  \frac{11}{140} - \frac{54}{5} a_{2} ,
 & a_{60} = \frac{71}{140} + \frac{36}{5} a_{2} ,  & a_{61} = - \frac{13}{140} - \frac{54}{5} a_{2} ,
    \\
 a_{62} =  - \frac{7}{60} - \frac{54}{5} a_{2} ,  & a_{63} =  \frac{9}{28} + \frac{36}{5} a_{2},
 & a_{64} = \frac{71}{420} + \frac{36}{5} a_{2} ,   & a_{65} =  -\frac{3}{140} + \frac{36}{5} a_{2},
   \\
 a_{66} =  \frac{59}{1680} + \frac{77}{10} a_{2} ,  &  a_{67} =  \frac{37}{5040} - \frac{33}{10}a_{2} ,
 & a_{68} = \frac{73}{1008} - \frac{33}{10} a_{2} ,  &   \\
 a_{69} =  - \frac{17}{210} - \frac{44}{5} a_{2} ,  & a_{70} =  \frac{17}{140} + \frac{11}{10} a_{2} ,
 & a_{71} = \frac{76}{35} + 144 a_{2} ,   & a_{72} =  \frac{44}{35} + 144 a_{2} ,    \\
 a_{73} =  \frac{44}{35} + 144 a_{2} ,  & a_{74} =  \frac{64}{35} + 144 a_{2} ,
 & a_{75} =  \frac{68}{35} + 144 a_{2} ,
  & a_{76} =-\frac{72}{35} - 144 a_{2} ,      \\
 a_{77} =  -\frac{136}{35} - 288 a_{2} ,   &  a_{78} =  - \frac{24}{7} - 288 a_{2} ,
 & a_{79} = \frac{22}{35} + 96 a_{2} ,  & a_{80} =  \frac{8}{15} - 24 a_{2} ,     \\
 a_{81} =  \frac{148}{105} - 24 a_{2} ,  & a_{82} =  \frac{81}{70} + 12 a_{2} ,
 & a_{83} = - \frac{38}{35} - 168 a_{2},  & a_{84} =  -\frac{13}{13440} - \frac{1}{20} a_{2} ,
  \\
 a_{85} =  \frac{89}{40320} + \frac{1}{5} a_{2} ,
 & a_{86} = \frac{97}{40320} + \frac{1}{5} a_{2} ,
 &  a_{87} = - \frac{61}{13440} - \frac{3}{10} a_{2} , &   \\
 a_{88} =  - \frac{1}{315} - \frac{1}{4} a_{2} ,  & a_{89} =  \frac{1}{315} + \frac{1}{4} a_{2} ,
 & a_{90} =  \frac{17}{5040} + \frac{1}{4} a_{2} ,
 & a_{91} = - \frac{23}{1680} - a_{2} ,     \\
 a_{92} =  \frac{1}{630} + \frac{2}{15} a_{2} ,  & a_{93} = - \frac{1}{270} - \frac{8}{15} a_{2} ,
 & a_{94} = - \frac{1}{378} - \frac{8}{15} a_{2} , & a_{95} =  \frac{1}{42} + \frac{4}{5} a_{2} ,
   \\
 a_{96} =  \frac{1}{1344} + \frac{3}{10} a_{2} ,  & a_{97} = - \frac{1}{448} - \frac{9}{20} a_{2} ,
  & a_{98} = \frac{3}{2240} - \frac{9}{20} a_{2} , &     \\
 a_{99} =  \frac{43}{2240} + \frac{3}{10} a_{2} ,  &  a_{100} =  \frac{1}{2688} + \frac{3}{20} a_{2} ,
 & a_{101} = - \frac{5}{8064} - \frac{1}{10} a_{2} , &    \\
 a_{102} =  \frac{23}{40320} - \frac{1}{10} a_{2} ,  &  a_{103} =  \frac{11}{2688} + \frac{1}{40} a_{2} ,
& a_{104} = \frac{1}{23040}, & a_{105} = - \frac{64}{35} - 144 a_{2}.
  \end{array}
\]
\end{lem}
\begin{rem} We have checked the validity of Lemma \ref{lem:coeff} in many additional cases by direct calculation.  Please see
Appendix \ref{sec:check}
for details of a typical verification.
\end{rem}
After plugging these formulas into equation \eqref{eqn:g3T2Ta},
the coefficient of $a_{2}$ is an expression totally symmetric with respect to
$\vw_1$ and $\vw_2$, which can be written as
\[ \Omega(\vw_1, \vw_2) + \Omega(\vw_2, \vw_1) \]
where $\Omega$ is a tensor described in the following way:
First we rewrite function $\Phi(\vw_1, \vw_2)$ defined in equation \eqref{eqn:g3T2Ta} as
\begin{equation*}  
\Phi(\vw_1, \vw_2)
\, = \,  -\gwiih{3}{T^{2}(\vw_{1}) T(\vw_{2})}
 + \sum_{i=1}^{105} \, a_{i} \, \Theta_i
\end{equation*}
where $\Theta_i$ is the coefficient of $a_i$ in
equation \eqref{eqn:g3T2Ta}.
For example,
\[ \Theta_1 = \gwiih{3}{T^{2}(\vw_{1} \qp \vw_{2})}, \hspace{20pt}
\Theta_4 = \gwiitwo{\vw_{1} \, T(\gua)} \gwii{\ga \, \vw_{2} \, \gub \, \gb} .
\]
The tensor $\Omega$ can then be defined by
\begin{eqnarray}
\Omega(\vw_1, \vw_2)
&:=&  \frac{1}{2} \, \Theta_2
-6 \, \Theta_{14}
   + 3 \, \Theta_{18}
\nonumber \\
&&    -4  \, \Theta_{4}  + 3 \, \Theta_{7}
    - 120 \, \Theta_{8}
    + 240 \, \Theta_{11}
     + 60 \, \Theta_{13}
    + 12 \, \Theta_{17}
     -120 \, \Theta_{20}
    + 60 \, \Theta_{22}
    \nonumber \\
&&
    + 4 \, \Theta_{23}
   -\frac{5}{2} \, \Theta_{25}
     + 2 \, \Theta_{26}
     -6 \, \Theta_{28} + 3 \, \Theta_{29}
      + 60 \, \Theta_{31}
    - \frac{3}{10} \, \Theta_{32}
     - \frac{8}{5} \, \Theta_{33}
    + \frac{22}{5} \, \Theta_{35}
    \nonumber \\
&&
    - \frac{11}{20} \, \Theta_{37}
    + \frac{49}{10} \, \Theta_{38}
       - \frac{71}{20} \, \Theta_{40}
    + 6  \, \Theta_{41}
    + 6 \, \Theta_{43}
    - 12 \, \Theta_{45}
      + 18 \, \Theta_{46}
    - 18 \, \Theta_{48}
   \nonumber \\
&&    - 9 \, \Theta_{49}
    - 6 \, \Theta_{50}
      + 3 \, \Theta_{51}
    + 6 \, \Theta_{52}
    - 3 \, \Theta_{54}
    - 12 \, \Theta_{55}
       + \frac{18}{5} \, \Theta_{57}
    - \frac{54}{5} \, \Theta_{58}
      + \frac{18}{5} \, \Theta_{60}
       \nonumber \\
&&
    - \frac{54}{5} \, \Theta_{61}
       + \frac{36}{5} \, \Theta_{63}
    + \frac{18}{5} \, \Theta_{65}
      + \frac{77}{20} \, \Theta_{66}
    - \frac{33}{10} \, \Theta_{67}
      - \frac{22}{5}  \, \Theta_{69}
    + \frac{11}{20} \, \Theta_{70}
      +  72 \, \Theta_{71}
      \nonumber \\
&&
    + 144 \, \Theta_{72}
     + 144 \, \Theta_{74}
       - 72 \, \Theta_{76}
      - 288 \, \Theta_{77}
    + 48 \, \Theta_{79}
      -24 \, \Theta_{80}
      + 6 \, \Theta_{82}
    - 84 \, \Theta_{83}
    \nonumber \\
&&
    - \frac{1}{40} \, \Theta_{84}
       + \frac{1}{5} \, \Theta_{85}
    - \frac{3}{20} \, \Theta_{87}
       - \frac{1}{8} \, \Theta_{88}
    + \frac{1}{4} \, \Theta_{89}
     - \frac{1}{2} \, \Theta_{91}
    + \frac{1}{15} \, \Theta_{92}
      - \frac{8}{15} \, \Theta_{93}
      \nonumber \\
&&
    + \frac{2}{5} \, \Theta_{95}
      + \frac{3}{20} \, \Theta_{96}
    - \frac{9}{20} \, \Theta_{97}
      + \frac{3}{20} \, \Theta_{99}
    + \frac{3}{40} \, \Theta_{100}
      - \frac{1}{10} \, \Theta_{101}
    + \frac{1}{80} \, \Theta_{103}
    \nonumber \\
&& - 72 \, \Theta_{105}.
    \label{eqn:a2Coef}
\end{eqnarray}
Note that $\Omega(\vw_1, \vw_2)$ only involves
data up to genus-$2$. The top order terms of $\Omega(\vw_1, \vw_2)$ is given by
\begin{eqnarray*}
 \frac{1}{2} \, \Theta_2
-6 \, \Theta_{14}
   + 3 \, \Theta_{18}
&=& \frac{1}{2} \gwiitwo{ \vw_{1} \, \vw_{2} \,  T(\ga \qp \gua)}
-6 \gwiitwo{\vw_{1} \, T(\gua \qp \vw_{2}) \, \ga }
    \nonumber \\
&&    + 3 \gwiitwo{T(\gua) \, \ga \, \left\{\vw_{1} \qp \vw_{2} \right\}}.
\end{eqnarray*}

Now we are ready to prove the two main results of this paper:

\begin{thm} \label{pro:g2uni}
For Gromov-Witten invariants of any compact symplectic manifold, we have
the following genus-2 universal equation
\begin{equation}  \label{eqn:a2Eqn?}
\Omega(\vw_1, \vw_2) + \Omega(\vw_2, \vw_1) = 0
\end{equation}
for all vector fields $\vw_1$ and $\vw_2$.
\end{thm}

\begin{thm} \label{thm:g3T2Ta2}
For $i=1, \cdots, 105$, let $\tilde{a}_i$ be the constant obtained from $a_i$ in Lemma \ref{lem:coeff} by setting $a_2 =0$.
For Gromov-Witten invariants of any compact symplectic manifold, we have
the following genus-3 universal equation
\allowdisplaybreaks
\begin{eqnarray}
  \gwiih{3}{T^{2}(\vw_{1}) T(\vw_{2})}
&=& \sum_{i=1}^{105} \, \tilde{a}_i \, \Theta_i
    \label{eqn:g3T2Ta2}
\end{eqnarray}
where $\vw_1$ and $\vw_2$ are arbitrary vector fields on the big phase space.
\end{thm}

{\bf Proof of Theorem \ref{pro:g2uni} and  Theorem \ref{thm:g3T2Ta2}}:
By Lemma~\ref{lem:coeff} and discussions following that lemma, we have
\begin{eqnarray}
\gwiih{3}{T^{2}(\vw_{1}) T(\vw_{2})}
&=&    a_2 \, \{ \Omega(\vw_1, \vw_2) + \Omega(\vw_2, \vw_1) \} + \sum_{i=1}^{105} \, \tilde{a}_i \, \Theta_i
    \label{eqn:g3T2T_with_a2}
\end{eqnarray}
where $a_2$ is a constant.

Each term in equation~\eqref{eqn:g3T2Ta} corresponds to an element
in $R^3(\overline{\cal M}_{3,2})$. So terms in equation~\eqref{eqn:g3T2Ta}
with coefficients $a_1, \ldots, a_{105}$
give 105 elements in the tautological group $R^3(\overline{\cal M}_{3,2})$.
It has been computed by Bergstr\"{o}m in \cite{B} that the corresponding cohomology group
on $\overline{\cal M}_{3,2}$ has dimension 104. Therefore
 the 105 cohomological classes corresponding to terms in equation~\eqref{eqn:g3T2Ta}
with coefficients $a_1, \ldots, a_{105}$ must be linearly dependent.
So there must exist a linear relation among these 105 terms.
This linear relation gives a universal equation for Gromov-Witten invariants.
We can add freely
any scalar multiplication of this equation to equation~\eqref{eqn:g3T2Ta}.
By equation \eqref{eqn:g3T2T_with_a2}, the only freedom in equation~\eqref{eqn:g3T2Ta}
is
\[ a_2 \, \left\{  \Omega(\vw_1, \vw_2) + \Omega(\vw_2, \vw_1) \right\}. \]
 Therefore we must have
\[ \Omega(\vw_1, \vw_2) + \Omega(\vw_2, \vw_1) =0. \]
 This proves
Theorem~\ref{pro:g2uni}.   Theorem \ref{thm:g3T2Ta2} then follows from
equation \eqref{eqn:g3T2T_with_a2}.
$\Box$

{\bf Proof of Theorem~\ref{thm:g3T2Tskew}}:
We can use Theorem~\ref{thm:g3T2Ta2} to compute
$\gwiih{3}{T^{2}(\vw_{1}) T(\vw_{2})}$ and $\gwiih{3}{T^{2}(\vw_{2}) T(\vw_{1})}$
separately. The difference of these two terms then proves Theorem~\ref{thm:g3T2Tskew}.
Note that all symmetric terms in Theorem~\ref{thm:g3T2Ta2} are canceled out
in
\[ \gwiih{3}{T^{2}(\vw_{1}) T(\vw_{2})}-\gwiih{3}{T^{2}(\vw_{2}) T(\vw_{1})}.\]
This makes Theorem~\ref{thm:g3T2Tskew} a much simpler formula than
Theorem~\ref{thm:g3T2Ta2}.
$\Box$

{\bf Remark}:
We would like to point out that
Theorem~\ref{thm:g3T2Tskew} can be proved by using equation \eqref{eqn:g3T2T_with_a2}
which is weaker than Theorem~\ref{thm:g3T2Ta2}. In particular, the proof of
Theorem~\ref{thm:g3T2Tskew} can be done without using results in \cite{B}.
The reason is that
\[ a_2 \, \{ \Omega(\vw_1, \vw_2)+ \Omega(\vw_2, \vw_1) \} \]
 is  symmetric with respect to
$\vw_1$ and $\vw_2$. So this term is canceled in
\[ \gwiih{3}{T^{2}(\vw_{1}) T(\vw_{2})}-\gwiih{3}{T^{2}(\vw_{2}) T(\vw_{1})} \]
when using equation \eqref{eqn:g3T2T_with_a2}.
$\Box$

\section{Comparison with known universal equations}
\label{sec:compare}

For a vector field $\vw$ with level of descendant $n$, $T(\vw)$ is a vector field with highest level of
descendant equal to $n+1$. So applying the operator $T$ will increase the level of descendant by $1$.
Universal equations for Gromov-Witten invariants often involve the operator $T$. We call the highest total degree
of $T$ in each term of a universal equation the {\it level of descendant for the universal equation}.
For example, the level of
descendant for the equation in Theorem~\ref{pro:g2uni} is equal to $1$ and the level of descendant for
equations in Theorem~\ref{thm:g3T2Tskew} and Theorem~\ref{thm:g3T2Ta2} is equal to $3$.
The level of descendant is a very important numerical quantity for studying relations among universal equations.
In general, it is very easy to increase the level of descendant by simply replacing a vector field $\vw$ by
$T(\vw)$ in a universal equation. However, it is much harder to decrease the level of descendant of a equation.
The only way to do this without sacrificing universality is to use the string equation. Let
\[ {\cal S} := - \sum_{m, \alpha} \tilde{t}^{\alpha}_{m}
        \grav{m-1}{\alpha} \]
be the {\it string vector field}, where
$\tilde{t}^{\alpha}_{m} = t^{\alpha}_{m} - \delta_{\alpha, 1} \delta_{m, 1}$.
Then the string equation implies
\[ \gwiig{\vs \vw_1 \cdots \vw_k} = \sum_{i=1}^{k} \gwiig{\vw_1 \cdots \tau_{-}(\vw_i) \cdots \vw_k} \]
for $g>0$. Sometimes, this formula can be used to decrease the level of descendant of a universal equation
when at least two arbitrary vector fields are involved in the equation. However it should be noticed that
taking the derivative of a universal equation along the string vector field and then applying the string equation
does not produce new equation. In many cases, setting an arbitrary vector field
in a universal equation to be the string vector field and then applying
the string equation often result in a trivial equation.

We say that a universal equation is {\it new} if it can not be
obtained from known universal equations by taking derivatives or other algebraic manipulations.
Below we discuss relations between previously known universal equations and the equations in
Theorem~\ref{pro:g2uni} and Theorem~\ref{thm:g3T2Ta2}.

We can see that the genus-2 equation in Theorem~\ref{pro:g2uni} is new by simply considering the top order terms
which has the following form:
\begin{eqnarray}
\gwiitwo{\vw_1 \, \vw_2 \, T(\ga \qp \gua) } + 6 \gwiitwo{ \{ \vw_1 \qp \vw_2 \} \, T(\gua) \, \ga}
 &&
\nonumber \\
 - 6 \gwiitwo{ \vw_1 \,  \gua  \, T(\vw_2 \qp\ga) } - 6 \gwiitwo{ \vw_2 \,  \gua \, T(\vw_1 \qp \ga) }
&=& {\rm L.O.T.} \label{eqn:g2unitop}
\end{eqnarray}
There are 3 genus-2 universal equations known before.
The 2 genus-2 equations given in \cite{Ge2} (one of them used Mumford relation) have level of descendant equal to $2$.
Applying string equation to these equations does not produce non-trivial equations.
Our genus-2 equation in Theorem~\ref{pro:g2uni} has level of descendant equal to $1$.
Therefore equations in
\cite{Ge2} do not help in deriving equation~\eqref{eqn:g2unitop}. Hence we only need to consider
the BP equation whose top order terms have the form as given in equation \eqref{eqn:BP}.
Since equation~\eqref{eqn:g2unitop} involves the third order derivatives of genus-2 generating
function and BP equation only involves the second order derivatives of genus-2 generating function,
we need to consider the derivative of the BP equation whose top order terms have the form given
in equation~\eqref{eqn:DBP}. To compare with equation~\eqref{eqn:g2unitop},
we need to set two of the vector fields in equation~\eqref{eqn:DBP} to $\gua$ and $\ga$ respectively.
Using symmetry of equation~\eqref{eqn:DBP}, we can only obtain equations \eqref{eqn:DBP-1} and
\eqref{eqn:DBP-2} after relabeling the remaining vector fields.
Observe that top order terms in equation~\eqref{eqn:g2unitop} consists of exactly half of top order terms
from equation \eqref{eqn:DBP-1} and half of top order terms from equation \eqref{eqn:DBP-2}. But we can not
derive equation~\eqref{eqn:g2unitop} from equations \eqref{eqn:DBP-1} and
\eqref{eqn:DBP-2}. To see this, we observe that we have to use equation \eqref{eqn:DBP-2}
to obtain the term $\gwiitwo{\vw_1 \, \vw_2 \, T(\gua \qp \ga)}$
in equation~\eqref{eqn:g2unitop}. But this would produce
a new term $\gwiitwo{\vw_2 \, T(\vw_1) \, \{\gua \qp \ga \}}$ which can not be derived from
equation \eqref{eqn:DBP-1}. This indicates that the genus-2 universal equation in
Theorem~\ref{pro:g2uni} is indeed new.

Now we consider the genus-3 universal equation given in Theorem~\ref{thm:g3T2Ta2} whose top
order terms has the form
\begin{equation} \label{eqn:g3unitop}
 \gwiih{3}{T^2(\vw_1) \, T(\vw_2) } \,\, = \,\, {\rm L.O.T.}
\end{equation}
This equation has level of
descendant equal to $3$. There are other genus-3 universal equations
proved in \cite{KL} and \cite{LP}. The levels of descendant for genus-3 equations in \cite{LP}
are at least $6$, therefore these equations can not be used to derive the equation in
Theorem~\ref{thm:g3T2Ta2}. The top order terms of genus-3 equation in \cite{KL} has the form
\begin{equation} \label{eqn:g3Ttop} \gwiih{3}{T^3(\vw)} \,\, =  \,\, {\rm L.O.T.}
 \end{equation}
The highest power of $T$ involved in this equation is $3$ while the highest power of $T$ applied
to each individual vector in
equation \eqref{eqn:g3unitop} is $2$.
Therefore the genus-3 equation in Theorem~\ref{thm:g3T2Ta2}
can not be derived from the genus-3 equation in \cite{KL}.

On the other hand, it turns out that the genus-$3$ topological recursion relation
in \cite{KL} can be derived using Theorem~\ref{thm:g3T2Ta2}. To see this,
we first observe that
\[ T(\vs) = \vd := - \sum_{m, \alpha} \tilde{t}^{\alpha}_{m} \, \,
\grava{m} \]
is the dilaton vector field (cf. \cite{L1} for a proof of this fact).
Applying Theorem~\ref{thm:g3T2Ta2} for $\vw_1 = T(\vw)$ and $\vw_2 = \vs$,
then getting rid of the string and dilaton vector fields by using
the string and dilaton equations, we can obtain the main result of \cite{KL}.
Note that in this derivation, we can also
use equation \eqref{eqn:g3T2T_with_a2} instead of Theorem~\ref{thm:g3T2Ta2}.
In fact we can prove
\[ \Omega(\vw, \vs)+ \Omega(\vs, \vw) = \Omega(\vw, \vd)+ \Omega(\vd, \vw) = 0 \]
directly using the string and dilaton equations.
Therefore equation \eqref{eqn:g3T2T_with_a2} and Theorem~\ref{thm:g3T2Ta2} give the
same result for this special case.
The proofs of these
facts are quite long and are omitted here.

{\bf Remark}:
Universal equations in Theorems~\ref{thm:g3T2Tskew}, \ref{pro:g2uni}, and \ref{thm:g3T2Ta2}
correspond
to 3 relations in the tautological ring of $\overline{\cal M}_{3,2}$ which do not involve
$\kappa$ classes.
Since the primary goal in this paper is to find universal equations for Gromov-Witten invariants,
we omit the explicit forms of relations corresponding to Theorems \ref{pro:g2uni}, and \ref{thm:g3T2Ta2}.
One can find the dual graph representation
for the relation corresponding to Theorem~\ref{thm:g3T2Tskew}
in Appendix \ref{sec:dualgraph}.
Since the proof of Theorem~\ref{pro:g2uni} and Theorem~\ref{thm:g3T2Ta2} needs a result in \cite{B},
the proofs of the corresponding relations are only at
the level of cohomology.
As explained in the remark
following the proof of Theorems~\ref{thm:g3T2Tskew}, the proof of the relation corresponding to
Theorems~\ref{thm:g3T2Tskew} is valid in the tautological Chow ring.
We also note that there is a codimension
2 relation involving $\kappa$ classes in the tautological ring of
$\overline{\cal M}_{3,2}$ due to Polito \cite{Po}.
It is apparent that the genus-2
relation corresponding to Theorem~\ref{pro:g2uni} is independent of the relation
in \cite{Po} and other known genus-2 relations as explained before.  At this point, it is not clear how the other two relations found in this paper are related to
Polito's result since all of these relations are very complicated.
After the first version of this paper was posted on the arXiv, Pandharipande, Pixton, and Zvonkine
have recently obtained a class of relations of all genera in the tautological cohomology rings in \cite{PPZ}. These relations again contain $\kappa$ classes. In order to get universal equations for Gromov-Witten invariants from these relations, $\kappa$ classes have to be eliminated.
It would be interesting to derive universal equations for Gromov-Witten invariants from Polito's result and relations in \cite{PPZ} and compare them with ours.
We may investigate these issues in the future.

\appendix
\vspace{30pt}
\hspace{160pt} {\bf \Large Appendix}

\section{Relations among constants $a_1, \ldots, a_{105}$}

The formulas in Lemma~\ref{lem:coeff} are obtained by solving
a system of 104 linearly independent relations among constants $a_1, \ldots, a_{105}$
in equation~\eqref{eqn:g3T2Ta}. These relations are obtained by using Gromov-Witten theory
of a point and $\mathbb{P}^1$. In this appendix, we list all these relations and indicate
how they are obtained.

\subsection{Relations from the Gromov-Witten invariants of a point}
\label{sec:point}

From  $\Phi(\tau_{0}, \tau_{5})|_{t=0} = 0$, we obtain
\begin{align}
 0 &= \frac{a_2}{288}+\frac{a_{15}}{1152}+\frac{a_{16}}{288}+\frac{a_{24}}{1152}+\frac{a_{84}}{24}+\frac{a_{90}}{24}+\frac{a_{92}}{24}+a_{104}\nonumber \\ &-\frac{77}{414720}.
\label{eqn:ptstart}
 \end{align}
From  $\Phi(\tau_{1}, \tau_{4})|_{t=0} = 0$, we obtain
\begin{align}
 0 &= \frac{a_2}{96}+\frac{a_5}{1152}+\frac{a_6}{384}+\frac{a_{84}}{6}+\frac{a_{86}}{24}+\frac{a_{92}}{6}+\frac{a_{94}}{24}+5 a_{104}\nonumber \\ &-\frac{503}{1451520}.
 \end{align}
From  $\Phi(\tau_{2}, \tau_{3})|_{t=0} = 0$, we obtain
\begin{align}
 0 &= \frac{29 a_2}{1440}+\frac{a_{33}}{576}+\frac{a_{62}}{576}+\frac{7 a_{84}}{24}+\frac{7 a_{92}}{24}+\frac{a_{98}}{24}+\frac{a_{101}}{24}+10 a_{104}\nonumber \\ &-\frac{607}{1451520}.
 \end{align}
From  $\Phi(\tau_{3}, \tau_{2})|_{t=0} = 0$, we obtain
\begin{align}
 0 &= \frac{29 a_2}{1440}+\frac{a_{34}}{576}+\frac{a_{61}}{576}+\frac{7 a_{84}}{24}+\frac{7 a_{92}}{24}+\frac{a_{97}}{24}+\frac{a_{102}}{24}+10 a_{104}\nonumber \\ &-\frac{503}{1451520}.
 \end{align}
From  $\Phi(\tau_{4}, \tau_{1})|_{t=0} = 0$, we obtain
\begin{align}
 0 &= \frac{a_2}{96}+\frac{a_3}{1152}+\frac{a_4}{384}+\frac{a_{84}}{6}+\frac{a_{85}}{24}+\frac{a_{92}}{6}+\frac{a_{93}}{24}+5 a_{104}\nonumber \\ &-\frac{77}{414720}.
 \end{align}
From  $\Phi(\tau_{5}, \tau_{0})|_{t=0} = 0$, we obtain
\begin{align}
 0 &= \frac{a_2}{288}+\frac{a_{14}}{288}+\frac{a_{23}}{1152}+\frac{a_{84}}{24}+\frac{a_{89}}{24}+\frac{a_{92}}{24}+a_{104}-\frac{5}{82944}.
 \end{align}
From  $\tau_{6}\Phi(\tau_{0}, \tau_{0})|_{t=0} = 0$, we obtain
\begin{align}
 0 &= \frac{77 a_1}{414720}+\frac{5 a_2}{1152}+\frac{5 a_{14}}{1152}+\frac{5 a_{15}}{1152}+\frac{5 a_{16}}{1152}+\frac{5 a_{18}}{1152}+\frac{5 a_{19}}{1152}+\frac{a_{23}}{1152}\nonumber \\ &+\frac{a_{24}}{1152}+\frac{a_{25}}{1152}+\frac{a_{29}}{1152}+\frac{a_{84}}{24}+\frac{a_{88}}{24}+\frac{a_{89}}{24}+\frac{a_{90}}{24}+\frac{a_{92}}{24}\nonumber \\ &+a_{104}-\frac{77}{69120}.
 \end{align}
From  $\tau_{5}\Phi(\tau_{0}, \tau_{1})|_{t=0} = 0$, we obtain
\begin{align}
 0 &= \frac{5 a_2}{288}+\frac{a_3}{288}+\frac{a_4}{288}+\frac{a_{15}}{90}+\frac{5 a_{16}}{288}+\frac{a_{17}}{288}+\frac{a_{24}}{288}+\frac{a_{28}}{1152}\nonumber \\ &+\frac{a_{30}}{1152}+\frac{5 a_{84}}{24}+\frac{a_{85}}{24}+\frac{5 a_{90}}{24}+\frac{a_{91}}{24}+\frac{5 a_{92}}{24}+\frac{a_{93}}{24}+6 a_{104}\nonumber \\ &-\frac{17}{5760}.
 \end{align}
From  $\tau_{4}\Phi(\tau_{0}, \tau_{2})|_{t=0} = 0$, we obtain
\begin{align}
 0 &= \frac{11 a_2}{288}+\frac{a_7}{384}+\frac{a_{11}}{9216}+\frac{17 a_{15}}{960}+\frac{11 a_{16}}{288}+\frac{11 a_{24}}{1440}+\frac{a_{34}}{576}+\frac{a_{42}}{576}\nonumber \\ &+\frac{a_{55}}{576}+\frac{a_{61}}{576}+\frac{11 a_{84}}{24}+\frac{a_{87}}{24}+\frac{11 a_{90}}{24}+\frac{11 a_{92}}{24}+\frac{a_{95}}{24}+\frac{a_{97}}{24}\nonumber \\ &+\frac{a_{102}}{24}+15 a_{104}-\frac{1121}{241920}.
 \end{align}
From  $\tau_{3}\Phi(\tau_{0}, \tau_{3})|_{t=0} = 0$, we obtain
\begin{align}
 0 &= \frac{29 a_2}{576}+\frac{17 a_{15}}{960}+\frac{29 a_{16}}{576}+\frac{29 a_{24}}{2880}+\frac{a_{32}}{576}+\frac{a_{33}}{576}+\frac{a_{47}}{576}+\frac{a_{62}}{576}\nonumber \\ &+\frac{a_{65}}{576}+\frac{7 a_{84}}{12}+\frac{7 a_{90}}{12}+\frac{7 a_{92}}{12}+\frac{a_{98}}{24}+\frac{a_{99}}{24}+\frac{a_{100}}{24}+\frac{a_{101}}{24}\nonumber \\ &+20 a_{104}-\frac{1121}{241920}.
 \end{align}
From  $\tau_{2}\Phi(\tau_{0}, \tau_{4})|_{t=0} = 0$, we obtain
\begin{align}
 0 &= \frac{11 a_2}{288}+\frac{a_5}{1152}+\frac{a_6}{384}+\frac{a_9}{27648}+\frac{a_{10}}{9216}+\frac{a_{15}}{90}+\frac{11 a_{16}}{288}+\frac{11 a_{24}}{1440}\nonumber \\ &+\frac{a_{37}}{576}+\frac{a_{44}}{576}+\frac{a_{53}}{576}+\frac{a_{57}}{576}+\frac{11 a_{84}}{24}+\frac{a_{86}}{24}+\frac{11 a_{90}}{24}+\frac{11 a_{92}}{24}\nonumber \\ &+\frac{a_{94}}{24}+\frac{a_{96}}{24}+\frac{a_{103}}{24}+15 a_{104}-\frac{17}{5760}.
 \end{align}
From  $\tau_{5}\Phi(\tau_{1}, \tau_{0})|_{t=0} = 0$, we obtain
\begin{align}
 0 &= \frac{5 a_2}{288}+\frac{a_5}{288}+\frac{a_6}{288}+\frac{5 a_{14}}{288}+\frac{a_{17}}{288}+\frac{a_{23}}{288}+\frac{a_{27}}{1152}+\frac{a_{28}}{1152}\nonumber \\ &+\frac{a_{30}}{1152}+\frac{5 a_{84}}{24}+\frac{a_{86}}{24}+\frac{5 a_{89}}{24}+\frac{a_{91}}{24}+\frac{5 a_{92}}{24}+\frac{a_{94}}{24}+6 a_{104}\nonumber \\ &-\frac{503}{241920}.
 \end{align}
From  $\tau_{4}\Phi(\tau_{1}, \tau_{1})|_{t=0} = 0$, we obtain
\begin{align}
 0 &= \frac{5 a_2}{96}+\frac{11 a_3}{1440}+\frac{a_4}{96}+\frac{11 a_5}{1440}+\frac{a_6}{96}+\frac{a_7}{192}+\frac{5 a_{84}}{6}+\frac{a_{85}}{6}\nonumber \\ &+\frac{a_{86}}{6}+\frac{a_{87}}{12}+\frac{5 a_{92}}{6}+\frac{a_{93}}{6}+\frac{a_{94}}{6}+\frac{a_{95}}{12}+30 a_{104}-\frac{1121}{241920}.
 \end{align}
From  $\tau_{3}\Phi(\tau_{1}, \tau_{2})|_{t=0} = 0$, we obtain
\begin{align}
 0 &= \frac{29 a_2}{288}+\frac{29 a_5}{2880}+\frac{29 a_6}{1440}+\frac{a_{32}}{288}+\frac{a_{34}}{192}+\frac{a_{36}}{576}+\frac{a_{61}}{192}+\frac{a_{63}}{576}\nonumber \\ &+\frac{a_{65}}{288}+\frac{35 a_{84}}{24}+\frac{7 a_{86}}{24}+\frac{35 a_{92}}{24}+\frac{7 a_{94}}{24}+\frac{a_{97}}{8}+\frac{a_{99}}{8}+\frac{a_{100}}{12}\nonumber \\ &+\frac{a_{102}}{6}+60 a_{104}-\frac{583}{96768}.
 \end{align}
From  $\tau_{2}\Phi(\tau_{1}, \tau_{3})|_{t=0} = 0$, we obtain
\begin{align}
 0 &= \frac{29 a_2}{288}+\frac{11 a_5}{1440}+\frac{29 a_6}{1440}+\frac{a_{33}}{288}+\frac{a_{37}}{192}+\frac{a_{39}}{576}+\frac{a_{57}}{192}+\frac{a_{59}}{576}\nonumber \\ &+\frac{a_{62}}{288}+\frac{35 a_{84}}{24}+\frac{7 a_{86}}{24}+\frac{35 a_{92}}{24}+\frac{7 a_{94}}{24}+\frac{a_{96}}{8}+\frac{a_{98}}{8}+\frac{a_{101}}{12}\nonumber \\ &+\frac{a_{103}}{6}+60 a_{104}-\frac{1121}{241920}.
 \end{align}
From  $\tau_{4}\Phi(\tau_{2}, \tau_{0})|_{t=0} = 0$, we obtain
\begin{align}
 0 &= \frac{11 a_2}{288}+\frac{a_7}{384}+\frac{a_{12}}{9216}+\frac{11 a_{14}}{288}+\frac{11 a_{23}}{1440}+\frac{a_{33}}{576}+\frac{a_{41}}{576}+\frac{a_{56}}{576}\nonumber \\ &+\frac{a_{62}}{576}+\frac{11 a_{84}}{24}+\frac{a_{87}}{24}+\frac{11 a_{89}}{24}+\frac{11 a_{92}}{24}+\frac{a_{95}}{24}+\frac{a_{98}}{24}+\frac{a_{101}}{24}\nonumber \\ &+15 a_{104}-\frac{607}{241920}.
 \end{align}
From  $\tau_{3}\Phi(\tau_{2}, \tau_{1})|_{t=0} = 0$, we obtain
\begin{align}
 0 &= \frac{29 a_2}{288}+\frac{29 a_3}{2880}+\frac{29 a_4}{1440}+\frac{a_{32}}{288}+\frac{a_{33}}{192}+\frac{a_{35}}{576}+\frac{a_{62}}{192}+\frac{a_{64}}{576}\nonumber \\ &+\frac{a_{65}}{288}+\frac{35 a_{84}}{24}+\frac{7 a_{85}}{24}+\frac{35 a_{92}}{24}+\frac{7 a_{93}}{24}+\frac{a_{98}}{8}+\frac{a_{99}}{8}+\frac{a_{100}}{12}\nonumber \\ &+\frac{a_{101}}{6}+60 a_{104}-\frac{1121}{241920}.
 \end{align}
From  $\tau_{2}\Phi(\tau_{2}, \tau_{2})|_{t=0} = 0$, we obtain
\begin{align}
 0 &= \frac{7 a_2}{48}+\frac{a_{33}}{144}+\frac{a_{34}}{144}+\frac{a_{37}}{144}+\frac{a_{57}}{144}+\frac{a_{61}}{144}+\frac{a_{62}}{144}+\frac{a_{67}}{576}\nonumber \\ &+\frac{a_{68}}{576}+\frac{a_{69}}{576}+\frac{a_{83}}{13824}+2 a_{84}+2 a_{92}+\frac{a_{96}}{6}+\frac{a_{97}}{6}+\frac{a_{98}}{6}\nonumber \\ &+\frac{a_{101}}{4}+\frac{a_{102}}{4}+\frac{a_{103}}{4}+90 a_{104}-\frac{1121}{241920}.
 \end{align}
From  $\tau_{3}\Phi(\tau_{3}, \tau_{0})|_{t=0} = 0$, we obtain
\begin{align}
 0 &= \frac{29 a_2}{576}+\frac{29 a_{14}}{576}+\frac{29 a_{23}}{2880}+\frac{a_{32}}{576}+\frac{a_{34}}{576}+\frac{a_{46}}{576}+\frac{a_{61}}{576}+\frac{a_{65}}{576}\nonumber \\ &+\frac{7 a_{84}}{12}+\frac{7 a_{89}}{12}+\frac{7 a_{92}}{12}+\frac{a_{97}}{24}+\frac{a_{99}}{24}+\frac{a_{100}}{24}+\frac{a_{102}}{24}+20 a_{104}\nonumber \\ &-\frac{503}{241920}.
 \end{align}
From  $\tau_{2}\Phi(\tau_{3}, \tau_{1})|_{t=0} = 0$, we obtain
\begin{align}
 0 &= \frac{29 a_2}{288}+\frac{11 a_3}{1440}+\frac{29 a_4}{1440}+\frac{a_{34}}{288}+\frac{a_{37}}{192}+\frac{a_{38}}{576}+\frac{a_{57}}{192}+\frac{a_{58}}{576}\nonumber \\ &+\frac{a_{61}}{288}+\frac{35 a_{84}}{24}+\frac{7 a_{85}}{24}+\frac{35 a_{92}}{24}+\frac{7 a_{93}}{24}+\frac{a_{96}}{8}+\frac{a_{97}}{8}+\frac{a_{102}}{12}\nonumber \\ &+\frac{a_{103}}{6}+60 a_{104}-\frac{17}{5760}.
 \end{align}
From  $\tau_{2}\Phi(\tau_{4}, \tau_{0})|_{t=0} = 0$, we obtain
\begin{align}
 0 &= \frac{11 a_2}{288}+\frac{a_3}{1152}+\frac{a_4}{384}+\frac{a_8}{9216}+\frac{11 a_{14}}{288}+\frac{11 a_{23}}{1440}+\frac{a_{37}}{576}+\frac{a_{43}}{576}\nonumber \\ &+\frac{a_{52}}{576}+\frac{a_{57}}{576}+\frac{11 a_{84}}{24}+\frac{a_{85}}{24}+\frac{11 a_{89}}{24}+\frac{11 a_{92}}{24}+\frac{a_{93}}{24}+\frac{a_{96}}{24}\nonumber \\ &+\frac{a_{103}}{24}+15 a_{104}-\frac{77}{69120}.
 \end{align}
From  $\tau_{3}\tau_{4}\Phi(\tau_{0}, \tau_{0})|_{t=0} = 0$, we obtain
\begin{align}
 0 &= \frac{1121 a_1}{241920}+\frac{17 a_2}{160}+\frac{a_7}{384}+\frac{a_{11}}{9216}+\frac{a_{12}}{9216}+\frac{17 a_{14}}{160}+\frac{17 a_{15}}{160}+\frac{17 a_{16}}{160}\nonumber \\ &+\frac{17 a_{18}}{160}+\frac{17 a_{19}}{160}+\frac{a_{20}}{9216}+\frac{17 a_{23}}{960}+\frac{17 a_{24}}{960}+\frac{17 a_{25}}{960}+\frac{17 a_{29}}{960}+\frac{a_{32}}{576}\nonumber \\ &+\frac{a_{33}}{576}+\frac{a_{34}}{576}+\frac{a_{41}}{576}+\frac{a_{42}}{576}+\frac{a_{46}}{576}+\frac{a_{47}}{576}+\frac{a_{49}}{576}+\frac{a_{51}}{576}\nonumber \\ &+\frac{a_{55}}{576}+\frac{a_{56}}{576}+\frac{a_{61}}{576}+\frac{a_{62}}{576}+\frac{a_{65}}{576}+\frac{25 a_{84}}{24}+\frac{a_{87}}{24}+\frac{25 a_{88}}{24}\nonumber \\ &+\frac{25 a_{89}}{24}+\frac{25 a_{90}}{24}+\frac{25 a_{92}}{24}+\frac{a_{95}}{24}+\frac{a_{97}}{24}+\frac{a_{98}}{24}+\frac{a_{99}}{24}+\frac{a_{100}}{24}\nonumber \\ &+\frac{a_{101}}{24}+\frac{a_{102}}{24}+35 a_{104}-\frac{1121}{34560}.
 \end{align}
From  $\tau_{2}\tau_{5}\Phi(\tau_{0}, \tau_{0})|_{t=0} = 0$, we obtain
\begin{align}
 0 &= \frac{17 a_1}{5760}+\frac{a_2}{15}+\frac{a_3}{288}+\frac{a_4}{288}+\frac{a_5}{288}+\frac{a_6}{288}+\frac{a_8}{6912}+\frac{a_9}{6912}\nonumber \\ &+\frac{a_{10}}{6912}+\frac{a_{14}}{15}+\frac{a_{15}}{15}+\frac{a_{16}}{15}+\frac{a_{17}}{288}+\frac{a_{18}}{15}+\frac{a_{19}}{15}+\frac{a_{21}}{6912}\nonumber \\ &+\frac{a_{22}}{6912}+\frac{a_{23}}{90}+\frac{a_{24}}{90}+\frac{a_{25}}{90}+\frac{a_{26}}{1152}+\frac{a_{27}}{1152}+\frac{a_{28}}{1152}+\frac{a_{29}}{90}\nonumber \\ &+\frac{a_{30}}{1152}+\frac{a_{31}}{27648}+\frac{a_{37}}{576}+\frac{a_{43}}{576}+\frac{a_{44}}{576}+\frac{a_{50}}{576}+\frac{a_{52}}{576}+\frac{a_{53}}{576}\nonumber \\ &+\frac{a_{57}}{576}+\frac{2 a_{84}}{3}+\frac{a_{85}}{24}+\frac{a_{86}}{24}+\frac{2 a_{88}}{3}+\frac{2 a_{89}}{3}+\frac{2 a_{90}}{3}+\frac{a_{91}}{24}\nonumber \\ &+\frac{2 a_{92}}{3}+\frac{a_{93}}{24}+\frac{a_{94}}{24}+\frac{a_{96}}{24}+\frac{a_{103}}{24}+21 a_{104}-\frac{119}{5760}.
 \end{align}
From  $\tau_{3}\tau_{3}\Phi(\tau_{0}, \tau_{1})|_{t=0} = 0$, we obtain
\begin{align}
 0 &= \frac{29 a_2}{96}+\frac{29 a_3}{576}+\frac{29 a_4}{576}+\frac{109 a_{15}}{576}+\frac{29 a_{16}}{96}+\frac{29 a_{17}}{576}+\frac{29 a_{24}}{576}+\frac{29 a_{28}}{2880}\nonumber \\ &+\frac{29 a_{30}}{2880}+\frac{a_{32}}{96}+\frac{a_{33}}{96}+\frac{a_{35}}{288}+\frac{a_{47}}{96}+\frac{a_{48}}{288}+\frac{a_{62}}{96}+\frac{a_{64}}{288}\nonumber \\ &+\frac{a_{65}}{96}+\frac{7 a_{84}}{2}+\frac{7 a_{85}}{12}+\frac{7 a_{90}}{2}+\frac{7 a_{91}}{12}+\frac{7 a_{92}}{2}+\frac{7 a_{93}}{12}+\frac{a_{98}}{4}\nonumber \\ &+\frac{a_{99}}{3}+\frac{a_{100}}{4}+\frac{a_{101}}{3}+140 a_{104}-\frac{205}{3456}.
 \end{align}
From  $\tau_{2}\tau_{4}\Phi(\tau_{0}, \tau_{1})|_{t=0} = 0$, we obtain
\begin{align}
 0 &= \frac{11 a_2}{48}+\frac{11 a_3}{288}+\frac{11 a_4}{288}+\frac{11 a_5}{1440}+\frac{a_6}{96}+\frac{a_7}{128}+\frac{11 a_9}{34560}+\frac{a_{10}}{2304}\nonumber \\ &+\frac{a_{11}}{4608}+\frac{a_{13}}{9216}+\frac{7 a_{15}}{48}+\frac{11 a_{16}}{48}+\frac{11 a_{17}}{288}+\frac{11 a_{24}}{288}+\frac{11 a_{28}}{1440}+\frac{11 a_{30}}{1440}\nonumber \\ &+\frac{a_{34}}{288}+\frac{a_{37}}{144}+\frac{a_{38}}{576}+\frac{a_{42}}{288}+\frac{a_{44}}{144}+\frac{a_{45}}{576}+\frac{a_{53}}{144}+\frac{a_{54}}{576}\nonumber \\ &+\frac{a_{55}}{288}+\frac{a_{57}}{144}+\frac{a_{58}}{576}+\frac{a_{61}}{288}+\frac{11 a_{84}}{4}+\frac{11 a_{85}}{24}+\frac{a_{86}}{6}+\frac{a_{87}}{8}\nonumber \\ &+\frac{11 a_{90}}{4}+\frac{11 a_{91}}{24}+\frac{11 a_{92}}{4}+\frac{11 a_{93}}{24}+\frac{a_{94}}{6}+\frac{a_{95}}{8}+\frac{a_{96}}{6}+\frac{a_{97}}{8}\nonumber \\ &+\frac{a_{102}}{12}+\frac{5 a_{103}}{24}+105 a_{104}-\frac{53}{1152}.
 \end{align}
From  $\tau_{2}\tau_{3}\Phi(\tau_{0}, \tau_{2})|_{t=0} = 0$, we obtain
\begin{align}
 0 &= \frac{5 a_2}{12}+\frac{29 a_5}{2880}+\frac{29 a_6}{1440}+\frac{29 a_7}{1440}+\frac{29 a_9}{69120}+\frac{29 a_{10}}{34560}+\frac{29 a_{11}}{34560}+\frac{109 a_{15}}{576}\nonumber \\ &+\frac{5 a_{16}}{12}+\frac{5 a_{24}}{72}+\frac{a_{32}}{144}+\frac{a_{33}}{144}+\frac{7 a_{34}}{576}+\frac{a_{36}}{576}+\frac{7 a_{37}}{576}+\frac{a_{40}}{576}\nonumber \\ &+\frac{7 a_{42}}{576}+\frac{7 a_{44}}{576}+\frac{a_{47}}{144}+\frac{7 a_{53}}{576}+\frac{7 a_{55}}{576}+\frac{7 a_{57}}{576}+\frac{a_{60}}{576}+\frac{7 a_{61}}{576}\nonumber \\ &+\frac{a_{62}}{144}+\frac{a_{63}}{576}+\frac{a_{65}}{144}+\frac{a_{67}}{576}+\frac{a_{68}}{576}+\frac{a_{69}}{576}+\frac{a_{73}}{13824}+\frac{a_{77}}{13824}\nonumber \\ &+\frac{a_{83}}{13824}+\frac{59 a_{84}}{12}+\frac{7 a_{86}}{24}+\frac{7 a_{87}}{24}+\frac{59 a_{90}}{12}+\frac{59 a_{92}}{12}+\frac{7 a_{94}}{24}+\frac{7 a_{95}}{24}\nonumber \\ &+\frac{7 a_{96}}{24}+\frac{7 a_{97}}{24}+\frac{a_{98}}{6}+\frac{a_{99}}{4}+\frac{a_{100}}{6}+\frac{a_{101}}{4}+\frac{5 a_{102}}{12}+\frac{5 a_{103}}{12}\nonumber \\ &+210 a_{104}-\frac{205}{3456}.
 \end{align}
From  $\tau_{2}\tau_{2}\Phi(\tau_{0}, \tau_{3})|_{t=0} = 0$, we obtain
\begin{align}
 0 &= \frac{5 a_2}{12}+\frac{11 a_5}{720}+\frac{29 a_6}{720}+\frac{11 a_9}{17280}+\frac{29 a_{10}}{17280}+\frac{7 a_{15}}{48}+\frac{5 a_{16}}{12}+\frac{5 a_{24}}{72}\nonumber \\ &+\frac{a_{32}}{144}+\frac{a_{33}}{144}+\frac{7 a_{37}}{288}+\frac{a_{39}}{288}+\frac{7 a_{44}}{288}+\frac{a_{47}}{144}+\frac{7 a_{53}}{288}+\frac{7 a_{57}}{288}\nonumber \\ &+\frac{a_{59}}{288}+\frac{a_{62}}{144}+\frac{a_{65}}{144}+\frac{a_{66}}{288}+\frac{a_{70}}{288}+\frac{a_{75}}{6912}+\frac{a_{79}}{6912}+\frac{59 a_{84}}{12}\nonumber \\ &+\frac{7 a_{86}}{12}+\frac{59 a_{90}}{12}+\frac{59 a_{92}}{12}+\frac{7 a_{94}}{12}+\frac{7 a_{96}}{12}+\frac{a_{98}}{4}+\frac{a_{99}}{6}+\frac{a_{100}}{4}\nonumber \\ &+\frac{a_{101}}{6}+\frac{5 a_{103}}{6}+210 a_{104}-\frac{53}{1152}.
 \end{align}
From  $\tau_{3}\tau_{3}\Phi(\tau_{1}, \tau_{0})|_{t=0} = 0$, we obtain
\begin{align}
 0 &= \frac{29 a_2}{96}+\frac{29 a_5}{576}+\frac{29 a_6}{576}+\frac{29 a_{14}}{96}+\frac{29 a_{17}}{576}+\frac{29 a_{23}}{576}+\frac{29 a_{27}}{2880}+\frac{29 a_{28}}{2880}\nonumber \\ &+\frac{29 a_{30}}{2880}+\frac{a_{32}}{96}+\frac{a_{34}}{96}+\frac{a_{36}}{288}+\frac{a_{46}}{96}+\frac{a_{48}}{288}+\frac{a_{61}}{96}+\frac{a_{63}}{288}\nonumber \\ &+\frac{a_{65}}{96}+\frac{7 a_{84}}{2}+\frac{7 a_{86}}{12}+\frac{7 a_{89}}{2}+\frac{7 a_{91}}{12}+\frac{7 a_{92}}{2}+\frac{7 a_{94}}{12}+\frac{a_{97}}{4}\nonumber \\ &+\frac{a_{99}}{3}+\frac{a_{100}}{4}+\frac{a_{102}}{3}+140 a_{104}-\frac{583}{13824}.
 \end{align}
From  $\tau_{2}\tau_{4}\Phi(\tau_{1}, \tau_{0})|_{t=0} = 0$, we obtain
\begin{align}
 0 &= \frac{11 a_2}{48}+\frac{11 a_3}{1440}+\frac{a_4}{96}+\frac{11 a_5}{288}+\frac{11 a_6}{288}+\frac{a_7}{128}+\frac{a_8}{2304}+\frac{a_{12}}{4608}\nonumber \\ &+\frac{a_{13}}{9216}+\frac{11 a_{14}}{48}+\frac{11 a_{17}}{288}+\frac{11 a_{23}}{288}+\frac{11 a_{27}}{1440}+\frac{11 a_{28}}{1440}+\frac{11 a_{30}}{1440}+\frac{a_{33}}{288}\nonumber \\ &+\frac{a_{37}}{144}+\frac{a_{39}}{576}+\frac{a_{41}}{288}+\frac{a_{43}}{144}+\frac{a_{45}}{576}+\frac{a_{52}}{144}+\frac{a_{54}}{576}+\frac{a_{56}}{288}\nonumber \\ &+\frac{a_{57}}{144}+\frac{a_{59}}{576}+\frac{a_{62}}{288}+\frac{11 a_{84}}{4}+\frac{a_{85}}{6}+\frac{11 a_{86}}{24}+\frac{a_{87}}{8}+\frac{11 a_{89}}{4}\nonumber \\ &+\frac{11 a_{91}}{24}+\frac{11 a_{92}}{4}+\frac{a_{93}}{6}+\frac{11 a_{94}}{24}+\frac{a_{95}}{8}+\frac{a_{96}}{6}+\frac{a_{98}}{8}+\frac{a_{101}}{12}\nonumber \\ &+\frac{5 a_{103}}{24}+105 a_{104}-\frac{1121}{34560}.
 \end{align}
From  $\tau_{2}\tau_{3}\Phi(\tau_{1}, \tau_{1})|_{t=0} = 0$, we obtain
\begin{align}
 0 &= \frac{29 a_2}{48}+\frac{5 a_3}{72}+\frac{29 a_4}{288}+\frac{5 a_5}{72}+\frac{29 a_6}{288}+\frac{29 a_7}{720}+\frac{a_{32}}{96}+\frac{a_{33}}{96}\nonumber \\ &+\frac{a_{34}}{96}+\frac{a_{35}}{288}+\frac{a_{36}}{288}+\frac{a_{37}}{48}+\frac{a_{38}}{192}+\frac{a_{39}}{192}+\frac{a_{40}}{288}+\frac{a_{57}}{48}\nonumber \\ &+\frac{a_{58}}{192}+\frac{a_{59}}{192}+\frac{a_{60}}{288}+\frac{a_{61}}{96}+\frac{a_{62}}{96}+\frac{a_{63}}{288}+\frac{a_{64}}{288}+\frac{a_{65}}{96}\nonumber \\ &+\frac{35 a_{84}}{4}+\frac{35 a_{85}}{24}+\frac{35 a_{86}}{24}+\frac{7 a_{87}}{12}+\frac{35 a_{92}}{4}+\frac{35 a_{93}}{24}+\frac{35 a_{94}}{24}+\frac{7 a_{95}}{12}\nonumber \\ &+\frac{a_{96}}{2}+\frac{3 a_{97}}{8}+\frac{3 a_{98}}{8}+\frac{a_{99}}{2}+\frac{a_{100}}{4}+\frac{a_{101}}{3}+\frac{a_{102}}{3}+\frac{5 a_{103}}{6}\nonumber \\ &+420 a_{104}-\frac{205}{3456}.
 \end{align}
From  $\tau_{2}\tau_{2}\Phi(\tau_{1}, \tau_{2})|_{t=0} = 0$, we obtain
\begin{align}
 0 &= \frac{7 a_2}{8}+\frac{5 a_5}{72}+\frac{7 a_6}{48}+\frac{a_{32}}{72}+\frac{a_{33}}{36}+\frac{a_{34}}{36}+\frac{a_{36}}{144}+\frac{a_{37}}{18}\nonumber \\ &+\frac{a_{39}}{72}+\frac{a_{57}}{18}+\frac{a_{59}}{72}+\frac{a_{61}}{36}+\frac{a_{62}}{36}+\frac{a_{63}}{144}+\frac{a_{65}}{72}+\frac{a_{66}}{144}\nonumber \\ &+\frac{a_{67}}{144}+\frac{a_{68}}{96}+\frac{a_{69}}{144}+\frac{a_{70}}{96}+\frac{a_{79}}{3456}+\frac{a_{81}}{6912}+\frac{a_{83}}{3456}+12 a_{84}\nonumber \\ &+2 a_{86}+12 a_{92}+2 a_{94}+\frac{4 a_{96}}{3}+\frac{2 a_{97}}{3}+a_{98}+\frac{a_{99}}{2}+\frac{a_{100}}{2}\nonumber \\ &+a_{101}+\frac{5 a_{102}}{4}+\frac{5 a_{103}}{2}+630 a_{104}-\frac{205}{3456}.
 \end{align}
From  $\tau_{2}\tau_{3}\Phi(\tau_{2}, \tau_{0})|_{t=0} = 0$, we obtain
\begin{align}
 0 &= \frac{5 a_2}{12}+\frac{29 a_3}{2880}+\frac{29 a_4}{1440}+\frac{29 a_7}{1440}+\frac{29 a_8}{34560}+\frac{29 a_{12}}{34560}+\frac{5 a_{14}}{12}+\frac{5 a_{23}}{72}\nonumber \\ &+\frac{a_{32}}{144}+\frac{7 a_{33}}{576}+\frac{a_{34}}{144}+\frac{a_{35}}{576}+\frac{7 a_{37}}{576}+\frac{a_{40}}{576}+\frac{7 a_{41}}{576}+\frac{7 a_{43}}{576}\nonumber \\ &+\frac{a_{46}}{144}+\frac{7 a_{52}}{576}+\frac{7 a_{56}}{576}+\frac{7 a_{57}}{576}+\frac{a_{60}}{576}+\frac{a_{61}}{144}+\frac{7 a_{62}}{576}+\frac{a_{64}}{576}\nonumber \\ &+\frac{a_{65}}{144}+\frac{a_{67}}{576}+\frac{a_{68}}{576}+\frac{a_{69}}{576}+\frac{a_{72}}{13824}+\frac{a_{78}}{13824}+\frac{a_{83}}{13824}+\frac{59 a_{84}}{12}\nonumber \\ &+\frac{7 a_{85}}{24}+\frac{7 a_{87}}{24}+\frac{59 a_{89}}{12}+\frac{59 a_{92}}{12}+\frac{7 a_{93}}{24}+\frac{7 a_{95}}{24}+\frac{7 a_{96}}{24}+\frac{a_{97}}{6}\nonumber \\ &+\frac{7 a_{98}}{24}+\frac{a_{99}}{4}+\frac{a_{100}}{6}+\frac{5 a_{101}}{12}+\frac{a_{102}}{4}+\frac{5 a_{103}}{12}+210 a_{104}-\frac{1121}{34560}.
 \end{align}
From  $\tau_{2}\tau_{2}\Phi(\tau_{2}, \tau_{1})|_{t=0} = 0$, we obtain
\begin{align}
 0 &= \frac{7 a_2}{8}+\frac{5 a_3}{72}+\frac{7 a_4}{48}+\frac{a_{32}}{72}+\frac{a_{33}}{36}+\frac{a_{34}}{36}+\frac{a_{35}}{144}+\frac{a_{37}}{18}\nonumber \\ &+\frac{a_{38}}{72}+\frac{a_{57}}{18}+\frac{a_{58}}{72}+\frac{a_{61}}{36}+\frac{a_{62}}{36}+\frac{a_{64}}{144}+\frac{a_{65}}{72}+\frac{a_{66}}{144}\nonumber \\ &+\frac{a_{67}}{96}+\frac{a_{68}}{144}+\frac{a_{69}}{144}+\frac{a_{70}}{96}+\frac{a_{79}}{3456}+\frac{a_{80}}{6912}+\frac{a_{83}}{3456}+12 a_{84}\nonumber \\ &+2 a_{85}+12 a_{92}+2 a_{93}+\frac{4 a_{96}}{3}+a_{97}+\frac{2 a_{98}}{3}+\frac{a_{99}}{2}+\frac{a_{100}}{2}\nonumber \\ &+\frac{5 a_{101}}{4}+a_{102}+\frac{5 a_{103}}{2}+630 a_{104}-\frac{53}{1152}.
 \end{align}
From  $\tau_{2}\tau_{2}\Phi(\tau_{3}, \tau_{0})|_{t=0} = 0$, we obtain
\begin{align}
 0 &= \frac{5 a_2}{12}+\frac{11 a_3}{720}+\frac{29 a_4}{720}+\frac{29 a_8}{17280}+\frac{5 a_{14}}{12}+\frac{5 a_{23}}{72}+\frac{a_{32}}{144}+\frac{a_{34}}{144}\nonumber \\ &+\frac{7 a_{37}}{288}+\frac{a_{38}}{288}+\frac{7 a_{43}}{288}+\frac{a_{46}}{144}+\frac{7 a_{52}}{288}+\frac{7 a_{57}}{288}+\frac{a_{58}}{288}+\frac{a_{61}}{144}\nonumber \\ &+\frac{a_{65}}{144}+\frac{a_{66}}{288}+\frac{a_{70}}{288}+\frac{a_{74}}{6912}+\frac{a_{79}}{6912}+\frac{59 a_{84}}{12}+\frac{7 a_{85}}{12}+\frac{59 a_{89}}{12}\nonumber \\ &+\frac{59 a_{92}}{12}+\frac{7 a_{93}}{12}+\frac{7 a_{96}}{12}+\frac{a_{97}}{4}+\frac{a_{99}}{6}+\frac{a_{100}}{4}+\frac{a_{102}}{6}+\frac{5 a_{103}}{6}\nonumber \\ &+210 a_{104}-\frac{119}{5760}.
 \end{align}
From  $\tau_{2}\tau_{3}\tau_{3}\Phi(\tau_{0}, \tau_{0})|_{t=0} = 0$, we obtain
\begin{align}
 0 &= \frac{205 a_1}{3456}+\frac{763 a_2}{576}+\frac{29 a_3}{576}+\frac{29 a_4}{576}+\frac{29 a_5}{576}+\frac{29 a_6}{576}+\frac{29 a_7}{720}+\frac{29 a_8}{13824}\nonumber \\ &+\frac{29 a_9}{13824}+\frac{29 a_{10}}{13824}+\frac{29 a_{11}}{17280}+\frac{29 a_{12}}{17280}+\frac{763 a_{14}}{576}+\frac{763 a_{15}}{576}+\frac{763 a_{16}}{576}+\frac{29 a_{17}}{576}\nonumber \\ &+\frac{763 a_{18}}{576}+\frac{763 a_{19}}{576}+\frac{29 a_{20}}{17280}+\frac{29 a_{21}}{13824}+\frac{29 a_{22}}{13824}+\frac{109 a_{23}}{576}+\frac{109 a_{24}}{576}+\frac{109 a_{25}}{576}\nonumber \\ &+\frac{29 a_{26}}{2880}+\frac{29 a_{27}}{2880}+\frac{29 a_{28}}{2880}+\frac{109 a_{29}}{576}+\frac{29 a_{30}}{2880}+\frac{29 a_{31}}{69120}+\frac{7 a_{32}}{288}+\frac{7 a_{33}}{288}\nonumber \\ &+\frac{7 a_{34}}{288}+\frac{a_{35}}{288}+\frac{a_{36}}{288}+\frac{7 a_{37}}{288}+\frac{a_{40}}{288}+\frac{7 a_{41}}{288}+\frac{7 a_{42}}{288}+\frac{7 a_{43}}{288}\nonumber \\ &+\frac{7 a_{44}}{288}+\frac{7 a_{46}}{288}+\frac{7 a_{47}}{288}+\frac{a_{48}}{288}+\frac{7 a_{49}}{288}+\frac{7 a_{50}}{288}+\frac{7 a_{51}}{288}+\frac{7 a_{52}}{288}\nonumber \\ &+\frac{7 a_{53}}{288}+\frac{7 a_{55}}{288}+\frac{7 a_{56}}{288}+\frac{7 a_{57}}{288}+\frac{a_{60}}{288}+\frac{7 a_{61}}{288}+\frac{7 a_{62}}{288}+\frac{a_{63}}{288}\nonumber \\ &+\frac{a_{64}}{288}+\frac{7 a_{65}}{288}+\frac{a_{67}}{288}+\frac{a_{68}}{288}+\frac{a_{69}}{288}+\frac{a_{71}}{6912}+\frac{a_{72}}{6912}+\frac{a_{73}}{6912}\nonumber \\ &+\frac{a_{77}}{6912}+\frac{a_{78}}{6912}+\frac{a_{83}}{6912}+\frac{40 a_{84}}{3}+\frac{7 a_{85}}{12}+\frac{7 a_{86}}{12}+\frac{7 a_{87}}{12}+\frac{40 a_{88}}{3}\nonumber \\ &+\frac{40 a_{89}}{3}+\frac{40 a_{90}}{3}+\frac{7 a_{91}}{12}+\frac{40 a_{92}}{3}+\frac{7 a_{93}}{12}+\frac{7 a_{94}}{12}+\frac{7 a_{95}}{12}+\frac{7 a_{96}}{12}\nonumber \\ &+\frac{7 a_{97}}{12}+\frac{7 a_{98}}{12}+\frac{5 a_{99}}{6}+\frac{7 a_{100}}{12}+\frac{5 a_{101}}{6}+\frac{5 a_{102}}{6}+\frac{5 a_{103}}{6}+560 a_{104}\nonumber \\ &-\frac{205}{432}.
 \end{align}
From  $\tau_{2}\tau_{2}\tau_{4}\Phi(\tau_{0}, \tau_{0})|_{t=0} = 0$, we obtain
\begin{align}
 0 &= \frac{53 a_1}{1152}+\frac{49 a_2}{48}+\frac{11 a_3}{144}+\frac{11 a_4}{144}+\frac{11 a_5}{144}+\frac{11 a_6}{144}+\frac{a_7}{64}+\frac{11 a_8}{3456}\nonumber \\ &+\frac{11 a_9}{3456}+\frac{11 a_{10}}{3456}+\frac{a_{11}}{2304}+\frac{a_{12}}{2304}+\frac{a_{13}}{4608}+\frac{49 a_{14}}{48}+\frac{49 a_{15}}{48}+\frac{49 a_{16}}{48}\nonumber \\ &+\frac{11 a_{17}}{144}+\frac{49 a_{18}}{48}+\frac{49 a_{19}}{48}+\frac{a_{20}}{2304}+\frac{11 a_{21}}{3456}+\frac{11 a_{22}}{3456}+\frac{7 a_{23}}{48}+\frac{7 a_{24}}{48}\nonumber \\ &+\frac{7 a_{25}}{48}+\frac{11 a_{26}}{720}+\frac{11 a_{27}}{720}+\frac{11 a_{28}}{720}+\frac{7 a_{29}}{48}+\frac{11 a_{30}}{720}+\frac{11 a_{31}}{17280}+\frac{a_{32}}{144}\nonumber \\ &+\frac{a_{33}}{144}+\frac{a_{34}}{144}+\frac{11 a_{37}}{288}+\frac{a_{38}}{288}+\frac{a_{39}}{288}+\frac{a_{41}}{144}+\frac{a_{42}}{144}+\frac{11 a_{43}}{288}\nonumber \\ &+\frac{11 a_{44}}{288}+\frac{a_{45}}{288}+\frac{a_{46}}{144}+\frac{a_{47}}{144}+\frac{a_{49}}{144}+\frac{11 a_{50}}{288}+\frac{a_{51}}{144}+\frac{11 a_{52}}{288}\nonumber \\ &+\frac{11 a_{53}}{288}+\frac{a_{54}}{288}+\frac{a_{55}}{144}+\frac{a_{56}}{144}+\frac{11 a_{57}}{288}+\frac{a_{58}}{288}+\frac{a_{59}}{288}+\frac{a_{61}}{144}\nonumber \\ &+\frac{a_{62}}{144}+\frac{a_{65}}{144}+\frac{a_{66}}{288}+\frac{a_{70}}{288}+\frac{a_{74}}{6912}+\frac{a_{75}}{6912}+\frac{a_{79}}{6912}+\frac{125 a_{84}}{12}\nonumber \\ &+\frac{11 a_{85}}{12}+\frac{11 a_{86}}{12}+\frac{a_{87}}{4}+\frac{125 a_{88}}{12}+\frac{125 a_{89}}{12}+\frac{125 a_{90}}{12}+\frac{11 a_{91}}{12}+\frac{125 a_{92}}{12}\nonumber \\ &+\frac{11 a_{93}}{12}+\frac{11 a_{94}}{12}+\frac{a_{95}}{4}+\frac{11 a_{96}}{12}+\frac{a_{97}}{4}+\frac{a_{98}}{4}+\frac{a_{99}}{6}+\frac{a_{100}}{4}\nonumber \\ &+\frac{a_{101}}{6}+\frac{a_{102}}{6}+\frac{5 a_{103}}{4}+420 a_{104}+\frac{a_{105}}{6912}-\frac{53}{144}.
 \end{align}
From  $\tau_{2}\tau_{2}\tau_{3}\Phi(\tau_{0}, \tau_{1})|_{t=0} = 0$, we obtain
\begin{align}
 0 &= \frac{35 a_2}{12}+\frac{5 a_3}{12}+\frac{5 a_4}{12}+\frac{5 a_5}{36}+\frac{29 a_6}{144}+\frac{29 a_7}{240}+\frac{5 a_9}{864}+\frac{29 a_{10}}{3456}\nonumber \\ &+\frac{29 a_{11}}{8640}+\frac{29 a_{13}}{17280}+\frac{11 a_{15}}{6}+\frac{35 a_{16}}{12}+\frac{5 a_{17}}{12}+\frac{5 a_{24}}{12}+\frac{5 a_{28}}{72}+\frac{5 a_{30}}{72}\nonumber \\ &+\frac{7 a_{32}}{144}+\frac{7 a_{33}}{144}+\frac{7 a_{34}}{144}+\frac{a_{35}}{72}+\frac{a_{36}}{144}+\frac{35 a_{37}}{288}+\frac{7 a_{38}}{288}+\frac{a_{39}}{96}\nonumber \\ &+\frac{a_{40}}{96}+\frac{7 a_{42}}{144}+\frac{35 a_{44}}{288}+\frac{7 a_{45}}{288}+\frac{7 a_{47}}{144}+\frac{a_{48}}{72}+\frac{35 a_{53}}{288}+\frac{7 a_{54}}{288}\nonumber \\ &+\frac{7 a_{55}}{144}+\frac{35 a_{57}}{288}+\frac{7 a_{58}}{288}+\frac{a_{59}}{96}+\frac{a_{60}}{96}+\frac{7 a_{61}}{144}+\frac{7 a_{62}}{144}+\frac{a_{63}}{144}\nonumber \\ &+\frac{a_{64}}{72}+\frac{7 a_{65}}{144}+\frac{a_{66}}{96}+\frac{a_{67}}{96}+\frac{a_{68}}{144}+\frac{a_{69}}{144}+\frac{a_{70}}{72}+\frac{a_{73}}{3456}\nonumber \\ &+\frac{a_{75}}{2304}+\frac{a_{76}}{6912}+\frac{a_{77}}{3456}+\frac{a_{79}}{2304}+\frac{a_{80}}{6912}+\frac{a_{83}}{3456}+\frac{413 a_{84}}{12}+\frac{59 a_{85}}{12}\nonumber \\ &+\frac{35 a_{86}}{12}+\frac{7 a_{87}}{4}+\frac{413 a_{90}}{12}+\frac{59 a_{91}}{12}+\frac{413 a_{92}}{12}+\frac{59 a_{93}}{12}+\frac{35 a_{94}}{12}+\frac{7 a_{95}}{4}\nonumber \\ &+\frac{35 a_{96}}{12}+\frac{7 a_{97}}{4}+\frac{17 a_{98}}{12}+\frac{23 a_{99}}{12}+\frac{17 a_{100}}{12}+\frac{23 a_{101}}{12}+\frac{5 a_{102}}{3}+5 a_{103}\nonumber \\ &+1680 a_{104}-\frac{193}{288}.
 \end{align}
From  $\tau_{2}\tau_{2}\tau_{2}\Phi(\tau_{0}, \tau_{2})|_{t=0} = 0$, we obtain
\begin{align}
 0 &= \frac{49 a_2}{12}+\frac{5 a_5}{24}+\frac{7 a_6}{16}+\frac{7 a_7}{48}+\frac{5 a_9}{576}+\frac{7 a_{10}}{384}+\frac{7 a_{11}}{1152}+\frac{11 a_{15}}{6}\nonumber \\ &+\frac{49 a_{16}}{12}+\frac{7 a_{24}}{12}+\frac{a_{32}}{12}+\frac{a_{33}}{12}+\frac{a_{34}}{12}+\frac{a_{36}}{48}+\frac{a_{37}}{4}+\frac{a_{39}}{24}\nonumber \\ &+\frac{a_{40}}{48}+\frac{a_{42}}{12}+\frac{a_{44}}{4}+\frac{a_{47}}{12}+\frac{a_{53}}{4}+\frac{a_{55}}{12}+\frac{a_{57}}{4}+\frac{a_{59}}{24}\nonumber \\ &+\frac{a_{60}}{48}+\frac{a_{61}}{12}+\frac{a_{62}}{12}+\frac{a_{63}}{48}+\frac{a_{65}}{12}+\frac{a_{66}}{24}+\frac{a_{67}}{48}+\frac{a_{68}}{32}\nonumber \\ &+\frac{a_{69}}{48}+\frac{a_{70}}{16}+\frac{a_{73}}{1152}+\frac{a_{75}}{576}+\frac{a_{77}}{1152}+\frac{a_{79}}{576}+\frac{a_{81}}{2304}+\frac{a_{82}}{2304}\nonumber \\ &+\frac{a_{83}}{1152}+48 a_{84}+6 a_{86}+2 a_{87}+48 a_{90}+48 a_{92}+6 a_{94}+2 a_{95}\nonumber \\ &+6 a_{96}+2 a_{97}+3 a_{98}+3 a_{99}+3 a_{100}+3 a_{101}+\frac{15 a_{102}}{4}+\frac{45 a_{103}}{4}\nonumber \\ &+2520 a_{104}-\frac{193}{288}.
 \end{align}
From  $\tau_{2}\tau_{2}\tau_{2}\Phi(\tau_{1}, \tau_{1})|_{t=0} = 0$, we obtain
\begin{align}
 0 &= \frac{49 a_2}{8}+\frac{7 a_3}{12}+\frac{7 a_4}{8}+\frac{7 a_5}{12}+\frac{7 a_6}{8}+\frac{7 a_7}{24}+\frac{a_{32}}{8}+\frac{a_{33}}{6}\nonumber \\ &+\frac{a_{34}}{6}+\frac{a_{35}}{24}+\frac{a_{36}}{24}+\frac{5 a_{37}}{12}+\frac{a_{38}}{12}+\frac{a_{39}}{12}+\frac{a_{40}}{24}+\frac{5 a_{57}}{12}\nonumber \\ &+\frac{a_{58}}{12}+\frac{a_{59}}{12}+\frac{a_{60}}{24}+\frac{a_{61}}{6}+\frac{a_{62}}{6}+\frac{a_{63}}{24}+\frac{a_{64}}{24}+\frac{a_{65}}{8}\nonumber \\ &+\frac{a_{66}}{16}+\frac{a_{67}}{16}+\frac{a_{68}}{16}+\frac{a_{69}}{24}+\frac{a_{70}}{8}+\frac{a_{79}}{384}+\frac{a_{80}}{1152}+\frac{a_{81}}{1152}\nonumber \\ &+\frac{a_{82}}{1152}+\frac{a_{83}}{576}+84 a_{84}+12 a_{85}+12 a_{86}+4 a_{87}+84 a_{92}+12 a_{93}\nonumber \\ &+12 a_{94}+4 a_{95}+10 a_{96}+6 a_{97}+6 a_{98}+6 a_{99}+\frac{9 a_{100}}{2}+\frac{15 a_{101}}{2}\nonumber \\ &+\frac{15 a_{102}}{2}+\frac{45 a_{103}}{2}+5040 a_{104}-\frac{193}{288}.
 \end{align}
From  $\tau_{2}\tau_{2}\tau_{2}\Phi(\tau_{2}, \tau_{0})|_{t=0} = 0$, we obtain
\begin{align}
 0 &= \frac{49 a_2}{12}+\frac{5 a_3}{24}+\frac{7 a_4}{16}+\frac{7 a_7}{48}+\frac{7 a_8}{384}+\frac{7 a_{12}}{1152}+\frac{49 a_{14}}{12}+\frac{7 a_{23}}{12}\nonumber \\ &+\frac{a_{32}}{12}+\frac{a_{33}}{12}+\frac{a_{34}}{12}+\frac{a_{35}}{48}+\frac{a_{37}}{4}+\frac{a_{38}}{24}+\frac{a_{40}}{48}+\frac{a_{41}}{12}\nonumber \\ &+\frac{a_{43}}{4}+\frac{a_{46}}{12}+\frac{a_{52}}{4}+\frac{a_{56}}{12}+\frac{a_{57}}{4}+\frac{a_{58}}{24}+\frac{a_{60}}{48}+\frac{a_{61}}{12}\nonumber \\ &+\frac{a_{62}}{12}+\frac{a_{64}}{48}+\frac{a_{65}}{12}+\frac{a_{66}}{24}+\frac{a_{67}}{32}+\frac{a_{68}}{48}+\frac{a_{69}}{48}+\frac{a_{70}}{16}\nonumber \\ &+\frac{a_{72}}{1152}+\frac{a_{74}}{576}+\frac{a_{78}}{1152}+\frac{a_{79}}{576}+\frac{a_{80}}{2304}+\frac{a_{82}}{2304}+\frac{a_{83}}{1152}+48 a_{84}\nonumber \\ &+6 a_{85}+2 a_{87}+48 a_{89}+48 a_{92}+6 a_{93}+2 a_{95}+6 a_{96}+3 a_{97}\nonumber \\ &+2 a_{98}+3 a_{99}+3 a_{100}+\frac{15 a_{101}}{4}+3 a_{102}+\frac{45 a_{103}}{4}+2520 a_{104}-\frac{53}{144}.
 \end{align}
From  $\tau_{2}\tau_{2}\tau_{2}\tau_{3}\Phi(\tau_{0}, \tau_{0})|_{t=0} = 0$, we obtain
\begin{align}
 0 &= \frac{193 a_1}{288}+\frac{44 a_2}{3}+\frac{5 a_3}{4}+\frac{5 a_4}{4}+\frac{5 a_5}{4}+\frac{5 a_6}{4}+\frac{61 a_7}{120}+\frac{5 a_8}{96}\nonumber \\ &+\frac{5 a_9}{96}+\frac{5 a_{10}}{96}+\frac{31 a_{11}}{1920}+\frac{31 a_{12}}{1920}+\frac{29 a_{13}}{5760}+\frac{44 a_{14}}{3}+\frac{44 a_{15}}{3}+\frac{44 a_{16}}{3}\nonumber \\ &+\frac{5 a_{17}}{4}+\frac{44 a_{18}}{3}+\frac{44 a_{19}}{3}+\frac{31 a_{20}}{1920}+\frac{5 a_{21}}{96}+\frac{5 a_{22}}{96}+\frac{11 a_{23}}{6}+\frac{11 a_{24}}{6}\nonumber \\ &+\frac{11 a_{25}}{6}+\frac{5 a_{26}}{24}+\frac{5 a_{27}}{24}+\frac{5 a_{28}}{24}+\frac{11 a_{29}}{6}+\frac{5 a_{30}}{24}+\frac{5 a_{31}}{576}+\frac{11 a_{32}}{48}\nonumber \\ &+\frac{11 a_{33}}{48}+\frac{11 a_{34}}{48}+\frac{a_{35}}{24}+\frac{a_{36}}{24}+\frac{59 a_{37}}{96}+\frac{7 a_{38}}{96}+\frac{7 a_{39}}{96}+\frac{5 a_{40}}{96}\nonumber \\ &+\frac{11 a_{41}}{48}+\frac{11 a_{42}}{48}+\frac{59 a_{43}}{96}+\frac{59 a_{44}}{96}+\frac{7 a_{45}}{96}+\frac{11 a_{46}}{48}+\frac{11 a_{47}}{48}+\frac{a_{48}}{24}\nonumber \\ &+\frac{11 a_{49}}{48}+\frac{59 a_{50}}{96}+\frac{11 a_{51}}{48}+\frac{59 a_{52}}{96}+\frac{59 a_{53}}{96}+\frac{7 a_{54}}{96}+\frac{11 a_{55}}{48}+\frac{11 a_{56}}{48}\nonumber \\ &+\frac{59 a_{57}}{96}+\frac{7 a_{58}}{96}+\frac{7 a_{59}}{96}+\frac{5 a_{60}}{96}+\frac{11 a_{61}}{48}+\frac{11 a_{62}}{48}+\frac{a_{63}}{24}+\frac{a_{64}}{24}\nonumber \\ &+\frac{11 a_{65}}{48}+\frac{7 a_{66}}{96}+\frac{5 a_{67}}{96}+\frac{5 a_{68}}{96}+\frac{a_{69}}{24}+\frac{5 a_{70}}{48}+\frac{a_{71}}{576}+\frac{a_{72}}{576}\nonumber \\ &+\frac{a_{73}}{576}+\frac{7 a_{74}}{2304}+\frac{7 a_{75}}{2304}+\frac{a_{76}}{2304}+\frac{a_{77}}{576}+\frac{a_{78}}{576}+\frac{7 a_{79}}{2304}+\frac{a_{80}}{2304}\nonumber \\ &+\frac{a_{81}}{2304}+\frac{a_{82}}{2304}+\frac{a_{83}}{576}+\frac{605 a_{84}}{4}+\frac{59 a_{85}}{4}+\frac{59 a_{86}}{4}+\frac{29 a_{87}}{4}+\frac{605 a_{88}}{4}\nonumber \\ &+\frac{605 a_{89}}{4}+\frac{605 a_{90}}{4}+\frac{59 a_{91}}{4}+\frac{605 a_{92}}{4}+\frac{59 a_{93}}{4}+\frac{59 a_{94}}{4}+\frac{29 a_{95}}{4}+\frac{59 a_{96}}{4}\nonumber \\ &+\frac{29 a_{97}}{4}+\frac{29 a_{98}}{4}+\frac{35 a_{99}}{4}+\frac{29 a_{100}}{4}+\frac{35 a_{101}}{4}+\frac{35 a_{102}}{4}+\frac{105 a_{103}}{4}+7560 a_{104}\nonumber \\ &+\frac{7 a_{105}}{2304}-\frac{193}{32}.
 \end{align}
From  $\tau_{2}\tau_{2}\tau_{2}\tau_{2}\Phi(\tau_{0}, \tau_{1})|_{t=0} = 0$, we obtain
\begin{align}
 0 &= \frac{98 a_2}{3}+\frac{49 a_3}{12}+\frac{49 a_4}{12}+\frac{7 a_5}{3}+\frac{7 a_6}{2}+\frac{7 a_7}{4}+\frac{7 a_9}{72}+\frac{7 a_{10}}{48}\nonumber \\ &+\frac{7 a_{11}}{144}+\frac{7 a_{13}}{288}+\frac{245 a_{15}}{12}+\frac{98 a_{16}}{3}+\frac{49 a_{17}}{12}+\frac{49 a_{24}}{12}+\frac{7 a_{28}}{12}+\frac{7 a_{30}}{12}\nonumber \\ &+\frac{2 a_{32}}{3}+\frac{2 a_{33}}{3}+\frac{2 a_{34}}{3}+\frac{a_{35}}{6}+\frac{a_{36}}{6}+2 a_{37}+\frac{a_{38}}{3}+\frac{a_{39}}{3}\nonumber \\ &+\frac{a_{40}}{4}+\frac{2 a_{42}}{3}+2 a_{44}+\frac{a_{45}}{3}+\frac{2 a_{47}}{3}+\frac{a_{48}}{6}+2 a_{53}+\frac{a_{54}}{3}\nonumber \\ &+\frac{2 a_{55}}{3}+2 a_{57}+\frac{a_{58}}{3}+\frac{a_{59}}{3}+\frac{a_{60}}{4}+\frac{2 a_{61}}{3}+\frac{2 a_{62}}{3}+\frac{a_{63}}{6}\nonumber \\ &+\frac{a_{64}}{6}+\frac{2 a_{65}}{3}+\frac{a_{66}}{3}+\frac{a_{67}}{4}+\frac{a_{68}}{4}+\frac{a_{69}}{6}+\frac{5 a_{70}}{8}+\frac{a_{73}}{144}\nonumber \\ &+\frac{a_{75}}{72}+\frac{a_{76}}{288}+\frac{a_{77}}{144}+\frac{a_{79}}{72}+\frac{a_{80}}{288}+\frac{a_{81}}{288}+\frac{a_{82}}{192}+\frac{a_{83}}{144}\nonumber \\ &+384 a_{84}+48 a_{85}+48 a_{86}+24 a_{87}+384 a_{90}+48 a_{91}+384 a_{92}+48 a_{93}\nonumber \\ &+48 a_{94}+24 a_{95}+48 a_{96}+24 a_{97}+24 a_{98}+30 a_{99}+24 a_{100}+30 a_{101}\nonumber \\ &+30 a_{102}+105 a_{103}+22680 a_{104}-\frac{1225}{144}.
 \end{align}
From  $\tau_{2}\tau_{2}\tau_{2}\tau_{2}\tau_{2}\Phi(\tau_{0}, \tau_{0})|_{t=0} = 0$, we obtain
\begin{align}
 0 &= \frac{1225 a_1}{144}+\frac{735 a_2}{4}+\frac{245 a_3}{12}+\frac{245 a_4}{12}+\frac{245 a_5}{12}+\frac{245 a_6}{12}+\frac{35 a_7}{4}+\frac{245 a_8}{288}\nonumber \\ &+\frac{245 a_9}{288}+\frac{245 a_{10}}{288}+\frac{35 a_{11}}{144}+\frac{35 a_{12}}{144}+\frac{35 a_{13}}{288}+\frac{735 a_{14}}{4}+\frac{735 a_{15}}{4}+\frac{735 a_{16}}{4}\nonumber \\ &+\frac{245 a_{17}}{12}+\frac{735 a_{18}}{4}+\frac{735 a_{19}}{4}+\frac{35 a_{20}}{144}+\frac{245 a_{21}}{288}+\frac{245 a_{22}}{288}+\frac{245 a_{23}}{12}+\frac{245 a_{24}}{12}\nonumber \\ &+\frac{245 a_{25}}{12}+\frac{35 a_{26}}{12}+\frac{35 a_{27}}{12}+\frac{35 a_{28}}{12}+\frac{245 a_{29}}{12}+\frac{35 a_{30}}{12}+\frac{35 a_{31}}{288}+\frac{10 a_{32}}{3}\nonumber \\ &+\frac{10 a_{33}}{3}+\frac{10 a_{34}}{3}+\frac{5 a_{35}}{6}+\frac{5 a_{36}}{6}+10 a_{37}+\frac{5 a_{38}}{3}+\frac{5 a_{39}}{3}+\frac{5 a_{40}}{4}\nonumber \\ &+\frac{10 a_{41}}{3}+\frac{10 a_{42}}{3}+10 a_{43}+10 a_{44}+\frac{5 a_{45}}{3}+\frac{10 a_{46}}{3}+\frac{10 a_{47}}{3}+\frac{5 a_{48}}{6}\nonumber \\ &+\frac{10 a_{49}}{3}+10 a_{50}+\frac{10 a_{51}}{3}+10 a_{52}+10 a_{53}+\frac{5 a_{54}}{3}+\frac{10 a_{55}}{3}+\frac{10 a_{56}}{3}\nonumber \\ &+10 a_{57}+\frac{5 a_{58}}{3}+\frac{5 a_{59}}{3}+\frac{5 a_{60}}{4}+\frac{10 a_{61}}{3}+\frac{10 a_{62}}{3}+\frac{5 a_{63}}{6}+\frac{5 a_{64}}{6}\nonumber \\ &+\frac{10 a_{65}}{3}+\frac{5 a_{66}}{3}+\frac{5 a_{67}}{4}+\frac{5 a_{68}}{4}+\frac{5 a_{69}}{6}+\frac{25 a_{70}}{8}+\frac{5 a_{71}}{144}+\frac{5 a_{72}}{144}\nonumber \\ &+\frac{5 a_{73}}{144}+\frac{5 a_{74}}{72}+\frac{5 a_{75}}{72}+\frac{5 a_{76}}{288}+\frac{5 a_{77}}{144}+\frac{5 a_{78}}{144}+\frac{5 a_{79}}{72}+\frac{5 a_{80}}{288}\nonumber \\ &+\frac{5 a_{81}}{288}+\frac{5 a_{82}}{192}+\frac{5 a_{83}}{144}+1920 a_{84}+240 a_{85}+240 a_{86}+120 a_{87}+1920 a_{88}\nonumber \\ &+1920 a_{89}+1920 a_{90}+240 a_{91}+1920 a_{92}+240 a_{93}+240 a_{94}+120 a_{95}+240 a_{96}\nonumber \\ &+120 a_{97}+120 a_{98}+150 a_{99}+120 a_{100}+150 a_{101}+150 a_{102}+525 a_{103}+113400 a_{104}\nonumber \\ &+\frac{5 a_{105}}{72}-\frac{6125}{72}.
\label{eqn:ptend}
 \end{align}

\subsection{Relations from the Gromov-Witten invariants of a $\mathbb{C}P^1$}
\label{sec:P1}

From the degree $0$ part of $\Phi(\tau_{0,0}, \tau_{2,1})|_{t=0} = 0$, we obtain
\begin{align}
 0 &= -\frac{7 a_9}{138240}-\frac{7 a_{10}}{46080}-\frac{a_{75}}{13824}-\frac{a_{79}}{13824}-\frac{a_{82}}{13824}+\frac{31}{96768}.
\label{eqn:P1start}
 \end{align}
From the degree $0$ part of $\Phi(\tau_{0,0}, \tau_{3,0})|_{t=0} = 0$, we obtain
\begin{align}
 0 &= \frac{7 a_2}{720}+\frac{a_9}{5760}+\frac{a_{10}}{1920}+\frac{7 a_{15}}{2880}+\frac{7 a_{16}}{720}+\frac{7 a_{24}}{2880}+\frac{a_{37}}{288}+\frac{a_{44}}{288}\nonumber \\ &+\frac{a_{53}}{288}+\frac{a_{57}}{288}+\frac{a_{70}}{288}+\frac{a_{75}}{6912}+\frac{a_{79}}{6912}-\frac{2329}{1451520}.
 \end{align}
From the degree $0$ part of $\Phi(\tau_{1,0}, \tau_{1,1})|_{t=0} = 0$, we obtain
\begin{align}
 0 &= -\frac{a_{79}}{6912}-\frac{a_{80}}{13824}-\frac{a_{81}}{13824}-\frac{a_{82}}{6912}-\frac{a_{83}}{13824}+\frac{31}{96768}.
 \end{align}
From the degree $0$ part of $\Phi(\tau_{1,0}, \tau_{2,0})|_{t=0} = 0$, we obtain
\begin{align}
 0 &= \frac{7 a_2}{240}+\frac{7 a_5}{2880}+\frac{7 a_6}{960}+\frac{a_{33}}{288}+\frac{a_{37}}{96}+\frac{a_{39}}{288}+\frac{a_{57}}{96}+\frac{a_{59}}{288}\nonumber \\ &+\frac{a_{62}}{288}+\frac{a_{67}}{288}+\frac{a_{70}}{96}+\frac{a_{79}}{3456}+\frac{a_{81}}{6912}+\frac{a_{83}}{6912}-\frac{1501}{725760}.
 \end{align}
From the degree $0$ part of $\Phi(\tau_{1,1}, \tau_{1,0})|_{t=0} = 0$, we obtain
\begin{align}
 0 &= -\frac{a_{79}}{6912}-\frac{a_{80}}{13824}-\frac{a_{81}}{13824}-\frac{a_{82}}{6912}-\frac{a_{83}}{13824}+\frac{31}{96768}.
 \end{align}
From the degree $0$ part of $\Phi(\tau_{2,0}, \tau_{0,1})|_{t=0} = 0$, we obtain
\begin{align}
 0 &= -\frac{7 a_8}{46080}-\frac{a_{74}}{13824}-\frac{a_{79}}{13824}-\frac{a_{82}}{13824}+\frac{31}{193536}.
 \end{align}
From the degree $0$ part of $\Phi(\tau_{3,0},\tau_{0,0})|_{t=0} = 0$, we obtain
\begin{align}
 0 &= \frac{7 a_2}{720}+\frac{a_8}{1920}+\frac{7 a_{14}}{720}+\frac{7 a_{23}}{2880}+\frac{a_{37}}{288}+\frac{a_{43}}{288}+\frac{a_{52}}{288}+\frac{a_{57}}{288}\nonumber \\ &+\frac{a_{70}}{288}+\frac{a_{74}}{6912}+\frac{a_{79}}{6912}-\frac{205}{290304}.
 \end{align}
From the degree $0$ part of $\tau_{3,1}\Phi(\tau_{0,0}, \tau_{0,0})|_{t=0} = 0$, we obtain
\begin{align}
 0 &= -\frac{31 a_1}{96768}-\frac{7 a_8}{34560}-\frac{7 a_9}{34560}-\frac{7 a_{10}}{34560}-\frac{7 a_{21}}{34560}-\frac{7 a_{22}}{34560}-\frac{7 a_{31}}{138240}\nonumber \\ &-\frac{a_{74}}{13824}-\frac{a_{75}}{13824}-\frac{a_{79}}{13824}-\frac{a_{82}}{13824}-\frac{a_{105}}{13824}+\frac{31}{16128}.
 \end{align}
From the degree $0$ part of $\tau_{4,0}\Phi(\tau_{0,0}, \tau_{0,0})|_{t=0} = 0$, we obtain
\begin{align}
 0 &= \frac{2329 a_1}{1451520}+\frac{7 a_2}{576}+\frac{a_8}{1440}+\frac{a_9}{1440}+\frac{a_{10}}{1440}+\frac{7 a_{14}}{576}+\frac{7 a_{15}}{576}+\frac{7 a_{16}}{576}\nonumber \\ &+\frac{7 a_{18}}{576}+\frac{7 a_{19}}{576}+\frac{a_{21}}{1440}+\frac{a_{22}}{1440}+\frac{7 a_{23}}{2880}+\frac{7 a_{24}}{2880}+\frac{7 a_{25}}{2880}+\frac{7 a_{29}}{2880}\nonumber \\ &+\frac{a_{31}}{5760}+\frac{a_{37}}{288}+\frac{a_{43}}{288}+\frac{a_{44}}{288}+\frac{a_{50}}{288}+\frac{a_{52}}{288}+\frac{a_{53}}{288}+\frac{a_{57}}{288}\nonumber \\ &+\frac{a_{70}}{288}+\frac{a_{74}}{6912}+\frac{a_{75}}{6912}+\frac{a_{79}}{6912}+\frac{a_{105}}{6912}-\frac{2329}{241920}.
 \end{align}
From the degree $0$ part of $\tau_{2,1}\Phi(\tau_{0,0}, \tau_{1,0})|_{t=0} = 0$, we obtain
\begin{align}
 0 &= -\frac{7 a_9}{23040}-\frac{7 a_{10}}{11520}-\frac{7 a_{11}}{46080}-\frac{7 a_{13}}{46080}-\frac{a_{73}}{13824}-\frac{a_{75}}{4608}-\frac{a_{76}}{13824}\nonumber \\ &-\frac{a_{77}}{13824}-\frac{a_{79}}{4608}-\frac{a_{80}}{13824}-\frac{a_{81}}{13824}-\frac{a_{82}}{4608}-\frac{a_{83}}{13824}+\frac{31}{10752}.
 \end{align}
From the degree $0$ part of $\tau_{3,0}\Phi(\tau_{0,0}, \tau_{1,0})|_{t=0} = 0$, we obtain
\begin{align}
 0 &= \frac{7 a_2}{144}+\frac{7 a_3}{720}+\frac{7 a_4}{720}+\frac{43 a_9}{34560}+\frac{a_{10}}{480}+\frac{a_{11}}{1920}+\frac{a_{13}}{1920}+\frac{7 a_{15}}{288}\nonumber \\ &+\frac{7 a_{16}}{144}+\frac{7 a_{17}}{720}+\frac{7 a_{24}}{720}+\frac{7 a_{28}}{2880}+\frac{7 a_{30}}{2880}+\frac{a_{34}}{288}+\frac{a_{37}}{72}+\frac{a_{38}}{288}\nonumber \\ &+\frac{a_{42}}{288}+\frac{a_{44}}{72}+\frac{a_{45}}{288}+\frac{a_{53}}{72}+\frac{a_{54}}{288}+\frac{a_{55}}{288}+\frac{a_{57}}{72}+\frac{a_{58}}{288}\nonumber \\ &+\frac{a_{61}}{288}+\frac{a_{68}}{288}+\frac{a_{70}}{72}+\frac{a_{73}}{6912}+\frac{a_{75}}{2304}+\frac{a_{76}}{6912}+\frac{a_{77}}{6912}+\frac{a_{79}}{2304}\nonumber \\ &+\frac{a_{80}}{6912}+\frac{a_{83}}{6912}-\frac{859}{48384}.
 \end{align}
From the degree $0$ part of $\tau_{2,0}\Phi(\tau_{0,0}, \tau_{2,0})|_{t=0} = 0$, we obtain
\begin{align}
 0 &= \frac{7 a_2}{96}+\frac{7 a_5}{2880}+\frac{7 a_6}{960}+\frac{7 a_7}{960}+\frac{31 a_9}{23040}+\frac{221 a_{10}}{69120}+\frac{7 a_{11}}{23040}+\frac{7 a_{15}}{288}\nonumber \\ &+\frac{7 a_{16}}{96}+\frac{7 a_{24}}{480}+\frac{a_{32}}{288}+\frac{a_{33}}{288}+\frac{a_{37}}{48}+\frac{a_{39}}{288}+\frac{a_{40}}{288}+\frac{a_{44}}{48}\nonumber \\ &+\frac{a_{47}}{288}+\frac{a_{53}}{48}+\frac{a_{57}}{48}+\frac{a_{59}}{288}+\frac{a_{60}}{288}+\frac{a_{62}}{288}+\frac{a_{65}}{288}+\frac{a_{66}}{288}\nonumber \\ &+\frac{a_{67}}{288}+\frac{a_{70}}{48}+\frac{a_{73}}{6912}+\frac{a_{75}}{1152}+\frac{a_{77}}{6912}+\frac{a_{79}}{1152}+\frac{a_{81}}{6912}+\frac{a_{82}}{2304}\nonumber \\ &+\frac{a_{83}}{6912}-\frac{859}{48384}.
 \end{align}
From the degree $0$ part of $\tau_{2,1}\Phi(\tau_{1,0}, \tau_{0,0})|_{t=0} = 0$, we obtain
\begin{align}
 0 &= -\frac{7 a_8}{11520}-\frac{7 a_{12}}{46080}-\frac{7 a_{13}}{46080}-\frac{a_{72}}{13824}-\frac{a_{74}}{4608}-\frac{a_{76}}{13824}-\frac{a_{78}}{13824}\nonumber \\ &-\frac{a_{79}}{4608}-\frac{a_{80}}{13824}-\frac{a_{81}}{13824}-\frac{a_{82}}{4608}-\frac{a_{83}}{13824}+\frac{31}{16128}.
 \end{align}
From the degree $0$ part of $\tau_{2,0}\Phi(\tau_{1,0}, \tau_{1,0})|_{t=0} = 0$, we obtain
\begin{align}
 0 &= \frac{7 a_2}{48}+\frac{7 a_3}{480}+\frac{7 a_4}{240}+\frac{7 a_5}{480}+\frac{7 a_6}{240}+\frac{7 a_7}{480}+\frac{a_{32}}{144}+\frac{a_{33}}{96}\nonumber \\ &+\frac{a_{34}}{96}+\frac{a_{35}}{288}+\frac{a_{36}}{288}+\frac{a_{37}}{24}+\frac{a_{38}}{96}+\frac{a_{39}}{96}+\frac{a_{40}}{144}+\frac{a_{57}}{24}\nonumber \\ &+\frac{a_{58}}{96}+\frac{a_{59}}{96}+\frac{a_{60}}{144}+\frac{a_{61}}{96}+\frac{a_{62}}{96}+\frac{a_{63}}{288}+\frac{a_{64}}{288}+\frac{a_{65}}{144}\nonumber \\ &+\frac{a_{66}}{144}+\frac{a_{67}}{96}+\frac{a_{68}}{96}+\frac{a_{69}}{288}+\frac{a_{70}}{24}+\frac{5 a_{79}}{3456}+\frac{a_{80}}{1728}+\frac{a_{81}}{1728}\nonumber \\ &+\frac{a_{82}}{1152}+\frac{5 a_{83}}{6912}-\frac{859}{48384}.
 \end{align}
From the degree $1$ part of $\tau_{5,1}\Phi(\tau_{0,0}, \tau_{0,0})|_{t=0} = 0$, we obtain
\begin{align}
 0 &= \frac{a_1}{30240}-\frac{a_2}{2880}-\frac{a_8}{11520}-\frac{a_9}{11520}-\frac{a_{10}}{11520}+\frac{7 a_{14}}{2880}+\frac{a_{15}}{192}+\frac{7 a_{16}}{2880}\nonumber \\ &+\frac{7 a_{18}}{2880}+\frac{a_{19}}{192}-\frac{a_{21}}{11520}-\frac{a_{22}}{11520}+\frac{a_{23}}{960}+\frac{a_{24}}{960}+\frac{a_{25}}{960}+\frac{a_{29}}{960}\nonumber \\ &-\frac{a_{31}}{46080}-\frac{a_{37}}{288}-\frac{a_{43}}{288}-\frac{a_{44}}{288}-\frac{a_{50}}{288}-\frac{a_{52}}{288}-\frac{a_{53}}{288}-\frac{a_{57}}{288}\nonumber \\ &+\frac{a_{70}}{288}+\frac{a_{74}}{13824}+\frac{a_{75}}{13824}+\frac{a_{79}}{13824}-\frac{a_{82}}{13824}+\frac{a_{84}}{6}+\frac{a_{88}}{6}+\frac{a_{89}}{6}\nonumber \\ &+\frac{a_{90}}{6}+\frac{a_{92}}{12}-\frac{a_{103}}{6}+8 a_{104}+\frac{a_{105}}{13824}-\frac{1}{5040}.
 \end{align}
From the degree $1$ part of $\tau_{4,1}\Phi(\tau_{0,0}, \tau_{0,1})|_{t=0} = 0$, we obtain
\begin{align}
 0 &= \frac{11 a_1}{60480}+\frac{a_2}{960}+\frac{a_8}{34560}+\frac{a_9}{34560}-\frac{a_{10}}{11520}+\frac{19 a_{14}}{2880}+\frac{a_{15}}{960}+\frac{a_{16}}{320}\nonumber \\ &+\frac{a_{18}}{320}+\frac{a_{19}}{960}+\frac{a_{21}}{34560}-\frac{a_{22}}{11520}+\frac{7 a_{23}}{2880}+\frac{a_{24}}{960}+\frac{a_{25}}{960}+\frac{a_{29}}{576}\nonumber \\ &-\frac{a_{31}}{46080}-\frac{a_{37}}{288}-\frac{a_{44}}{288}-\frac{a_{50}}{288}-\frac{a_{52}}{288}-\frac{a_{53}}{288}-\frac{a_{57}}{288}+\frac{a_{70}}{288}\nonumber \\ &+\frac{a_{74}}{13824}+\frac{a_{75}}{13824}+\frac{a_{79}}{13824}-\frac{a_{82}}{13824}+\frac{a_{84}}{6}+\frac{a_{88}}{6}+\frac{a_{89}}{12}+\frac{a_{90}}{6}\nonumber \\ &+\frac{a_{92}}{8}-\frac{a_{93}}{24}-\frac{a_{103}}{6}+8 a_{104}+\frac{a_{105}}{13824}-\frac{11}{20160}.
 \end{align}
From the degree $1$ part of $\tau_{5,0}\Phi(\tau_{0,0}, \tau_{0,1})|_{t=0} = 0$, we obtain
\begin{align}
 0 &= -\frac{779 a_1}{1451520}+\frac{a_2}{120}-\frac{13 a_8}{34560}-\frac{13 a_9}{34560}-\frac{13 a_{10}}{34560}-\frac{a_{14}}{60}+\frac{a_{15}}{120}\nonumber \\ &+\frac{13 a_{16}}{720}+\frac{13 a_{18}}{720}+\frac{a_{19}}{120}-\frac{13 a_{21}}{34560}-\frac{13 a_{22}}{34560}-\frac{a_{23}}{120}+\frac{a_{24}}{180}+\frac{a_{25}}{180}\nonumber \\ &-\frac{a_{29}}{720}-\frac{13 a_{31}}{138240}-\frac{a_{37}}{144}-\frac{a_{43}}{144}-\frac{a_{44}}{144}-\frac{a_{50}}{144}-\frac{a_{52}}{144}-\frac{a_{53}}{144}\nonumber \\ &-\frac{a_{57}}{144}+\frac{a_{74}}{13824}+\frac{a_{75}}{13824}+\frac{a_{79}}{13824}+\frac{a_{82}}{13824}+\frac{a_{84}}{2}+\frac{a_{88}}{2}+\frac{a_{89}}{6}\nonumber \\ &+\frac{a_{90}}{2}+\frac{a_{92}}{3}+\frac{a_{93}}{12}-\frac{a_{103}}{6}+16 a_{104}+\frac{a_{105}}{13824}+\frac{61}{34560}.
 \end{align}
From the degree $1$ part of $\tau_{4,1}\Phi(\tau_{0,0}, \tau_{1,0})|_{t=0} = 0$, we obtain
\begin{align}
 0 &= -\frac{a_2}{144}+\frac{a_3}{240}-\frac{a_4}{720}-\frac{23 a_9}{69120}-\frac{a_{10}}{3840}-\frac{a_{11}}{15360}-\frac{a_{13}}{15360}\nonumber \\ &+\frac{7 a_{15}}{480}+\frac{a_{16}}{144}+\frac{a_{17}}{720}+\frac{a_{24}}{240}+\frac{a_{28}}{960}+\frac{a_{30}}{960}-\frac{a_{34}}{288}-\frac{a_{37}}{72}\nonumber \\ &-\frac{a_{38}}{288}-\frac{a_{42}}{288}-\frac{a_{44}}{72}-\frac{a_{45}}{288}-\frac{a_{53}}{72}-\frac{a_{54}}{288}-\frac{a_{55}}{288}-\frac{a_{57}}{72}\nonumber \\ &-\frac{a_{58}}{288}-\frac{a_{61}}{288}+\frac{a_{68}}{288}+\frac{a_{70}}{72}+\frac{a_{73}}{13824}+\frac{a_{75}}{4608}+\frac{a_{76}}{13824}+\frac{a_{77}}{13824}\nonumber \\ &+\frac{a_{79}}{4608}+\frac{a_{80}}{13824}-\frac{a_{81}}{13824}-\frac{a_{82}}{4608}+\frac{a_{83}}{13824}+\frac{5 a_{84}}{6}+\frac{a_{85}}{6}+\frac{5 a_{90}}{6}\nonumber \\ &+\frac{a_{91}}{6}+\frac{5 a_{92}}{12}+\frac{a_{93}}{12}-\frac{a_{102}}{6}-\frac{5 a_{103}}{6}+48 a_{104}-\frac{41}{161280}.
 \end{align}
From the degree $1$ part of $\tau_{3,1}\Phi(\tau_{0,0}, \tau_{1,1})|_{t=0} = 0$, we obtain
\begin{align}
 0 &= -\frac{31 a_1}{96768}+\frac{7 a_3}{720}+\frac{7 a_4}{720}-\frac{7 a_8}{34560}-\frac{a_9}{13824}+\frac{7 a_{11}}{138240}+\frac{7 a_{13}}{138240}\nonumber \\ &+\frac{a_{15}}{288}+\frac{7 a_{17}}{1440}-\frac{7 a_{21}}{34560}-\frac{7 a_{22}}{34560}+\frac{7 a_{28}}{2880}+\frac{7 a_{30}}{2880}-\frac{7 a_{31}}{138240}-\frac{a_{34}}{288}\nonumber \\ &-\frac{a_{37}}{144}+\frac{a_{38}}{288}-\frac{a_{42}}{288}-\frac{a_{44}}{144}-\frac{a_{53}}{144}-\frac{a_{54}}{288}-\frac{a_{55}}{288}-\frac{a_{57}}{144}\nonumber \\ &+\frac{a_{68}}{288}+\frac{a_{70}}{72}+\frac{a_{73}}{13824}-\frac{a_{74}}{13824}+\frac{a_{75}}{6912}+\frac{a_{76}}{13824}+\frac{a_{77}}{13824}+\frac{a_{79}}{6912}\nonumber \\ &+\frac{a_{80}}{13824}-\frac{a_{81}}{13824}-\frac{a_{82}}{3456}+\frac{a_{83}}{13824}+\frac{a_{84}}{3}+\frac{a_{90}}{3}+\frac{a_{91}}{12}+\frac{a_{92}}{4}\nonumber \\ &+\frac{a_{93}}{24}-\frac{a_{95}}{24}-\frac{a_{97}}{12}-\frac{a_{102}}{6}-\frac{2 a_{103}}{3}+32 a_{104}-\frac{a_{105}}{13824}-\frac{1}{4608}.
 \end{align}
From the degree $1$ part of $\tau_{4,0}\Phi(\tau_{0,0}, \tau_{1,1})|_{t=0} = 0$, we obtain
\begin{align}
 0 &= \frac{2329 a_1}{1451520}-\frac{a_2}{576}-\frac{a_3}{30}-\frac{a_4}{30}+\frac{a_8}{1440}-\frac{19 a_9}{69120}-\frac{a_{10}}{2304}-\frac{13 a_{11}}{46080}\nonumber \\ &-\frac{13 a_{13}}{46080}+\frac{7 a_{14}}{576}+\frac{11 a_{15}}{576}+\frac{7 a_{16}}{576}-\frac{11 a_{17}}{720}+\frac{7 a_{18}}{576}+\frac{7 a_{19}}{576}+\frac{a_{21}}{1440}\nonumber \\ &+\frac{a_{22}}{1440}+\frac{7 a_{23}}{2880}+\frac{a_{24}}{192}+\frac{7 a_{25}}{2880}-\frac{a_{28}}{120}+\frac{7 a_{29}}{2880}-\frac{a_{30}}{120}+\frac{a_{31}}{5760}\nonumber \\ &-\frac{a_{34}}{144}-\frac{7 a_{37}}{288}-\frac{a_{38}}{144}-\frac{a_{42}}{144}+\frac{a_{43}}{288}-\frac{7 a_{44}}{288}-\frac{a_{45}}{144}+\frac{a_{50}}{288}\nonumber \\ &+\frac{a_{52}}{288}-\frac{7 a_{53}}{288}-\frac{a_{54}}{144}-\frac{a_{55}}{144}-\frac{7 a_{57}}{288}-\frac{a_{58}}{144}-\frac{a_{61}}{144}+\frac{a_{70}}{288}\nonumber \\ &+\frac{a_{73}}{13824}+\frac{a_{74}}{6912}+\frac{5 a_{75}}{13824}+\frac{a_{76}}{13824}+\frac{a_{77}}{13824}+\frac{5 a_{79}}{13824}+\frac{a_{80}}{13824}+\frac{a_{81}}{13824}\nonumber \\ &+\frac{a_{82}}{4608}+\frac{a_{83}}{13824}+\frac{3 a_{84}}{2}-\frac{a_{85}}{6}+\frac{3 a_{90}}{2}+\frac{a_{91}}{6}+\frac{11 a_{92}}{12}+\frac{a_{95}}{12}\nonumber \\ &+\frac{a_{97}}{6}-\frac{a_{102}}{6}-\frac{5 a_{103}}{6}+72 a_{104}+\frac{a_{105}}{6912}-\frac{131}{69120}.
 \end{align}
From the degree $1$ part of $\tau_{3,1}\Phi(\tau_{0,0}, \tau_{2,0})|_{t=0} = 0$, we obtain
\begin{align}
 0 &= \frac{31 a_1}{96768}-\frac{a_2}{288}+\frac{7 a_5}{2880}+\frac{7 a_6}{960}-\frac{7 a_7}{2880}+\frac{7 a_8}{34560}-\frac{a_9}{2304}-\frac{a_{10}}{1152}\nonumber \\ &-\frac{7 a_{11}}{69120}+\frac{7 a_{15}}{288}+\frac{7 a_{16}}{288}+\frac{7 a_{21}}{34560}+\frac{7 a_{22}}{34560}+\frac{a_{24}}{96}+\frac{7 a_{31}}{138240}-\frac{a_{32}}{288}\nonumber \\ &-\frac{a_{33}}{288}+\frac{a_{34}}{144}-\frac{a_{37}}{36}+\frac{a_{39}}{288}-\frac{a_{40}}{288}+\frac{a_{42}}{144}-\frac{a_{44}}{36}-\frac{a_{47}}{288}\nonumber \\ &-\frac{a_{53}}{36}+\frac{a_{55}}{144}-\frac{a_{57}}{36}+\frac{a_{59}}{288}-\frac{a_{60}}{288}-\frac{a_{62}}{288}-\frac{a_{65}}{288}+\frac{a_{66}}{288}\nonumber \\ &-\frac{a_{67}}{288}-\frac{a_{68}}{144}+\frac{a_{70}}{48}-\frac{a_{73}}{6912}+\frac{a_{74}}{13824}+\frac{7 a_{75}}{13824}-\frac{a_{77}}{6912}+\frac{7 a_{79}}{13824}\nonumber \\ &+\frac{a_{81}}{6912}+\frac{a_{82}}{13824}-\frac{a_{83}}{6912}+\frac{5 a_{84}}{3}+\frac{a_{87}}{6}+\frac{5 a_{90}}{3}+\frac{5 a_{92}}{6}+\frac{a_{95}}{12}\nonumber \\ &-\frac{a_{98}}{6}-\frac{a_{100}}{6}+\frac{a_{102}}{3}-\frac{5 a_{103}}{3}+112 a_{104}+\frac{a_{105}}{13824}-\frac{53}{48384}.
 \end{align}
From the degree $1$ part of $\tau_{2,1}\Phi(\tau_{0,0}, \tau_{2,1})|_{t=0} = 0$, we obtain
\begin{align}
 0 &= \frac{a_2}{288}+\frac{7 a_5}{2880}+\frac{7 a_6}{960}+\frac{7 a_7}{960}-\frac{a_9}{13824}-\frac{a_{10}}{3456}+\frac{7 a_{11}}{46080}-\frac{7 a_{13}}{46080}\nonumber \\ &+\frac{a_{15}}{288}+\frac{a_{16}}{96}+\frac{a_{24}}{288}-\frac{a_{32}}{288}-\frac{a_{33}}{288}-\frac{a_{37}}{144}+\frac{a_{39}}{288}+\frac{a_{40}}{288}\nonumber \\ &-\frac{a_{44}}{144}-\frac{a_{47}}{288}-\frac{a_{53}}{144}-\frac{a_{57}}{144}+\frac{a_{66}}{288}+\frac{a_{67}}{288}+\frac{a_{70}}{48}+\frac{a_{73}}{13824}\nonumber \\ &+\frac{a_{75}}{4608}-\frac{a_{76}}{13824}+\frac{a_{77}}{13824}+\frac{a_{79}}{4608}-\frac{a_{80}}{13824}+\frac{a_{81}}{13824}-\frac{a_{82}}{4608}+\frac{a_{83}}{13824}\nonumber \\ &+\frac{a_{84}}{3}+\frac{a_{90}}{3}+\frac{a_{92}}{4}+\frac{a_{94}}{24}+\frac{a_{95}}{24}-\frac{a_{98}}{12}-\frac{a_{99}}{12}-\frac{a_{100}}{6}\nonumber \\ &-\frac{a_{101}}{6}-a_{103}+48 a_{104}-\frac{1}{4608}.
 \end{align}
From the degree $1$ part of $\tau_{3,0}\Phi(\tau_{0,0}, \tau_{2,1})|_{t=0} = 0$, we obtain
\begin{align}
 0 &= -\frac{31 a_1}{96768}+\frac{7 a_3}{720}+\frac{7 a_4}{720}-\frac{a_7}{40}-\frac{7 a_8}{34560}-\frac{23 a_9}{69120}-\frac{13 a_{10}}{23040}\nonumber \\ &-\frac{a_{11}}{1920}+\frac{a_{13}}{1920}+\frac{7 a_{15}}{288}+\frac{5 a_{16}}{144}+\frac{7 a_{17}}{720}-\frac{7 a_{21}}{34560}-\frac{7 a_{22}}{34560}+\frac{a_{24}}{72}\nonumber \\ &+\frac{7 a_{28}}{2880}+\frac{7 a_{30}}{2880}-\frac{7 a_{31}}{138240}-\frac{a_{32}}{144}+\frac{a_{33}}{144}-\frac{a_{34}}{288}-\frac{5 a_{37}}{144}+\frac{a_{38}}{288}\nonumber \\ &-\frac{a_{40}}{144}-\frac{a_{42}}{288}-\frac{5 a_{44}}{144}+\frac{a_{45}}{288}-\frac{5 a_{53}}{144}+\frac{a_{54}}{288}-\frac{a_{55}}{288}-\frac{5 a_{57}}{144}\nonumber \\ &+\frac{a_{58}}{288}-\frac{a_{60}}{144}-\frac{a_{61}}{288}-\frac{a_{65}}{144}-\frac{a_{67}}{144}-\frac{a_{68}}{288}+\frac{a_{70}}{72}-\frac{a_{73}}{6912}\nonumber \\ &-\frac{a_{74}}{13824}+\frac{11 a_{75}}{13824}+\frac{a_{76}}{6912}-\frac{a_{77}}{6912}+\frac{11 a_{79}}{13824}+\frac{a_{80}}{6912}+\frac{5 a_{82}}{13824}-\frac{a_{83}}{6912}\nonumber \\ &+2 a_{84}-\frac{a_{87}}{6}+2 a_{90}+\frac{13 a_{92}}{12}-\frac{a_{97}}{6}+\frac{a_{99}}{6}-\frac{a_{100}}{6}+\frac{a_{101}}{3}\nonumber \\ &-\frac{5 a_{103}}{3}+128 a_{104}-\frac{a_{105}}{13824}-\frac{53}{48384}.
 \end{align}
From the degree $1$ part of $\tau_{2,1}\Phi(\tau_{0,0}, \tau_{3,0})|_{t=0} = 0$, we obtain
\begin{align}
 0 &= \frac{7 a_2}{720}-\frac{a_5}{120}-\frac{a_6}{40}-\frac{a_9}{3072}-\frac{113 a_{10}}{138240}+\frac{7 a_{11}}{46080}+\frac{7 a_{13}}{46080}+\frac{31 a_{15}}{1440}\nonumber \\ &+\frac{2 a_{16}}{45}+\frac{3 a_{24}}{160}+\frac{a_{32}}{144}-\frac{a_{33}}{144}-\frac{5 a_{37}}{144}-\frac{a_{39}}{144}-\frac{5 a_{44}}{144}-\frac{5 a_{53}}{144}\nonumber \\ &-\frac{5 a_{57}}{144}-\frac{a_{59}}{144}-\frac{a_{62}}{144}-\frac{a_{66}}{144}+\frac{a_{70}}{144}+\frac{a_{73}}{13824}+\frac{a_{75}}{2304}+\frac{a_{76}}{13824}\nonumber \\ &+\frac{a_{77}}{13824}+\frac{a_{79}}{2304}+\frac{a_{80}}{13824}+\frac{a_{81}}{13824}+\frac{a_{82}}{6912}+\frac{a_{83}}{13824}+2 a_{84}-\frac{a_{86}}{6}\nonumber \\ &+2 a_{90}+\frac{13 a_{92}}{12}-\frac{a_{96}}{6}+\frac{a_{98}}{6}+\frac{a_{100}}{3}-\frac{a_{101}}{6}-\frac{5 a_{103}}{3}+128 a_{104}\nonumber \\ &-\frac{127}{80640}.
 \end{align}
From the degree $1$ part of $\tau_{2,0}\Phi(\tau_{0,0}, \tau_{3,1})|_{t=0} = 0$, we obtain
\begin{align}
 0 &= -\frac{a_2}{288}+\frac{a_5}{960}-\frac{7 a_6}{2880}+\frac{7 a_7}{960}-\frac{11 a_9}{46080}-\frac{47 a_{10}}{46080}+\frac{7 a_{11}}{46080}\nonumber \\ &-\frac{7 a_{13}}{46080}+\frac{7 a_{15}}{480}+\frac{7 a_{16}}{288}+\frac{a_{24}}{96}-\frac{a_{32}}{288}-\frac{a_{33}}{288}-\frac{a_{37}}{48}-\frac{a_{39}}{288}\nonumber \\ &+\frac{a_{40}}{288}-\frac{a_{44}}{48}-\frac{a_{47}}{288}-\frac{a_{53}}{48}-\frac{a_{57}}{36}-\frac{a_{59}}{288}+\frac{a_{60}}{288}-\frac{a_{62}}{288}\nonumber \\ &-\frac{a_{65}}{288}-\frac{a_{66}}{288}+\frac{a_{67}}{288}+\frac{a_{70}}{144}+\frac{a_{73}}{13824}+\frac{a_{75}}{6912}-\frac{a_{76}}{13824}+\frac{a_{77}}{13824}\nonumber \\ &+\frac{a_{79}}{6912}-\frac{a_{80}}{13824}+\frac{a_{81}}{13824}+\frac{a_{82}}{3456}+\frac{a_{83}}{13824}+\frac{5 a_{84}}{3}+\frac{a_{86}}{6}+\frac{5 a_{90}}{3}\nonumber \\ &+\frac{5 a_{92}}{6}+\frac{a_{94}}{12}-\frac{a_{99}}{6}-\frac{a_{101}}{6}-\frac{4 a_{103}}{3}+112 a_{104}-\frac{139}{241920}.
 \end{align}
From the degree $1$ part of $\tau_{4,1}\Phi(\tau_{0,1}, \tau_{0,0})|_{t=0} = 0$, we obtain
\begin{align}
 0 &= \frac{11 a_1}{60480}+\frac{a_2}{960}-\frac{a_8}{11520}-\frac{a_9}{11520}+\frac{a_{10}}{34560}+\frac{a_{14}}{320}+\frac{5 a_{15}}{576}+\frac{19 a_{16}}{2880}\nonumber \\ &+\frac{a_{18}}{320}+\frac{a_{19}}{960}+\frac{a_{21}}{34560}-\frac{a_{22}}{11520}+\frac{a_{23}}{960}+\frac{7 a_{24}}{2880}+\frac{a_{25}}{960}+\frac{a_{29}}{576}\nonumber \\ &-\frac{a_{31}}{46080}-\frac{a_{37}}{288}-\frac{a_{43}}{288}-\frac{a_{50}}{288}-\frac{a_{52}}{288}-\frac{a_{53}}{288}-\frac{a_{57}}{288}+\frac{a_{70}}{288}\nonumber \\ &+\frac{a_{74}}{13824}+\frac{a_{75}}{13824}+\frac{a_{79}}{13824}-\frac{a_{82}}{13824}+\frac{a_{84}}{6}+\frac{a_{88}}{6}+\frac{a_{89}}{6}+\frac{a_{90}}{12}\nonumber \\ &+\frac{a_{92}}{8}-\frac{a_{94}}{24}-\frac{a_{103}}{6}+8 a_{104}+\frac{a_{105}}{13824}-\frac{11}{10080}.
 \end{align}
From the degree $1$ part of $\tau_{5,0}\Phi(\tau_{0,1}, \tau_{0,0})|_{t=0} = 0$, we obtain
\begin{align}
 0 &= -\frac{779 a_1}{1451520}+\frac{a_2}{120}-\frac{13 a_8}{34560}-\frac{13 a_9}{34560}-\frac{13 a_{10}}{34560}+\frac{13 a_{14}}{720}-\frac{a_{15}}{144}\nonumber \\ &-\frac{a_{16}}{60}+\frac{13 a_{18}}{720}+\frac{a_{19}}{120}-\frac{13 a_{21}}{34560}-\frac{13 a_{22}}{34560}+\frac{a_{23}}{180}-\frac{a_{24}}{120}+\frac{a_{25}}{180}\nonumber \\ &-\frac{a_{29}}{720}-\frac{13 a_{31}}{138240}-\frac{a_{37}}{144}-\frac{a_{43}}{144}-\frac{a_{44}}{144}-\frac{a_{50}}{144}-\frac{a_{52}}{144}-\frac{a_{53}}{144}\nonumber \\ &-\frac{a_{57}}{144}+\frac{a_{74}}{13824}+\frac{a_{75}}{13824}+\frac{a_{79}}{13824}+\frac{a_{82}}{13824}+\frac{a_{84}}{2}+\frac{a_{88}}{2}+\frac{a_{89}}{2}\nonumber \\ &+\frac{a_{90}}{6}+\frac{a_{92}}{3}+\frac{a_{94}}{12}-\frac{a_{103}}{6}+16 a_{104}+\frac{a_{105}}{13824}+\frac{779}{241920}.
 \end{align}
From the degree $1$ part of $\tau_{3,1}\Phi(\tau_{0,1}, \tau_{0,1})|_{t=0} = 0$, we obtain
\begin{align}
 0 &= -\frac{31 a_1}{96768}+\frac{7 a_{17}}{1440}-\frac{7 a_{21}}{34560}-\frac{7 a_{22}}{34560}+\frac{7 a_{27}}{5760}+\frac{7 a_{28}}{5760}+\frac{7 a_{30}}{2880}\nonumber \\ &-\frac{7 a_{31}}{138240}-\frac{a_{37}}{288}-\frac{a_{43}}{576}-\frac{a_{44}}{576}+\frac{a_{45}}{576}-\frac{a_{52}}{288}-\frac{a_{53}}{288}-\frac{a_{57}}{288}\nonumber \\ &+\frac{a_{70}}{288}+\frac{a_{74}}{13824}+\frac{a_{75}}{13824}+\frac{a_{79}}{13824}-\frac{a_{82}}{13824}+\frac{a_{84}}{6}+\frac{a_{89}}{12}+\frac{a_{90}}{12}\nonumber \\ &+\frac{a_{92}}{8}-\frac{a_{95}}{24}-\frac{a_{103}}{6}+8 a_{104}-\frac{a_{105}}{13824}.
 \end{align}
From the degree $1$ part of $\tau_{4,0}\Phi(\tau_{0,1}, \tau_{0,1})|_{t=0} = 0$, we obtain
\begin{align}
 0 &= \frac{2329 a_1}{1451520}+\frac{a_2}{192}-\frac{a_8}{17280}-\frac{a_9}{17280}-\frac{a_{10}}{17280}+\frac{5 a_{14}}{576}+\frac{5 a_{15}}{576}+\frac{5 a_{16}}{576}\nonumber \\ &-\frac{a_{17}}{60}+\frac{7 a_{18}}{576}+\frac{7 a_{19}}{576}+\frac{a_{21}}{1440}+\frac{a_{22}}{1440}+\frac{7 a_{23}}{2880}+\frac{7 a_{24}}{2880}+\frac{7 a_{25}}{2880}\nonumber \\ &-\frac{a_{27}}{240}-\frac{a_{28}}{240}+\frac{7 a_{29}}{2880}-\frac{a_{30}}{120}+\frac{a_{31}}{5760}-\frac{a_{37}}{96}-\frac{a_{43}}{288}-\frac{a_{44}}{288}\nonumber \\ &-\frac{a_{45}}{288}+\frac{a_{50}}{288}-\frac{a_{52}}{96}-\frac{a_{53}}{96}-\frac{a_{57}}{96}+\frac{a_{70}}{288}+\frac{a_{74}}{6912}+\frac{a_{75}}{6912}\nonumber \\ &+\frac{a_{79}}{6912}+\frac{2 a_{84}}{3}+\frac{a_{89}}{3}+\frac{a_{90}}{3}-\frac{a_{91}}{12}+\frac{a_{92}}{2}-\frac{a_{93}}{24}-\frac{a_{94}}{24}\nonumber \\ &+\frac{a_{95}}{12}-\frac{a_{103}}{3}+24 a_{104}+\frac{a_{105}}{6912}-\frac{25}{16128}.
 \end{align}
From the degree $1$ part of $\tau_{3,1}\Phi(\tau_{0,1}, \tau_{1,0})|_{t=0} = 0$, we obtain
\begin{align}
 0 &= -\frac{a_9}{4608}+\frac{7 a_{10}}{34560}+\frac{7 a_{11}}{138240}+\frac{7 a_{13}}{138240}+\frac{5 a_{15}}{288}+\frac{7 a_{16}}{288}+\frac{7 a_{17}}{1440}\nonumber \\ &+\frac{7 a_{24}}{720}+\frac{7 a_{28}}{2880}+\frac{7 a_{30}}{2880}-\frac{a_{34}}{288}-\frac{a_{37}}{72}-\frac{a_{38}}{288}-\frac{a_{53}}{72}-\frac{a_{54}}{288}\nonumber \\ &-\frac{a_{55}}{288}-\frac{a_{57}}{72}-\frac{a_{58}}{288}-\frac{a_{61}}{288}+\frac{a_{68}}{288}+\frac{a_{70}}{72}+\frac{a_{73}}{13824}+\frac{a_{75}}{4608}\nonumber \\ &+\frac{a_{76}}{13824}+\frac{a_{77}}{13824}+\frac{a_{79}}{4608}+\frac{a_{80}}{13824}-\frac{a_{81}}{13824}-\frac{a_{82}}{4608}+\frac{a_{83}}{13824}+\frac{5 a_{84}}{6}\nonumber \\ &+\frac{a_{85}}{6}+\frac{5 a_{90}}{12}+\frac{a_{91}}{12}+\frac{5 a_{92}}{8}+\frac{a_{93}}{8}-\frac{a_{94}}{6}-\frac{a_{95}}{12}-\frac{a_{102}}{6}\nonumber \\ &-\frac{5 a_{103}}{6}+48 a_{104}-\frac{23}{10752}.
 \end{align}
From the degree $1$ part of $\tau_{4,0}\Phi(\tau_{0,1}, \tau_{1,0})|_{t=0} = 0$, we obtain
\begin{align}
 0 &= \frac{a_2}{72}+\frac{a_3}{360}+\frac{a_4}{360}-\frac{67 a_9}{69120}-\frac{13 a_{10}}{11520}-\frac{13 a_{11}}{46080}-\frac{13 a_{13}}{46080}-\frac{a_{15}}{60}\nonumber \\ &-\frac{11 a_{16}}{144}-\frac{11 a_{17}}{720}-\frac{a_{24}}{30}-\frac{a_{28}}{120}-\frac{a_{30}}{120}-\frac{a_{34}}{144}-\frac{a_{37}}{36}-\frac{a_{38}}{144}\nonumber \\ &-\frac{a_{42}}{144}-\frac{a_{44}}{36}-\frac{a_{45}}{144}-\frac{a_{53}}{36}-\frac{a_{54}}{144}-\frac{a_{55}}{144}-\frac{a_{57}}{36}-\frac{a_{58}}{144}\nonumber \\ &-\frac{a_{61}}{144}+\frac{a_{73}}{13824}+\frac{a_{75}}{4608}+\frac{a_{76}}{13824}+\frac{a_{77}}{13824}+\frac{a_{79}}{4608}+\frac{a_{80}}{13824}+\frac{a_{81}}{13824}\nonumber \\ &+\frac{a_{82}}{4608}+\frac{a_{83}}{13824}+\frac{5 a_{84}}{2}+\frac{a_{85}}{2}+\frac{5 a_{90}}{6}+\frac{a_{91}}{6}+\frac{5 a_{92}}{3}+\frac{a_{93}}{3}\nonumber \\ &+\frac{a_{94}}{3}+\frac{a_{95}}{6}-\frac{a_{102}}{6}-\frac{5 a_{103}}{6}+96 a_{104}+\frac{1247}{161280}.
 \end{align}
From the degree $1$ part of $\tau_{2,1}\Phi(\tau_{0,1}, \tau_{1,1})|_{t=0} = 0$, we obtain
\begin{align}
 0 &= \frac{7 a_7}{960}-\frac{a_9}{13824}-\frac{7 a_{11}}{46080}-\frac{7 a_{13}}{46080}+\frac{a_{15}}{576}-\frac{a_{34}}{288}-\frac{a_{37}}{288}+\frac{a_{40}}{288}\nonumber \\ &-\frac{a_{42}}{576}-\frac{a_{44}}{576}-\frac{a_{53}}{288}-\frac{a_{57}}{288}-\frac{a_{58}}{576}+\frac{a_{60}}{576}-\frac{a_{61}}{576}+\frac{a_{63}}{576}\nonumber \\ &+\frac{a_{68}}{288}+\frac{a_{70}}{96}+\frac{a_{73}}{13824}+\frac{a_{75}}{13824}-\frac{a_{76}}{13824}-\frac{a_{77}}{13824}+\frac{a_{79}}{13824}+\frac{a_{80}}{13824}\nonumber \\ &-\frac{a_{81}}{13824}-\frac{a_{82}}{4608}+\frac{a_{83}}{13824}+\frac{a_{84}}{6}+\frac{a_{90}}{12}+\frac{a_{92}}{8}+\frac{a_{93}}{24}-\frac{a_{99}}{12}\nonumber \\ &-\frac{a_{102}}{6}-\frac{a_{103}}{2}+24 a_{104}-\frac{1}{13824}.
 \end{align}
From the degree $1$ part of $\tau_{3,0}\Phi(\tau_{0,1}, \tau_{1,1})|_{t=0} = 0$, we obtain
\begin{align}
 0 &= -\frac{31 a_1}{96768}+\frac{7 a_3}{720}+\frac{7 a_4}{720}-\frac{a_7}{40}-\frac{7 a_8}{34560}-\frac{a_9}{3456}+\frac{7 a_{10}}{34560}\nonumber \\ &+\frac{a_{11}}{1920}+\frac{a_{13}}{1920}+\frac{a_{15}}{48}+\frac{7 a_{16}}{288}+\frac{7 a_{17}}{720}-\frac{7 a_{21}}{34560}-\frac{7 a_{22}}{34560}+\frac{7 a_{24}}{720}\nonumber \\ &+\frac{7 a_{28}}{2880}+\frac{7 a_{30}}{2880}-\frac{7 a_{31}}{138240}-\frac{a_{34}}{96}-\frac{a_{37}}{48}+\frac{a_{38}}{288}-\frac{a_{40}}{144}-\frac{a_{42}}{288}\nonumber \\ &-\frac{a_{44}}{288}+\frac{a_{45}}{288}-\frac{a_{53}}{48}+\frac{a_{54}}{288}+\frac{a_{55}}{288}-\frac{a_{57}}{48}-\frac{a_{58}}{288}-\frac{a_{60}}{288}\nonumber \\ &-\frac{a_{61}}{288}-\frac{a_{63}}{288}+\frac{a_{68}}{288}+\frac{a_{70}}{72}+\frac{a_{73}}{6912}-\frac{a_{74}}{13824}+\frac{5 a_{75}}{13824}+\frac{a_{76}}{6912}\nonumber \\ &+\frac{a_{77}}{6912}+\frac{5 a_{79}}{13824}+\frac{a_{80}}{6912}-\frac{a_{82}}{13824}+\frac{a_{83}}{6912}+\frac{7 a_{84}}{6}-\frac{a_{87}}{6}+\frac{7 a_{90}}{12}\nonumber \\ &+\frac{7 a_{92}}{8}+\frac{a_{93}}{6}-\frac{a_{94}}{6}-\frac{a_{95}}{8}-\frac{a_{97}}{12}+\frac{a_{99}}{6}-\frac{a_{102}}{3}-\frac{7 a_{103}}{6}\nonumber \\ &+80 a_{104}-\frac{a_{105}}{13824}-\frac{19}{8064}.
 \end{align}
From the degree $1$ part of $\tau_{2,1}\Phi(\tau_{0,1}, \tau_{2,0})|_{t=0} = 0$, we obtain
\begin{align}
 0 &= \frac{a_2}{96}+\frac{7 a_5}{2880}+\frac{7 a_6}{960}+\frac{7 a_7}{960}-\frac{a_9}{3456}+\frac{a_{10}}{11520}+\frac{7 a_{11}}{23040}+\frac{a_{15}}{48}\nonumber \\ &+\frac{a_{16}}{24}+\frac{7 a_{24}}{480}-\frac{a_{32}}{288}-\frac{a_{33}}{288}+\frac{a_{34}}{144}-\frac{a_{37}}{48}+\frac{a_{39}}{288}+\frac{a_{40}}{288}\nonumber \\ &+\frac{a_{42}}{288}-\frac{a_{53}}{48}-\frac{a_{57}}{48}+\frac{a_{59}}{288}+\frac{a_{61}}{288}-\frac{a_{62}}{288}-\frac{a_{63}}{288}+\frac{a_{66}}{288}\nonumber \\ &-\frac{a_{67}}{288}-\frac{a_{68}}{144}+\frac{a_{70}}{48}-\frac{a_{73}}{6912}+\frac{a_{75}}{2304}+\frac{a_{77}}{6912}+\frac{a_{79}}{2304}+\frac{a_{81}}{6912}\nonumber \\ &-\frac{a_{83}}{6912}+\frac{7 a_{84}}{6}+\frac{7 a_{90}}{12}+\frac{7 a_{92}}{8}-\frac{a_{94}}{4}+\frac{a_{95}}{24}+\frac{a_{97}}{12}-\frac{a_{98}}{6}\nonumber \\ &-\frac{a_{99}}{6}-\frac{a_{100}}{6}+\frac{a_{102}}{3}-\frac{3 a_{103}}{2}+96 a_{104}-\frac{19}{8064}.
 \end{align}
From the degree $1$ part of $\tau_{3,0}\Phi(\tau_{0,1}, \tau_{2,0})|_{t=0} = 0$, we obtain
\begin{align}
 0 &= \frac{31 a_1}{96768}+\frac{7 a_2}{144}-\frac{a_7}{40}+\frac{7 a_8}{34560}-\frac{3 a_9}{2560}-\frac{31 a_{10}}{13824}-\frac{a_{11}}{960}-\frac{11 a_{15}}{720}\nonumber \\ &-\frac{a_{16}}{8}+\frac{7 a_{21}}{34560}+\frac{7 a_{22}}{34560}-\frac{43 a_{24}}{720}+\frac{7 a_{31}}{138240}-\frac{a_{32}}{144}-\frac{a_{33}}{144}+\frac{a_{34}}{48}\nonumber \\ &-\frac{7 a_{37}}{144}-\frac{a_{40}}{144}+\frac{a_{42}}{144}-\frac{7 a_{44}}{144}-\frac{a_{47}}{144}-\frac{7 a_{53}}{144}-\frac{a_{55}}{144}-\frac{7 a_{57}}{144}\nonumber \\ &-\frac{a_{60}}{144}+\frac{a_{61}}{144}-\frac{a_{62}}{144}+\frac{a_{63}}{144}-\frac{a_{65}}{144}-\frac{a_{67}}{144}-\frac{a_{68}}{144}-\frac{a_{73}}{3456}\nonumber \\ &+\frac{a_{74}}{13824}+\frac{7 a_{75}}{13824}-\frac{a_{77}}{3456}+\frac{7 a_{79}}{13824}+\frac{7 a_{82}}{13824}-\frac{a_{83}}{3456}+\frac{13 a_{84}}{3}-\frac{a_{87}}{6}\nonumber \\ &+\frac{4 a_{90}}{3}+\frac{17 a_{92}}{6}+\frac{7 a_{94}}{12}+\frac{a_{97}}{6}-\frac{a_{98}}{6}+\frac{a_{99}}{3}-\frac{a_{100}}{6}+\frac{2 a_{102}}{3}\nonumber \\ &-\frac{5 a_{103}}{3}+208 a_{104}+\frac{a_{105}}{13824}+\frac{139}{15120}.
 \end{align}
From the degree $1$ part of $\tau_{2,0}\Phi(\tau_{0,1}, \tau_{2,1})|_{t=0} = 0$, we obtain
\begin{align}
 0 &= \frac{a_2}{96}+\frac{7 a_5}{2880}+\frac{7 a_6}{960}+\frac{7 a_7}{960}-\frac{a_9}{8640}+\frac{a_{10}}{2560}+\frac{7 a_{11}}{46080}-\frac{7 a_{13}}{46080}\nonumber \\ &+\frac{5 a_{15}}{288}+\frac{a_{16}}{24}+\frac{7 a_{24}}{480}-\frac{a_{32}}{288}-\frac{a_{33}}{288}-\frac{a_{37}}{72}+\frac{a_{39}}{288}+\frac{a_{40}}{288}\nonumber \\ &+\frac{a_{44}}{288}-\frac{a_{53}}{48}-\frac{5 a_{57}}{288}-\frac{a_{59}}{288}+\frac{a_{60}}{288}-\frac{a_{65}}{288}-\frac{a_{66}}{288}+\frac{a_{67}}{288}\nonumber \\ &+\frac{a_{70}}{144}+\frac{a_{73}}{13824}+\frac{5 a_{75}}{13824}-\frac{a_{76}}{13824}+\frac{a_{77}}{13824}+\frac{a_{79}}{13824}-\frac{a_{80}}{13824}+\frac{a_{81}}{13824}\nonumber \\ &+\frac{a_{82}}{4608}+\frac{a_{83}}{13824}+\frac{7 a_{84}}{6}+\frac{7 a_{90}}{12}+\frac{7 a_{92}}{8}-\frac{5 a_{94}}{24}+\frac{a_{96}}{12}-\frac{a_{98}}{6}\nonumber \\ &-\frac{a_{99}}{6}-\frac{a_{101}}{6}-\frac{7 a_{103}}{6}+96 a_{104}-\frac{23}{10752}.
 \end{align}
From the degree $1$ part of $\tau_{2,0}\Phi(\tau_{0,1}, \tau_{3,0})|_{t=0} = 0$, we obtain
\begin{align}
 0 &= \frac{7 a_2}{144}-\frac{a_5}{120}-\frac{a_6}{40}-\frac{35 a_9}{27648}-\frac{503 a_{10}}{138240}+\frac{7 a_{11}}{46080}+\frac{7 a_{13}}{46080}-\frac{a_{15}}{60}\nonumber \\ &-\frac{a_{16}}{8}-\frac{43 a_{24}}{720}-\frac{a_{32}}{144}-\frac{a_{33}}{144}-\frac{a_{37}}{36}-\frac{a_{39}}{144}-\frac{a_{44}}{24}-\frac{a_{47}}{144}\nonumber \\ &-\frac{a_{53}}{18}-\frac{a_{57}}{24}-\frac{a_{62}}{144}-\frac{a_{65}}{144}-\frac{a_{66}}{144}-\frac{a_{70}}{72}+\frac{a_{73}}{13824}-\frac{a_{75}}{3456}\nonumber \\ &+\frac{a_{76}}{13824}+\frac{a_{77}}{13824}-\frac{a_{79}}{3456}+\frac{a_{80}}{13824}+\frac{a_{81}}{13824}+\frac{a_{82}}{3456}+\frac{a_{83}}{13824}+\frac{13 a_{84}}{3}\nonumber \\ &-\frac{a_{86}}{6}+\frac{4 a_{90}}{3}+\frac{17 a_{92}}{6}+\frac{7 a_{94}}{12}+\frac{a_{96}}{6}+\frac{a_{98}}{3}-\frac{a_{99}}{6}-\frac{a_{101}}{6}\nonumber \\ &-a_{103}+208 a_{104}+\frac{1793}{241920}.
 \end{align}
From the degree $1$ part of $\tau_{4,1}\Phi(\tau_{1,0}, \tau_{0,0})|_{t=0} = 0$, we obtain
\begin{align}
 0 &= -\frac{a_2}{144}+\frac{a_5}{240}-\frac{a_6}{720}-\frac{a_8}{3840}-\frac{a_{12}}{15360}-\frac{a_{13}}{15360}+\frac{a_{14}}{144}\nonumber \\ &+\frac{a_{17}}{720}+\frac{a_{23}}{240}+\frac{a_{27}}{960}+\frac{a_{28}}{960}+\frac{a_{30}}{960}-\frac{a_{33}}{288}-\frac{a_{37}}{72}-\frac{a_{39}}{288}\nonumber \\ &-\frac{a_{41}}{288}-\frac{a_{43}}{72}-\frac{a_{45}}{288}-\frac{a_{52}}{72}-\frac{a_{54}}{288}-\frac{a_{56}}{288}-\frac{a_{57}}{72}-\frac{a_{59}}{288}\nonumber \\ &-\frac{a_{62}}{288}+\frac{a_{67}}{288}+\frac{a_{70}}{72}+\frac{a_{72}}{13824}+\frac{a_{74}}{4608}+\frac{a_{76}}{13824}+\frac{a_{78}}{13824}+\frac{a_{79}}{4608}\nonumber \\ &-\frac{a_{80}}{13824}+\frac{a_{81}}{13824}-\frac{a_{82}}{4608}+\frac{a_{83}}{13824}+\frac{5 a_{84}}{6}+\frac{a_{86}}{6}+\frac{5 a_{89}}{6}+\frac{a_{91}}{6}\nonumber \\ &+\frac{5 a_{92}}{12}+\frac{a_{94}}{12}-\frac{a_{101}}{6}-\frac{5 a_{103}}{6}+48 a_{104}+\frac{5}{16128}.
 \end{align}
From the degree $1$ part of $\tau_{3,1}\Phi(\tau_{1,0}, \tau_{0,1})|_{t=0} = 0$, we obtain
\begin{align}
 0 &= \frac{7 a_8}{34560}+\frac{7 a_{12}}{138240}+\frac{7 a_{13}}{138240}+\frac{7 a_{14}}{288}+\frac{7 a_{17}}{1440}+\frac{7 a_{23}}{720}+\frac{7 a_{27}}{2880}+\frac{7 a_{28}}{2880}\nonumber \\ &+\frac{7 a_{30}}{2880}-\frac{a_{33}}{288}-\frac{a_{37}}{72}-\frac{a_{39}}{288}-\frac{a_{52}}{72}-\frac{a_{54}}{288}-\frac{a_{56}}{288}-\frac{a_{57}}{72}\nonumber \\ &-\frac{a_{59}}{288}-\frac{a_{62}}{288}+\frac{a_{67}}{288}+\frac{a_{70}}{72}+\frac{a_{72}}{13824}+\frac{a_{74}}{4608}+\frac{a_{76}}{13824}+\frac{a_{78}}{13824}\nonumber \\ &+\frac{a_{79}}{4608}-\frac{a_{80}}{13824}+\frac{a_{81}}{13824}-\frac{a_{82}}{4608}+\frac{a_{83}}{13824}+\frac{5 a_{84}}{6}+\frac{a_{86}}{6}+\frac{5 a_{89}}{12}\nonumber \\ &+\frac{a_{91}}{12}+\frac{5 a_{92}}{8}-\frac{a_{93}}{6}+\frac{a_{94}}{8}-\frac{a_{95}}{12}-\frac{a_{101}}{6}-\frac{5 a_{103}}{6}+48 a_{104}\nonumber \\ &-\frac{31}{48384}.
 \end{align}
From the degree $1$ part of $\tau_{4,0}\Phi(\tau_{1,0}, \tau_{0,1})|_{t=0} = 0$, we obtain
\begin{align}
 0 &= \frac{a_2}{72}+\frac{a_5}{360}+\frac{a_6}{360}-\frac{13 a_8}{11520}-\frac{13 a_{12}}{46080}-\frac{13 a_{13}}{46080}-\frac{11 a_{14}}{144}-\frac{11 a_{17}}{720}\nonumber \\ &-\frac{a_{23}}{30}-\frac{a_{27}}{120}-\frac{a_{28}}{120}-\frac{a_{30}}{120}-\frac{a_{33}}{144}-\frac{a_{37}}{36}-\frac{a_{39}}{144}-\frac{a_{41}}{144}\nonumber \\ &-\frac{a_{43}}{36}-\frac{a_{45}}{144}-\frac{a_{52}}{36}-\frac{a_{54}}{144}-\frac{a_{56}}{144}-\frac{a_{57}}{36}-\frac{a_{59}}{144}-\frac{a_{62}}{144}\nonumber \\ &+\frac{a_{72}}{13824}+\frac{a_{74}}{4608}+\frac{a_{76}}{13824}+\frac{a_{78}}{13824}+\frac{a_{79}}{4608}+\frac{a_{80}}{13824}+\frac{a_{81}}{13824}+\frac{a_{82}}{4608}\nonumber \\ &+\frac{a_{83}}{13824}+\frac{5 a_{84}}{2}+\frac{a_{86}}{2}+\frac{5 a_{89}}{6}+\frac{a_{91}}{6}+\frac{5 a_{92}}{3}+\frac{a_{93}}{3}+\frac{a_{94}}{3}\nonumber \\ &+\frac{a_{95}}{6}-\frac{a_{101}}{6}-\frac{5 a_{103}}{6}+96 a_{104}+\frac{911}{241920}.
 \end{align}
From the degree $1$ part of $\tau_{3,1}\Phi(\tau_{1,0}, \tau_{1,0})|_{t=0} = 0$, we obtain
\begin{align}
 0 &= -\frac{7 a_2}{144}+\frac{a_3}{96}-\frac{7 a_4}{720}+\frac{a_5}{96}-\frac{7 a_6}{720}-\frac{7 a_7}{1440}-\frac{a_{32}}{144}\nonumber \\ &-\frac{a_{33}}{96}-\frac{a_{34}}{96}-\frac{a_{35}}{288}-\frac{a_{36}}{288}-\frac{a_{37}}{24}-\frac{a_{38}}{96}-\frac{a_{39}}{96}-\frac{a_{40}}{144}\nonumber \\ &-\frac{a_{57}}{24}-\frac{a_{58}}{96}-\frac{a_{59}}{96}-\frac{a_{60}}{144}-\frac{a_{61}}{96}-\frac{a_{62}}{96}-\frac{a_{63}}{288}-\frac{a_{64}}{288}\nonumber \\ &-\frac{a_{65}}{144}+\frac{a_{66}}{144}+\frac{a_{67}}{96}+\frac{a_{68}}{96}+\frac{a_{69}}{288}+\frac{a_{70}}{24}+\frac{a_{79}}{1728}+\frac{a_{80}}{6912}\nonumber \\ &+\frac{a_{81}}{6912}+\frac{a_{83}}{3456}+\frac{10 a_{84}}{3}+\frac{2 a_{85}}{3}+\frac{2 a_{86}}{3}+\frac{a_{87}}{3}+\frac{5 a_{92}}{3}+\frac{a_{93}}{3}\nonumber \\ &+\frac{a_{94}}{3}+\frac{a_{95}}{6}-\frac{a_{100}}{3}-\frac{2 a_{101}}{3}-\frac{2 a_{102}}{3}-\frac{10 a_{103}}{3}+240 a_{104}+\frac{5}{6048}.
 \end{align}
From the degree $1$ part of $\tau_{2,1}\Phi(\tau_{1,0}, \tau_{1,1})|_{t=0} = 0$, we obtain
\begin{align}
 0 &= \frac{7 a_3}{480}+\frac{7 a_4}{240}+\frac{a_5}{288}+\frac{7 a_7}{480}-\frac{7 a_8}{11520}-\frac{7 a_{12}}{46080}-\frac{7 a_{13}}{46080}-\frac{a_{32}}{144}\nonumber \\ &-\frac{a_{33}}{288}-\frac{a_{34}}{96}+\frac{a_{35}}{288}-\frac{a_{36}}{288}-\frac{a_{37}}{72}+\frac{a_{38}}{96}-\frac{a_{39}}{288}+\frac{a_{40}}{144}\nonumber \\ &-\frac{a_{57}}{72}-\frac{a_{59}}{288}-\frac{a_{62}}{288}+\frac{a_{66}}{144}+\frac{a_{67}}{96}+\frac{a_{68}}{96}+\frac{a_{69}}{288}+\frac{a_{70}}{24}\nonumber \\ &-\frac{a_{72}}{13824}-\frac{a_{74}}{4608}-\frac{a_{76}}{13824}-\frac{a_{78}}{13824}+\frac{5 a_{79}}{13824}+\frac{a_{80}}{13824}+\frac{a_{81}}{13824}-\frac{a_{82}}{4608}\nonumber \\ &+\frac{a_{83}}{4608}+\frac{5 a_{84}}{6}+\frac{a_{86}}{6}+\frac{5 a_{92}}{8}+\frac{a_{93}}{6}+\frac{a_{94}}{8}+\frac{a_{95}}{12}-\frac{a_{97}}{4}\nonumber \\ &-\frac{a_{99}}{4}-\frac{a_{100}}{3}-\frac{a_{101}}{2}-\frac{2 a_{102}}{3}-\frac{5 a_{103}}{2}+144 a_{104}-\frac{1}{2304}.
 \end{align}
From the degree $1$ part of $\tau_{3,0}\Phi(\tau_{1,0}, \tau_{1,1})|_{t=0} = 0$, we obtain
\begin{align}
 0 &= -\frac{7 a_2}{144}-\frac{43 a_3}{720}-\frac{a_4}{10}+\frac{a_5}{72}-\frac{7 a_6}{720}-\frac{a_7}{20}+\frac{a_8}{480}\nonumber \\ &+\frac{a_{12}}{1920}+\frac{a_{13}}{1920}+\frac{7 a_{14}}{144}+\frac{7 a_{17}}{720}+\frac{7 a_{23}}{720}+\frac{7 a_{27}}{2880}+\frac{7 a_{28}}{2880}+\frac{7 a_{30}}{2880}\nonumber \\ &-\frac{a_{32}}{72}-\frac{5 a_{33}}{288}-\frac{a_{34}}{48}-\frac{a_{35}}{144}-\frac{a_{36}}{144}-\frac{5 a_{37}}{72}-\frac{a_{38}}{48}-\frac{5 a_{39}}{288}\nonumber \\ &-\frac{a_{40}}{72}+\frac{a_{41}}{288}+\frac{a_{43}}{72}+\frac{a_{45}}{288}+\frac{a_{52}}{72}+\frac{a_{54}}{288}+\frac{a_{56}}{288}-\frac{5 a_{57}}{72}\nonumber \\ &-\frac{a_{58}}{48}-\frac{5 a_{59}}{288}-\frac{a_{60}}{72}-\frac{a_{61}}{48}-\frac{5 a_{62}}{288}-\frac{a_{63}}{144}-\frac{a_{64}}{144}-\frac{a_{65}}{72}\nonumber \\ &+\frac{a_{67}}{288}+\frac{a_{70}}{72}+\frac{a_{72}}{6912}+\frac{a_{74}}{2304}+\frac{a_{76}}{6912}+\frac{a_{78}}{6912}+\frac{a_{79}}{768}+\frac{a_{80}}{2304}\nonumber \\ &+\frac{a_{81}}{1728}+\frac{a_{82}}{1152}+\frac{a_{83}}{1728}+5 a_{84}-\frac{2 a_{85}}{3}+a_{86}-\frac{a_{87}}{3}+\frac{35 a_{92}}{12}\nonumber \\ &+\frac{7 a_{94}}{12}+\frac{a_{97}}{2}+\frac{a_{99}}{2}-\frac{a_{100}}{3}-\frac{2 a_{101}}{3}-\frac{2 a_{102}}{3}-\frac{10 a_{103}}{3}+336 a_{104}\nonumber \\ &-\frac{313}{241920}.
 \end{align}
From the degree $1$ part of $\tau_{2,1}\Phi(\tau_{1,0}, \tau_{2,0})|_{t=0} = 0$, we obtain
\begin{align}
 0 &= \frac{a_2}{120}+\frac{3 a_5}{160}+\frac{a_6}{96}+\frac{7 a_8}{11520}-\frac{7 a_9}{138240}-\frac{7 a_{10}}{46080}+\frac{7 a_{12}}{46080}+\frac{7 a_{13}}{46080}\nonumber \\ &+\frac{a_{32}}{72}-\frac{a_{33}}{48}+\frac{a_{34}}{48}+\frac{a_{36}}{144}-\frac{5 a_{37}}{72}-\frac{a_{39}}{72}-\frac{5 a_{57}}{72}-\frac{a_{59}}{72}\nonumber \\ &-\frac{a_{62}}{48}-\frac{a_{66}}{72}-\frac{a_{68}}{48}-\frac{a_{69}}{144}+\frac{a_{70}}{48}+\frac{a_{72}}{13824}+\frac{a_{74}}{4608}-\frac{a_{75}}{13824}\nonumber \\ &+\frac{a_{76}}{13824}+\frac{a_{78}}{13824}+\frac{a_{79}}{1728}+\frac{a_{81}}{3456}+5 a_{84}+a_{86}+\frac{5 a_{92}}{2}+\frac{a_{94}}{2}\nonumber \\ &-\frac{a_{96}}{2}-\frac{a_{98}}{2}+\frac{2 a_{100}}{3}-a_{101}+\frac{4 a_{102}}{3}-5 a_{103}+432 a_{104}-\frac{79}{120960}.
 \end{align}
From the degree $1$ part of $\tau_{2,0}\Phi(\tau_{1,0}, \tau_{2,1})|_{t=0} = 0$, we obtain
\begin{align}
 0 &= -\frac{a_2}{48}+\frac{7 a_3}{480}+\frac{7 a_4}{240}+\frac{a_5}{96}-\frac{a_6}{240}+\frac{7 a_7}{480}-\frac{7 a_8}{11520}\nonumber \\ &+\frac{7 a_9}{138240}+\frac{7 a_{10}}{46080}-\frac{7 a_{12}}{46080}-\frac{7 a_{13}}{46080}-\frac{a_{32}}{144}-\frac{a_{33}}{288}-\frac{a_{34}}{96}+\frac{a_{35}}{288}\nonumber \\ &-\frac{a_{36}}{288}-\frac{7 a_{37}}{144}+\frac{a_{38}}{96}-\frac{a_{39}}{96}+\frac{a_{40}}{144}-\frac{5 a_{57}}{72}+\frac{a_{58}}{96}-\frac{5 a_{59}}{288}\nonumber \\ &+\frac{a_{60}}{144}-\frac{a_{61}}{96}-\frac{5 a_{62}}{288}-\frac{a_{63}}{288}+\frac{a_{64}}{288}-\frac{a_{65}}{144}-\frac{a_{66}}{144}-\frac{a_{67}}{96}\nonumber \\ &-\frac{a_{68}}{96}-\frac{a_{69}}{288}-\frac{a_{72}}{13824}-\frac{a_{74}}{4608}+\frac{a_{75}}{13824}-\frac{a_{76}}{13824}-\frac{a_{78}}{13824}+\frac{a_{79}}{3456}\nonumber \\ &+\frac{a_{80}}{1728}+\frac{a_{82}}{1152}-\frac{a_{83}}{6912}+5 a_{84}+a_{86}+\frac{5 a_{92}}{2}+\frac{a_{94}}{2}-\frac{a_{97}}{2}\nonumber \\ &-\frac{a_{99}}{2}-\frac{a_{101}}{3}-\frac{11 a_{103}}{3}+432 a_{104}+\frac{1}{5376}.
 \end{align}
From the degree $1$ part of $\tau_{3,1}\Phi(\tau_{1,1}, \tau_{0,0})|_{t=0} = 0$, we obtain
\begin{align}
 0 &= -\frac{31 a_1}{96768}+\frac{7 a_5}{720}+\frac{7 a_6}{720}-\frac{7 a_9}{34560}-\frac{7 a_{10}}{34560}+\frac{7 a_{12}}{138240}+\frac{7 a_{13}}{138240}\nonumber \\ &+\frac{7 a_{17}}{1440}-\frac{7 a_{21}}{34560}-\frac{7 a_{22}}{34560}+\frac{7 a_{27}}{2880}+\frac{7 a_{28}}{2880}+\frac{7 a_{30}}{2880}-\frac{7 a_{31}}{138240}-\frac{a_{33}}{288}\nonumber \\ &-\frac{a_{37}}{144}+\frac{a_{39}}{288}-\frac{a_{41}}{288}-\frac{a_{43}}{144}-\frac{a_{52}}{144}-\frac{a_{54}}{288}-\frac{a_{56}}{288}-\frac{a_{57}}{144}\nonumber \\ &+\frac{a_{67}}{288}+\frac{a_{70}}{72}+\frac{a_{72}}{13824}+\frac{a_{74}}{6912}-\frac{a_{75}}{13824}+\frac{a_{76}}{13824}+\frac{a_{78}}{13824}+\frac{a_{79}}{6912}\nonumber \\ &-\frac{a_{80}}{13824}+\frac{a_{81}}{13824}-\frac{a_{82}}{3456}+\frac{a_{83}}{13824}+\frac{a_{84}}{3}+\frac{a_{89}}{3}+\frac{a_{91}}{12}+\frac{a_{92}}{4}\nonumber \\ &+\frac{a_{94}}{24}-\frac{a_{95}}{24}-\frac{a_{98}}{12}-\frac{a_{101}}{6}-\frac{2 a_{103}}{3}+32 a_{104}-\frac{a_{105}}{13824}.
 \end{align}
From the degree $1$ part of $\tau_{4,0}\Phi(\tau_{1,1}, \tau_{0,0})|_{t=0} = 0$, we obtain
\begin{align}
 0 &= \frac{2329 a_1}{1451520}-\frac{a_2}{576}-\frac{a_5}{30}-\frac{a_6}{30}-\frac{a_8}{2304}+\frac{a_9}{1440}+\frac{a_{10}}{1440}-\frac{13 a_{12}}{46080}\nonumber \\ &-\frac{13 a_{13}}{46080}+\frac{7 a_{14}}{576}+\frac{7 a_{15}}{576}+\frac{7 a_{16}}{576}-\frac{11 a_{17}}{720}+\frac{7 a_{18}}{576}+\frac{7 a_{19}}{576}+\frac{a_{21}}{1440}\nonumber \\ &+\frac{a_{22}}{1440}+\frac{a_{23}}{192}+\frac{7 a_{24}}{2880}+\frac{7 a_{25}}{2880}-\frac{a_{27}}{120}-\frac{a_{28}}{120}+\frac{7 a_{29}}{2880}-\frac{a_{30}}{120}\nonumber \\ &+\frac{a_{31}}{5760}-\frac{a_{33}}{144}-\frac{7 a_{37}}{288}-\frac{a_{39}}{144}-\frac{a_{41}}{144}-\frac{7 a_{43}}{288}+\frac{a_{44}}{288}-\frac{a_{45}}{144}\nonumber \\ &+\frac{a_{50}}{288}-\frac{7 a_{52}}{288}+\frac{a_{53}}{288}-\frac{a_{54}}{144}-\frac{a_{56}}{144}-\frac{7 a_{57}}{288}-\frac{a_{59}}{144}-\frac{a_{62}}{144}\nonumber \\ &+\frac{a_{70}}{288}+\frac{a_{72}}{13824}+\frac{5 a_{74}}{13824}+\frac{a_{75}}{6912}+\frac{a_{76}}{13824}+\frac{a_{78}}{13824}+\frac{5 a_{79}}{13824}+\frac{a_{80}}{13824}\nonumber \\ &+\frac{a_{81}}{13824}+\frac{a_{82}}{4608}+\frac{a_{83}}{13824}+\frac{3 a_{84}}{2}-\frac{a_{86}}{6}+\frac{3 a_{89}}{2}+\frac{a_{91}}{6}+\frac{11 a_{92}}{12}\nonumber \\ &+\frac{a_{95}}{12}+\frac{a_{98}}{6}-\frac{a_{101}}{6}-\frac{5 a_{103}}{6}+72 a_{104}+\frac{a_{105}}{6912}-\frac{13}{8064}.
 \end{align}
From the degree $1$ part of $\tau_{2,1}\Phi(\tau_{1,1}, \tau_{0,1})|_{t=0} = 0$, we obtain
\begin{align}
 0 &= \frac{7 a_7}{960}-\frac{7 a_{12}}{46080}-\frac{7 a_{13}}{46080}-\frac{a_{33}}{288}-\frac{a_{37}}{288}+\frac{a_{40}}{288}-\frac{a_{41}}{576}-\frac{a_{43}}{576}\nonumber \\ &-\frac{a_{52}}{288}-\frac{a_{57}}{288}-\frac{a_{59}}{576}+\frac{a_{60}}{576}-\frac{a_{62}}{576}+\frac{a_{64}}{576}+\frac{a_{67}}{288}+\frac{a_{70}}{96}\nonumber \\ &+\frac{a_{72}}{13824}+\frac{a_{74}}{13824}-\frac{a_{76}}{13824}-\frac{a_{78}}{13824}+\frac{a_{79}}{13824}-\frac{a_{80}}{13824}+\frac{a_{81}}{13824}-\frac{a_{82}}{4608}\nonumber \\ &+\frac{a_{83}}{13824}+\frac{a_{84}}{6}+\frac{a_{89}}{12}+\frac{a_{92}}{8}+\frac{a_{94}}{24}-\frac{a_{99}}{12}-\frac{a_{101}}{6}-\frac{a_{103}}{2}\nonumber \\ &+24 a_{104}.
 \end{align}
From the degree $1$ part of $\tau_{3,0}\Phi(\tau_{1,1}, \tau_{0,1})|_{t=0} = 0$, we obtain
\begin{align}
 0 &= -\frac{31 a_1}{96768}+\frac{7 a_5}{720}+\frac{7 a_6}{720}-\frac{a_7}{40}+\frac{7 a_8}{34560}-\frac{7 a_9}{34560}-\frac{7 a_{10}}{34560}\nonumber \\ &+\frac{a_{12}}{1920}+\frac{a_{13}}{1920}+\frac{7 a_{14}}{288}+\frac{7 a_{17}}{720}-\frac{7 a_{21}}{34560}-\frac{7 a_{22}}{34560}+\frac{7 a_{23}}{720}+\frac{7 a_{27}}{2880}\nonumber \\ &+\frac{7 a_{28}}{2880}+\frac{7 a_{30}}{2880}-\frac{7 a_{31}}{138240}-\frac{a_{33}}{96}-\frac{a_{37}}{48}+\frac{a_{39}}{288}-\frac{a_{40}}{144}-\frac{a_{41}}{288}\nonumber \\ &-\frac{a_{43}}{288}+\frac{a_{45}}{288}-\frac{a_{52}}{48}+\frac{a_{54}}{288}+\frac{a_{56}}{288}-\frac{a_{57}}{48}-\frac{a_{59}}{288}-\frac{a_{60}}{288}\nonumber \\ &-\frac{a_{62}}{288}-\frac{a_{64}}{288}+\frac{a_{67}}{288}+\frac{a_{70}}{72}+\frac{a_{72}}{6912}+\frac{5 a_{74}}{13824}-\frac{a_{75}}{13824}+\frac{a_{76}}{6912}\nonumber \\ &+\frac{a_{78}}{6912}+\frac{5 a_{79}}{13824}+\frac{a_{81}}{6912}-\frac{a_{82}}{13824}+\frac{a_{83}}{6912}+\frac{7 a_{84}}{6}-\frac{a_{87}}{6}+\frac{7 a_{89}}{12}\nonumber \\ &+\frac{7 a_{92}}{8}-\frac{a_{93}}{6}+\frac{a_{94}}{6}-\frac{a_{95}}{8}-\frac{a_{98}}{12}+\frac{a_{99}}{6}-\frac{a_{101}}{3}-\frac{7 a_{103}}{6}\nonumber \\ &+80 a_{104}-\frac{a_{105}}{13824}-\frac{31}{48384}.
 \end{align}
From the degree $1$ part of $\tau_{2,1}\Phi(\tau_{1,1}, \tau_{1,0})|_{t=0} = 0$, we obtain
\begin{align}
 0 &= \frac{a_3}{288}+\frac{7 a_5}{480}+\frac{7 a_6}{240}+\frac{7 a_7}{480}-\frac{7 a_9}{23040}-\frac{7 a_{10}}{11520}-\frac{7 a_{11}}{46080}-\frac{7 a_{13}}{46080}\nonumber \\ &-\frac{a_{32}}{144}-\frac{a_{33}}{96}-\frac{a_{34}}{288}-\frac{a_{35}}{288}+\frac{a_{36}}{288}-\frac{a_{37}}{72}-\frac{a_{38}}{288}+\frac{a_{39}}{96}\nonumber \\ &+\frac{a_{40}}{144}-\frac{a_{57}}{72}-\frac{a_{58}}{288}-\frac{a_{61}}{288}+\frac{a_{66}}{144}+\frac{a_{67}}{96}+\frac{a_{68}}{96}+\frac{a_{69}}{288}\nonumber \\ &+\frac{a_{70}}{24}-\frac{a_{73}}{13824}-\frac{a_{75}}{4608}-\frac{a_{76}}{13824}-\frac{a_{77}}{13824}+\frac{5 a_{79}}{13824}+\frac{a_{80}}{13824}+\frac{a_{81}}{13824}\nonumber \\ &-\frac{a_{82}}{4608}+\frac{a_{83}}{4608}+\frac{5 a_{84}}{6}+\frac{a_{85}}{6}+\frac{5 a_{92}}{8}+\frac{a_{93}}{8}+\frac{a_{94}}{6}+\frac{a_{95}}{12}\nonumber \\ &-\frac{a_{98}}{4}-\frac{a_{99}}{4}-\frac{a_{100}}{3}-\frac{2 a_{101}}{3}-\frac{a_{102}}{2}-\frac{5 a_{103}}{2}+144 a_{104}-\frac{1}{4608}.
 \end{align}
From the degree $1$ part of $\tau_{2,0}\Phi(\tau_{1,1}, \tau_{1,1})|_{t=0} = 0$, we obtain
\begin{align}
 0 &= \frac{7 a_3}{480}+\frac{7 a_4}{240}+\frac{7 a_5}{480}+\frac{7 a_6}{240}+\frac{7 a_7}{480}-\frac{7 a_8}{11520}-\frac{7 a_9}{23040}-\frac{7 a_{10}}{11520}\nonumber \\ &-\frac{7 a_{11}}{46080}-\frac{7 a_{12}}{46080}-\frac{7 a_{13}}{23040}-\frac{a_{33}}{96}-\frac{a_{34}}{96}+\frac{a_{35}}{288}+\frac{a_{36}}{288}-\frac{a_{37}}{48}\nonumber \\ &+\frac{a_{38}}{96}+\frac{a_{39}}{96}+\frac{a_{40}}{144}-\frac{a_{57}}{48}+\frac{a_{60}}{144}+\frac{a_{63}}{288}+\frac{a_{64}}{288}+\frac{a_{67}}{288}\nonumber \\ &+\frac{a_{68}}{288}+\frac{a_{69}}{288}+\frac{a_{70}}{72}-\frac{a_{72}}{13824}-\frac{a_{73}}{13824}-\frac{a_{74}}{4608}-\frac{a_{75}}{4608}-\frac{a_{76}}{6912}\nonumber \\ &-\frac{a_{77}}{13824}-\frac{a_{78}}{13824}+\frac{a_{79}}{2304}+\frac{a_{80}}{2304}+\frac{a_{81}}{2304}+\frac{a_{82}}{2304}+\frac{a_{83}}{1728}+a_{84}\nonumber \\ &+\frac{3 a_{92}}{4}+\frac{a_{93}}{8}+\frac{a_{94}}{8}-\frac{a_{97}}{4}-\frac{a_{98}}{4}-\frac{a_{99}}{3}-\frac{5 a_{101}}{6}-\frac{5 a_{102}}{6}\nonumber \\ &-\frac{7 a_{103}}{3}+192 a_{104}-\frac{1}{4608}.
 \end{align}
From the degree $1$ part of $\tau_{2,0}\Phi(\tau_{1,1}, \tau_{2,0})|_{t=0} = 0$, we obtain
\begin{align}
 0 &= -\frac{a_2}{96}-\frac{11 a_5}{192}-\frac{379 a_6}{2880}+\frac{7 a_7}{960}+\frac{7 a_8}{11520}+\frac{193 a_9}{138240}+\frac{463 a_{10}}{138240}\nonumber \\ &+\frac{7 a_{11}}{23040}+\frac{7 a_{12}}{46080}+\frac{7 a_{13}}{46080}+\frac{7 a_{15}}{288}+\frac{7 a_{16}}{96}+\frac{7 a_{24}}{480}-\frac{a_{32}}{96}-\frac{11 a_{33}}{288}\nonumber \\ &+\frac{a_{34}}{48}-\frac{a_{36}}{144}-\frac{5 a_{37}}{72}-\frac{a_{39}}{32}+\frac{a_{40}}{288}+\frac{a_{44}}{48}+\frac{a_{47}}{288}+\frac{a_{53}}{48}\nonumber \\ &-\frac{13 a_{57}}{144}-\frac{a_{59}}{32}+\frac{a_{60}}{288}-\frac{11 a_{62}}{288}-\frac{a_{63}}{144}-\frac{a_{65}}{96}-\frac{a_{66}}{96}-\frac{5 a_{67}}{288}\nonumber \\ &-\frac{a_{68}}{144}-\frac{a_{69}}{144}-\frac{a_{70}}{48}+\frac{a_{72}}{13824}+\frac{a_{73}}{6912}+\frac{a_{74}}{4608}+\frac{13 a_{75}}{13824}+\frac{a_{76}}{13824}\nonumber \\ &+\frac{a_{77}}{6912}+\frac{a_{78}}{13824}+\frac{a_{79}}{6912}+\frac{a_{80}}{6912}-\frac{a_{81}}{3456}+\frac{a_{82}}{1152}-\frac{a_{83}}{1728}+\frac{19 a_{84}}{3}\nonumber \\ &-a_{86}+\frac{7 a_{92}}{2}-\frac{a_{94}}{6}+\frac{a_{98}}{6}-\frac{a_{99}}{2}-a_{101}+\frac{5 a_{102}}{3}-\frac{10 a_{103}}{3}\nonumber \\ &+528 a_{104}-\frac{59}{16128}.
 \end{align}
From the degree $1$ part of $\tau_{3,1}\Phi(\tau_{2,0}, \tau_{0,0})|_{t=0} = 0$, we obtain
\begin{align}
 0 &= \frac{31 a_1}{96768}-\frac{a_2}{288}+\frac{7 a_3}{2880}+\frac{7 a_4}{960}-\frac{7 a_7}{2880}-\frac{a_8}{1152}+\frac{7 a_9}{34560}+\frac{7 a_{10}}{34560}\nonumber \\ &-\frac{7 a_{12}}{69120}+\frac{7 a_{14}}{288}+\frac{7 a_{21}}{34560}+\frac{7 a_{22}}{34560}+\frac{a_{23}}{96}+\frac{7 a_{31}}{138240}-\frac{a_{32}}{288}+\frac{a_{33}}{144}\nonumber \\ &-\frac{a_{34}}{288}-\frac{a_{37}}{36}+\frac{a_{38}}{288}-\frac{a_{40}}{288}+\frac{a_{41}}{144}-\frac{a_{43}}{36}-\frac{a_{46}}{288}-\frac{a_{52}}{36}\nonumber \\ &+\frac{a_{56}}{144}-\frac{a_{57}}{36}+\frac{a_{58}}{288}-\frac{a_{60}}{288}-\frac{a_{61}}{288}-\frac{a_{65}}{288}+\frac{a_{66}}{288}-\frac{a_{67}}{144}\nonumber \\ &-\frac{a_{68}}{288}+\frac{a_{70}}{48}-\frac{a_{72}}{6912}+\frac{7 a_{74}}{13824}+\frac{a_{75}}{13824}-\frac{a_{78}}{6912}+\frac{7 a_{79}}{13824}+\frac{a_{80}}{6912}\nonumber \\ &+\frac{a_{82}}{13824}-\frac{a_{83}}{6912}+\frac{5 a_{84}}{3}+\frac{a_{87}}{6}+\frac{5 a_{89}}{3}+\frac{5 a_{92}}{6}+\frac{a_{95}}{12}-\frac{a_{97}}{6}\nonumber \\ &-\frac{a_{100}}{6}+\frac{a_{101}}{3}-\frac{5 a_{103}}{3}+112 a_{104}+\frac{a_{105}}{13824}-\frac{1}{1536}.
 \end{align}
From the degree $0$ part of $\tau_{2,1}\tau_{2,0}\Phi(\tau_{0,0}, \tau_{0,0})|_{t=0} = 0$, we obtain
\begin{align}
 0 &= -\frac{31 a_1}{10752}-\frac{7 a_8}{4608}-\frac{7 a_9}{4608}-\frac{7 a_{10}}{4608}-\frac{7 a_{11}}{23040}-\frac{7 a_{12}}{23040}-\frac{7 a_{13}}{23040}\nonumber \\ &-\frac{7 a_{20}}{23040}-\frac{7 a_{21}}{4608}-\frac{7 a_{22}}{4608}-\frac{7 a_{31}}{23040}-\frac{a_{71}}{6912}-\frac{a_{72}}{6912}-\frac{a_{73}}{6912}-\frac{a_{74}}{2304}\nonumber \\ &-\frac{a_{75}}{2304}-\frac{a_{76}}{6912}-\frac{a_{77}}{6912}-\frac{a_{78}}{6912}-\frac{a_{79}}{2304}-\frac{a_{80}}{6912}-\frac{a_{81}}{6912}-\frac{a_{82}}{2304}\nonumber \\ &-\frac{a_{83}}{6912}-\frac{a_{105}}{2304}+\frac{31}{1536}.
 \end{align}
From the degree $0$ part of $\tau_{2,0}\tau_{3,0}\Phi(\tau_{0,0}, \tau_{0,0})|_{t=0} = 0$, we obtain
\begin{align}
 0 &= \frac{859 a_1}{48384}+\frac{7 a_2}{48}+\frac{7 a_3}{720}+\frac{7 a_4}{720}+\frac{7 a_5}{720}+\frac{7 a_6}{720}+\frac{7 a_7}{960}+\frac{229 a_8}{34560}\nonumber \\ &+\frac{229 a_9}{34560}+\frac{229 a_{10}}{34560}+\frac{19 a_{11}}{23040}+\frac{19 a_{12}}{23040}+\frac{a_{13}}{1920}+\frac{7 a_{14}}{48}+\frac{7 a_{15}}{48}+\frac{7 a_{16}}{48}\nonumber \\ &+\frac{7 a_{17}}{720}+\frac{7 a_{18}}{48}+\frac{7 a_{19}}{48}+\frac{19 a_{20}}{23040}+\frac{229 a_{21}}{34560}+\frac{229 a_{22}}{34560}+\frac{7 a_{23}}{288}+\frac{7 a_{24}}{288}\nonumber \\ &+\frac{7 a_{25}}{288}+\frac{7 a_{26}}{2880}+\frac{7 a_{27}}{2880}+\frac{7 a_{28}}{2880}+\frac{7 a_{29}}{288}+\frac{7 a_{30}}{2880}+\frac{31 a_{31}}{23040}+\frac{a_{32}}{288}\nonumber \\ &+\frac{a_{33}}{288}+\frac{a_{34}}{288}+\frac{5 a_{37}}{144}+\frac{a_{38}}{288}+\frac{a_{39}}{288}+\frac{a_{40}}{288}+\frac{a_{41}}{288}+\frac{a_{42}}{288}\nonumber \\ &+\frac{5 a_{43}}{144}+\frac{5 a_{44}}{144}+\frac{a_{45}}{288}+\frac{a_{46}}{288}+\frac{a_{47}}{288}+\frac{a_{49}}{288}+\frac{5 a_{50}}{144}+\frac{a_{51}}{288}\nonumber \\ &+\frac{5 a_{52}}{144}+\frac{5 a_{53}}{144}+\frac{a_{54}}{288}+\frac{a_{55}}{288}+\frac{a_{56}}{288}+\frac{5 a_{57}}{144}+\frac{a_{58}}{288}+\frac{a_{59}}{288}\nonumber \\ &+\frac{a_{60}}{288}+\frac{a_{61}}{288}+\frac{a_{62}}{288}+\frac{a_{65}}{288}+\frac{a_{66}}{288}+\frac{a_{67}}{288}+\frac{a_{68}}{288}+\frac{5 a_{70}}{144}\nonumber \\ &+\frac{a_{71}}{3456}+\frac{a_{72}}{3456}+\frac{a_{73}}{3456}+\frac{a_{74}}{768}+\frac{a_{75}}{768}+\frac{a_{76}}{6912}+\frac{a_{77}}{3456}+\frac{a_{78}}{3456}\nonumber \\ &+\frac{a_{79}}{768}+\frac{a_{80}}{6912}+\frac{a_{81}}{6912}+\frac{a_{82}}{2304}+\frac{a_{83}}{3456}+\frac{a_{105}}{768}-\frac{859}{6912}.
 \end{align}
From the degree $1$ part of $\tau_{3,1}\tau_{3,0}\Phi(\tau_{0,0}, \tau_{0,0})|_{t=0} = 0$, we obtain
\begin{align}
 0 &= \frac{53 a_1}{48384}+\frac{a_2}{48}+\frac{7 a_3}{720}+\frac{7 a_4}{720}+\frac{7 a_5}{720}+\frac{7 a_6}{720}-\frac{79 a_7}{2880}-\frac{61 a_8}{34560}\nonumber \\ &-\frac{61 a_9}{34560}-\frac{61 a_{10}}{34560}-\frac{43 a_{11}}{69120}-\frac{43 a_{12}}{69120}+\frac{a_{13}}{1920}+\frac{a_{14}}{12}+\frac{7 a_{15}}{48}+\frac{a_{16}}{12}\nonumber \\ &+\frac{7 a_{17}}{720}+\frac{a_{18}}{12}+\frac{7 a_{19}}{48}-\frac{43 a_{20}}{69120}-\frac{61 a_{21}}{34560}-\frac{61 a_{22}}{34560}+\frac{7 a_{23}}{288}+\frac{7 a_{24}}{288}\nonumber \\ &+\frac{7 a_{25}}{288}+\frac{7 a_{26}}{2880}+\frac{7 a_{27}}{2880}+\frac{7 a_{28}}{2880}+\frac{7 a_{29}}{288}+\frac{7 a_{30}}{2880}-\frac{23 a_{31}}{69120}-\frac{a_{32}}{96}\nonumber \\ &+\frac{a_{33}}{288}+\frac{a_{34}}{288}-\frac{a_{37}}{16}+\frac{a_{38}}{288}+\frac{a_{39}}{288}-\frac{a_{40}}{96}+\frac{a_{41}}{288}+\frac{a_{42}}{288}\nonumber \\ &-\frac{a_{43}}{16}-\frac{a_{44}}{16}+\frac{a_{45}}{288}-\frac{a_{46}}{288}-\frac{a_{47}}{288}-\frac{a_{49}}{288}-\frac{a_{50}}{16}+\frac{a_{51}}{288}\nonumber \\ &-\frac{a_{52}}{16}-\frac{a_{53}}{16}+\frac{a_{54}}{288}+\frac{a_{55}}{288}+\frac{a_{56}}{288}-\frac{a_{57}}{16}+\frac{a_{58}}{288}+\frac{a_{59}}{288}\nonumber \\ &-\frac{a_{60}}{96}-\frac{a_{61}}{288}-\frac{a_{62}}{288}-\frac{a_{65}}{96}+\frac{a_{66}}{288}-\frac{a_{67}}{96}-\frac{a_{68}}{96}+\frac{5 a_{70}}{144}\nonumber \\ &-\frac{a_{71}}{3456}-\frac{a_{72}}{3456}-\frac{a_{73}}{3456}+\frac{a_{74}}{768}+\frac{a_{75}}{768}+\frac{a_{76}}{6912}-\frac{a_{77}}{3456}-\frac{a_{78}}{3456}\nonumber \\ &+\frac{a_{79}}{768}+\frac{a_{80}}{6912}+\frac{a_{81}}{6912}+\frac{a_{82}}{2304}-\frac{a_{83}}{3456}+\frac{11 a_{84}}{3}+\frac{11 a_{88}}{3}+\frac{11 a_{89}}{3}\nonumber \\ &+\frac{11 a_{90}}{3}+\frac{23 a_{92}}{12}+\frac{a_{95}}{12}-\frac{a_{97}}{6}-\frac{a_{98}}{6}+\frac{a_{99}}{6}-\frac{a_{100}}{3}+\frac{a_{101}}{3}\nonumber \\ &+\frac{a_{102}}{3}-\frac{10 a_{103}}{3}+240 a_{104}+\frac{a_{105}}{768}-\frac{53}{6912}.
 \end{align}
From the degree $1$ part of $\tau_{2,1}\tau_{3,1}\Phi(\tau_{0,0}, \tau_{0,0})|_{t=0} = 0$, we obtain
\begin{align}
 0 &= -\frac{5 a_1}{48384}+\frac{a_2}{144}+\frac{7 a_3}{720}+\frac{7 a_4}{720}+\frac{7 a_5}{720}+\frac{7 a_6}{720}+\frac{7 a_7}{960}\nonumber \\ &-\frac{5 a_8}{13824}-\frac{5 a_9}{13824}-\frac{5 a_{10}}{13824}+\frac{7 a_{11}}{34560}+\frac{7 a_{12}}{34560}-\frac{7 a_{13}}{69120}+\frac{a_{14}}{72}+\frac{a_{15}}{48}\nonumber \\ &+\frac{a_{16}}{72}+\frac{7 a_{17}}{1440}+\frac{a_{18}}{72}+\frac{a_{19}}{48}+\frac{7 a_{20}}{34560}-\frac{5 a_{21}}{13824}-\frac{5 a_{22}}{13824}+\frac{a_{23}}{288}\nonumber \\ &+\frac{a_{24}}{288}+\frac{a_{25}}{288}+\frac{7 a_{26}}{2880}+\frac{7 a_{27}}{2880}+\frac{7 a_{28}}{2880}+\frac{a_{29}}{288}+\frac{7 a_{30}}{2880}-\frac{17 a_{31}}{138240}\nonumber \\ &-\frac{a_{32}}{288}-\frac{a_{33}}{288}-\frac{a_{34}}{288}-\frac{a_{37}}{72}+\frac{a_{38}}{288}+\frac{a_{39}}{288}+\frac{a_{40}}{288}-\frac{a_{41}}{288}\nonumber \\ &-\frac{a_{42}}{288}-\frac{a_{43}}{72}-\frac{a_{44}}{72}-\frac{a_{46}}{288}-\frac{a_{47}}{288}-\frac{a_{49}}{288}-\frac{a_{50}}{72}-\frac{a_{51}}{288}\nonumber \\ &-\frac{a_{52}}{72}-\frac{a_{53}}{72}-\frac{a_{54}}{288}-\frac{a_{55}}{288}-\frac{a_{56}}{288}-\frac{a_{57}}{72}+\frac{a_{66}}{288}+\frac{a_{67}}{288}\nonumber \\ &+\frac{a_{68}}{288}+\frac{5 a_{70}}{144}+\frac{a_{71}}{6912}+\frac{a_{72}}{6912}+\frac{a_{73}}{6912}+\frac{5 a_{74}}{13824}+\frac{5 a_{75}}{13824}+\frac{a_{77}}{6912}\nonumber \\ &+\frac{a_{78}}{6912}+\frac{5 a_{79}}{13824}-\frac{7 a_{82}}{13824}+\frac{a_{83}}{6912}+\frac{2 a_{84}}{3}+\frac{2 a_{88}}{3}+\frac{2 a_{89}}{3}+\frac{2 a_{90}}{3}\nonumber \\ &+\frac{a_{91}}{12}+\frac{a_{92}}{2}+\frac{a_{93}}{24}+\frac{a_{94}}{24}-\frac{a_{97}}{12}-\frac{a_{98}}{12}-\frac{a_{99}}{12}-\frac{a_{100}}{6}\nonumber \\ &-\frac{a_{101}}{6}-\frac{a_{102}}{6}-\frac{5 a_{103}}{3}+80 a_{104}+\frac{5 a_{105}}{13824}-\frac{7}{4608}.
 \end{align}
From the degree $1$ part of $\tau_{2,1}\tau_{4,0}\Phi(\tau_{0,0}, \tau_{0,0})|_{t=0} = 0$, we obtain
\begin{align}
 0 &= \frac{923 a_1}{290304}+\frac{17 a_2}{576}-\frac{a_3}{30}-\frac{a_4}{30}-\frac{a_5}{30}-\frac{a_6}{30}-\frac{109 a_8}{69120}-\frac{109 a_9}{69120}\nonumber \\ &-\frac{109 a_{10}}{69120}-\frac{a_{11}}{7680}-\frac{a_{12}}{7680}-\frac{a_{13}}{7680}+\frac{5 a_{14}}{64}+\frac{73 a_{15}}{576}+\frac{5 a_{16}}{64}-\frac{11 a_{17}}{720}\nonumber \\ &+\frac{5 a_{18}}{64}+\frac{73 a_{19}}{576}-\frac{a_{20}}{7680}-\frac{109 a_{21}}{69120}-\frac{109 a_{22}}{69120}+\frac{23 a_{23}}{960}+\frac{23 a_{24}}{960}+\frac{23 a_{25}}{960}\nonumber \\ &-\frac{a_{26}}{120}-\frac{a_{27}}{120}-\frac{a_{28}}{120}+\frac{23 a_{29}}{960}-\frac{a_{30}}{120}-\frac{7 a_{31}}{46080}+\frac{a_{32}}{144}-\frac{a_{33}}{144}\nonumber \\ &-\frac{a_{34}}{144}-\frac{17 a_{37}}{288}-\frac{a_{38}}{144}-\frac{a_{39}}{144}-\frac{a_{41}}{144}-\frac{a_{42}}{144}-\frac{17 a_{43}}{288}-\frac{17 a_{44}}{288}\nonumber \\ &-\frac{a_{45}}{144}-\frac{17 a_{50}}{288}-\frac{a_{51}}{144}-\frac{17 a_{52}}{288}-\frac{17 a_{53}}{288}-\frac{a_{54}}{144}-\frac{a_{55}}{144}-\frac{a_{56}}{144}\nonumber \\ &-\frac{17 a_{57}}{288}-\frac{a_{58}}{144}-\frac{a_{59}}{144}-\frac{a_{61}}{144}-\frac{a_{62}}{144}-\frac{a_{66}}{144}+\frac{a_{70}}{96}+\frac{a_{71}}{6912}\nonumber \\ &+\frac{a_{72}}{6912}+\frac{a_{73}}{6912}+\frac{11 a_{74}}{13824}+\frac{11 a_{75}}{13824}+\frac{a_{76}}{6912}+\frac{a_{77}}{6912}+\frac{a_{78}}{6912}+\frac{11 a_{79}}{13824}\nonumber \\ &+\frac{a_{80}}{6912}+\frac{a_{81}}{6912}+\frac{5 a_{82}}{13824}+\frac{a_{83}}{6912}+\frac{7 a_{84}}{2}-\frac{a_{85}}{6}-\frac{a_{86}}{6}+\frac{7 a_{88}}{2}\nonumber \\ &+\frac{7 a_{89}}{2}+\frac{7 a_{90}}{2}+\frac{a_{91}}{6}+2 a_{92}+\frac{a_{95}}{12}-\frac{a_{96}}{6}+\frac{a_{97}}{6}+\frac{a_{98}}{6}\nonumber \\ &+\frac{a_{100}}{3}-\frac{a_{101}}{6}-\frac{a_{102}}{6}-\frac{5 a_{103}}{2}+200 a_{104}+\frac{11 a_{105}}{13824}-\frac{5489}{483840}.
 \end{align}
From the degree $1$ part of $\tau_{2,0}\tau_{4,1}\Phi(\tau_{0,0}, \tau_{0,0})|_{t=0} = 0$, we obtain
\begin{align}
 0 &= \frac{139 a_1}{241920}+\frac{a_2}{240}+\frac{a_3}{240}-\frac{a_4}{720}+\frac{a_5}{240}-\frac{a_6}{720}+\frac{7 a_7}{960}-\frac{7 a_8}{4608}\nonumber \\ &-\frac{89 a_9}{69120}-\frac{7 a_{10}}{4608}+\frac{a_{11}}{11520}+\frac{a_{12}}{11520}-\frac{a_{13}}{4608}+\frac{11 a_{14}}{240}+\frac{7 a_{15}}{80}+\frac{11 a_{16}}{240}\nonumber \\ &+\frac{a_{17}}{720}+\frac{11 a_{18}}{240}+\frac{7 a_{19}}{80}+\frac{a_{20}}{11520}-\frac{89 a_{21}}{69120}-\frac{7 a_{22}}{4608}+\frac{7 a_{23}}{480}+\frac{7 a_{24}}{480}\nonumber \\ &+\frac{7 a_{25}}{480}+\frac{a_{26}}{960}+\frac{a_{27}}{960}+\frac{a_{28}}{960}+\frac{7 a_{29}}{480}+\frac{a_{30}}{960}-\frac{11 a_{31}}{46080}-\frac{a_{32}}{288}\nonumber \\ &-\frac{a_{33}}{288}-\frac{a_{34}}{288}-\frac{5 a_{37}}{144}-\frac{a_{38}}{288}-\frac{a_{39}}{288}+\frac{a_{40}}{288}-\frac{a_{41}}{288}-\frac{a_{42}}{288}\nonumber \\ &-\frac{5 a_{43}}{144}-\frac{5 a_{44}}{144}-\frac{a_{45}}{288}-\frac{a_{46}}{288}-\frac{a_{47}}{288}-\frac{a_{49}}{288}-\frac{5 a_{50}}{144}-\frac{a_{51}}{288}\nonumber \\ &-\frac{5 a_{52}}{144}-\frac{5 a_{53}}{144}-\frac{a_{54}}{288}-\frac{a_{55}}{288}-\frac{a_{56}}{288}-\frac{a_{57}}{24}-\frac{a_{58}}{288}-\frac{a_{59}}{288}\nonumber \\ &+\frac{a_{60}}{288}-\frac{a_{61}}{288}-\frac{a_{62}}{288}-\frac{a_{65}}{288}-\frac{a_{66}}{288}+\frac{a_{67}}{288}+\frac{a_{68}}{288}+\frac{a_{70}}{48}\nonumber \\ &+\frac{a_{71}}{6912}+\frac{a_{72}}{6912}+\frac{a_{73}}{6912}+\frac{5 a_{74}}{13824}+\frac{5 a_{75}}{13824}+\frac{a_{77}}{6912}+\frac{a_{78}}{6912}+\frac{5 a_{79}}{13824}\nonumber \\ &+\frac{a_{82}}{13824}+\frac{a_{83}}{6912}+\frac{5 a_{84}}{2}+\frac{a_{85}}{6}+\frac{a_{86}}{6}+\frac{5 a_{88}}{2}+\frac{5 a_{89}}{2}+\frac{5 a_{90}}{2}\nonumber \\ &+\frac{a_{91}}{6}+\frac{5 a_{92}}{4}+\frac{a_{93}}{12}+\frac{a_{94}}{12}-\frac{a_{99}}{6}-\frac{a_{101}}{6}-\frac{a_{102}}{6}-\frac{13 a_{103}}{6}\nonumber \\ &+160 a_{104}+\frac{5 a_{105}}{13824}-\frac{199}{53760}.
 \end{align}
From the degree $1$ part of $\tau_{2,1}\tau_{2,1}\Phi(\tau_{0,0}, \tau_{0,1})|_{t=0} = 0$, we obtain
\begin{align}
 0 &= \frac{a_1}{13824}+\frac{a_2}{288}+\frac{7 a_7}{480}-\frac{a_8}{13824}-\frac{a_9}{13824}-\frac{5 a_{10}}{13824}+\frac{7 a_{11}}{23040}-\frac{7 a_{12}}{23040}\nonumber \\ &-\frac{7 a_{13}}{23040}+\frac{a_{14}}{576}+\frac{a_{15}}{288}+\frac{7 a_{16}}{576}+\frac{7 a_{18}}{576}+\frac{a_{19}}{288}+\frac{7 a_{20}}{23040}-\frac{a_{21}}{13824}\nonumber \\ &-\frac{5 a_{22}}{13824}+\frac{a_{24}}{288}+\frac{a_{25}}{288}+\frac{a_{29}}{576}-\frac{a_{31}}{13824}-\frac{a_{33}}{144}-\frac{a_{37}}{144}+\frac{a_{40}}{144}\nonumber \\ &-\frac{a_{41}}{288}-\frac{a_{43}}{288}-\frac{a_{44}}{144}-\frac{a_{47}}{288}-\frac{a_{49}}{288}-\frac{a_{50}}{144}-\frac{a_{52}}{144}-\frac{a_{53}}{144}\nonumber \\ &-\frac{a_{57}}{144}-\frac{a_{59}}{288}+\frac{a_{60}}{288}-\frac{a_{62}}{288}+\frac{a_{64}}{288}+\frac{a_{67}}{144}+\frac{a_{70}}{48}+\frac{a_{71}}{6912}\nonumber \\ &+\frac{a_{72}}{6912}+\frac{a_{73}}{6912}+\frac{a_{74}}{6912}+\frac{a_{75}}{6912}-\frac{a_{76}}{6912}+\frac{a_{77}}{6912}-\frac{a_{78}}{6912}+\frac{a_{79}}{6912}\nonumber \\ &-\frac{a_{80}}{6912}+\frac{a_{81}}{6912}-\frac{a_{82}}{2304}+\frac{a_{83}}{6912}+\frac{a_{84}}{3}+\frac{a_{88}}{3}+\frac{a_{89}}{6}+\frac{a_{90}}{3}\nonumber \\ &+\frac{a_{92}}{4}+\frac{a_{94}}{12}-\frac{a_{99}}{6}-\frac{a_{101}}{3}-a_{103}+48 a_{104}+\frac{a_{105}}{6912}-\frac{1}{4608}.
 \end{align}
From the degree $1$ part of $\tau_{2,1}\tau_{2,1}\Phi(\tau_{0,0}, \tau_{1,0})|_{t=0} = 0$, we obtain
\begin{align}
 0 &= \frac{a_2}{48}+\frac{5 a_3}{288}+\frac{a_4}{288}+\frac{7 a_5}{240}+\frac{7 a_6}{120}+\frac{7 a_7}{160}-\frac{5 a_9}{6912}-\frac{5 a_{10}}{3456}\nonumber \\ &+\frac{a_{11}}{69120}-\frac{41 a_{13}}{69120}+\frac{11 a_{15}}{288}+\frac{a_{16}}{16}+\frac{a_{17}}{96}+\frac{5 a_{24}}{288}+\frac{a_{28}}{288}+\frac{a_{30}}{288}\nonumber \\ &-\frac{a_{32}}{48}-\frac{a_{33}}{48}-\frac{a_{34}}{144}-\frac{a_{35}}{144}+\frac{a_{36}}{144}-\frac{5 a_{37}}{144}-\frac{a_{38}}{144}+\frac{a_{39}}{48}\nonumber \\ &+\frac{a_{40}}{48}-\frac{a_{42}}{144}-\frac{5 a_{44}}{144}-\frac{a_{45}}{144}-\frac{a_{47}}{48}-\frac{a_{48}}{144}-\frac{5 a_{53}}{144}-\frac{a_{54}}{144}\nonumber \\ &-\frac{a_{55}}{144}-\frac{5 a_{57}}{144}-\frac{a_{58}}{144}-\frac{a_{61}}{144}+\frac{a_{66}}{48}+\frac{a_{67}}{48}+\frac{a_{68}}{48}+\frac{a_{69}}{144}\nonumber \\ &+\frac{5 a_{70}}{48}+\frac{a_{73}}{2304}+\frac{7 a_{75}}{6912}+\frac{a_{76}}{6912}+\frac{a_{77}}{2304}+\frac{7 a_{79}}{6912}+\frac{a_{80}}{6912}+\frac{a_{81}}{6912}\nonumber \\ &-\frac{a_{82}}{2304}+\frac{a_{83}}{2304}+2 a_{84}+\frac{a_{85}}{3}+2 a_{90}+\frac{a_{91}}{3}+\frac{3 a_{92}}{2}+\frac{a_{93}}{4}\nonumber \\ &+\frac{a_{94}}{3}+\frac{a_{95}}{4}-\frac{a_{98}}{2}-\frac{2 a_{99}}{3}-a_{100}-\frac{4 a_{101}}{3}-a_{102}-6 a_{103}\nonumber \\ &+336 a_{104}-\frac{25}{6912}.
\label{eqn:P1end}
 \end{align}

\section{An Explicit Check of the Topological Recursion Relation}
\label{sec:check}

In this appendix, we perform some extra checks for the coefficients appearing in Lemma 2.1
using the Gromov-Witten invariants of $CP^1$ in order to provide the reader with a better idea of our methods.

Let $\Phi(\tau_{0,0}\tau_{5,0})$ be defined by Equation (\ref{eq:PhiDef}) for the Gromov-Witten theory of $CP^1$. Let $Q$ be the degree $1$ part of $\Phi(\tau_{0,0}\tau_{5,0})$ after setting all variables $t_n^a = 0$ for all $n\geq 0$ and $a=0,1$. We will verify below by direct calculation that the equations in Lemma \ref{lem:coeff} imply that $Q=0$.  We have also performed thousands of  computer verifications for the Gromov-Witten theory of $CP^1$ by showing that the similar restriction of the degree $d$ part of $k$-derivatives of $\Phi(\tau_{m_1 a_1}\tau_{m_2 a_2})$ is $0$ in various cases where $d\leq 4$ and $k\leq 3$.

Setting to zero those correlators in $Q$ which vanish due to dimensional considerations, we obtain
\begin{eqnarray*}
Q&=& a_{82}  \left< \tau_{0,0}^4 \tau_{5,0} \right>_{0,1}  \left< \tau_{0,1} \right>_{1,0}^3+ a_{60}  \left< \tau_{0,0}^4 \tau_{5,0} \right>_{0,1}  \left< \tau_{0,1} \right>_{1,0}  \left< \tau_{0,1}^2 \right>_{1,0}+ a_{95}  \left< \tau_{0,0}^4 \tau_{5,0} \right>_{0,1}  \left< \tau_{0,1}^3 \right>_{1,0}\\
&+& a_{59}  \left< \tau_{0,0}^2 \tau_{0,1}^2 \right>_{0,0}  \left< \tau_{0,1} \right>_{1,0}  \left< \tau_{0,0}^2 \tau_{5,0} \right>_{1,1}+ a_{75}  \left< \tau_{0,0}^2 \tau_{0,1} \right>_{0,0}  \left< \tau_{0,0}^2 \tau_{5,0} \right>_{1,1}  \left< \tau_{0,1} \right>_{1,0}^2\\
&+& a_{79}  \left< \tau_{0,0}^2 \tau_{0,1} \right>_{0,0}  \left< \tau_{0,0}^2 \tau_{5,0} \right>_{1,1}  \left< \tau_{0,1} \right>_{1,0}^2+ a_{47}  \left< \tau_{0,0}^2 \tau_{0,1} \right>_{0,0}  \left< \tau_{0,1}^2 \right>_{1,0}  \left< \tau_{0,0}^2 \tau_{5,0} \right>_{1,1}\\
&+& a_{65}  \left< \tau_{0,0}^2 \tau_{0,1} \right>_{0,0}  \left< \tau_{0,1}^2 \right>_{1,0}  \left< \tau_{0,0}^2 \tau_{5,0} \right>_{1,1}+  a_{92} \{  \left< \tau_{0,1}^3 \right>_{0,1}  \left< \tau_{0,0}^4 \tau_{5,0} \right>_{1,0}
+ 3  \left< \tau_{0,0}^2 \tau_{0,1} \right>_{0,0}  \left< \tau_{0,0}^2 \tau_{0,1}^2 \tau_{5,0} \right>_{1,1} \}\\
&+&  a_{10}  \left< \tau_{0,0}^2 \tau_{0,1} \right>_{0,0}  \left< \tau_{0,1} \right>_{1,0}\{  \left< \tau_{1,0}\tau_{5,0} \right>_{2,1}- \left< \tau_{0,0}^2 \tau_{1,0} \right>_{0,0}  \left< \tau_{0,1}\tau_{5,0} \right>_{2,1}\} \\
&+&  a_{16} \{   \left< \tau_{0,0}^2 \tau_{0,1} \right>_{0,0}\{  \left< \tau_{0,1}\tau_{1,0}\tau_{5,0} \right>_{2,1}- \left< \tau_{0,0}^2 \tau_{1,0} \right>_{0,0}  \left< \tau_{0,1}^2 \tau_{5,0} \right>_{2,1}\}\\
&+&  \left< \tau_{0,0}^2 \tau_{0,1} \right>_{0,0}\{  \left< \tau_{0,0}\tau_{1,1}\tau_{5,0} \right>_{2,1}- \left< \tau_{0,0}\tau_{0,1}\tau_{1,1} \right>_{0,1}  \left< \tau_{0,0}^2 \tau_{5,0} \right>_{2,0}\}
+ a_{9}  \left< \tau_{0,0}^2 \tau_{0,1} \right>_{0,0}  \left< \tau_{0,1} \right>_{1,0}  \left< \tau_{0,0}\tau_{6,0} \right>_{2,1}\\
&+& \left< \tau_{0,0}\tau_{0,1}\tau_{2,0} \right>_{0,1}  \left< \tau_{0,0}\tau_{6,0} \right>_{3,0}+ \left< \tau_{0,0}^2 \tau_{1,0} \right>_{0,0}  \left< \tau_{1,1}\tau_{6,0} \right>_{3,1} \\
&-& \left< \tau_{2,0}\tau_{6,0} \right>_{3,1}- \left< \tau_{0,0}^2 \tau_{1,0} \right>_{0,0}  \left< \tau_{0,0}\tau_{0,1}\tau_{1,1} \right>_{0,1}  \left< \tau_{0,0}\tau_{6,0} \right>_{3,0}\\
&+& 2  a_{15}  \left< \tau_{0,0}^2 \tau_{0,1} \right>_{0,0}  \left< \tau_{0,0}\tau_{0,1}\tau_{6,0} \right>_{2,1} + 2  a_{5}  \left< \tau_{0,0}^2 \tau_{0,1}^2 \right>_{0,0}  \left< \tau_{0,0}\tau_{6,0} \right>_{2,1}
 +  2  a_{2}  \left< \tau_{0,0}^2 \tau_{0,1} \right>_{0,0} \{  \left< \tau_{0,0}\tau_{1,1}\tau_{5,0} \right>_{2,1}\\
 &-& \left< \tau_{0,0}\tau_{0,1}\tau_{1,1} \right>_{0,1}  \left< \tau_{0,0}^2 \tau_{5,0} \right>_{2,0}\}+  2  a_{6}  \left< \tau_{0,0}^2 \tau_{0,1}^2 \right>_{0,0} \{  \left< \tau_{1,0}\tau_{5,0} \right>_{2,1}- \left< \tau_{0,0}^2 \tau_{1,0} \right>_{0,0}  \left< \tau_{0,1}\tau_{5,0} \right>_{2,1}\}\\
&+& 2  a_{24}  \left< \tau_{0,1}\tau_{5,0} \right>_{2,1}  \left< \tau_{0,0}^2 \tau_{0,1} \right>_{0,0}^2 + 4  a_{90}  \left< \tau_{0,0}^2 \tau_{0,1} \right>_{0,0}  \left< \tau_{0,0}^2 \tau_{0,1}^2 \tau_{5,0} \right>_{1,1} + 4  a_{84}  \left< \tau_{0,0}^2 \tau_{0,1} \right>_{0,0}  \left< \tau_{0,0}^2 \tau_{0,1}^2 \tau_{5,0} \right>_{1,1} \\
&+& 2  a_{57}  \left< \tau_{0,0}^2 \tau_{0,1} \right>_{0,0}  \left< \tau_{0,1} \right>_{1,0}  \left< \tau_{0,0}^2 \tau_{0,1}\tau_{5,0} \right>_{1,1} + 2  a_{53}  \left< \tau_{0,0}^2 \tau_{0,1} \right>_{0,0}  \left< \tau_{0,1} \right>_{1,0}  \left< \tau_{0,0}^2 \tau_{0,1}\tau_{5,0} \right>_{1,1} \\
&+& 2  a_{44}  \left< \tau_{0,0}^2 \tau_{0,1} \right>_{0,0}  \left< \tau_{0,1} \right>_{1,0}  \left< \tau_{0,0}^2 \tau_{0,1}\tau_{5,0} \right>_{1,1} + 2  a_{37}  \left< \tau_{0,0}^2 \tau_{0,1} \right>_{0,0}  \left< \tau_{0,1} \right>_{1,0}  \left< \tau_{0,0}^2 \tau_{0,1}\tau_{5,0} \right>_{1,1}\\
& +& 4  a_{96}  \left< \tau_{0,0}^2 \tau_{0,1}^2 \right>_{0,0}  \left< \tau_{0,0}^2 \tau_{0,1}\tau_{5,0} \right>_{1,1} + 3  a_{94}  \left< \tau_{0,0}^2 \tau_{0,1}^2 \right>_{0,0}  \left< \tau_{0,0}^2 \tau_{0,1}\tau_{5,0} \right>_{1,1} + 4  a_{86}  \left< \tau_{0,0}^2 \tau_{0,1}^2 \right>_{0,0}  \left< \tau_{0,0}^2 \tau_{0,1}\tau_{5,0} \right>_{1,1}\\
& +& 2  a_{32}  \left< \tau_{0,0}^2 \tau_{0,1} \right>_{0,0}  \left< \tau_{0,1}^2 \right>_{1,0}  \left< \tau_{0,0}^2 \tau_{5,0} \right>_{1,1} + 2  a_{66}  \left< \tau_{0,0}^2 \tau_{0,1}^2 \right>_{0,0}  \left< \tau_{0,1} \right>_{1,0}  \left< \tau_{0,0}^2 \tau_{5,0} \right>_{1,1}\\
& +& 2  a_{39}  \left< \tau_{0,0}^2 \tau_{0,1}^2 \right>_{0,0}  \left< \tau_{0,1} \right>_{1,0}  \left< \tau_{0,0}^2 \tau_{5,0} \right>_{1,1}
 + 4  a_{100}  \left< \tau_{0,0}^2 \tau_{0,1}^3 \right>_{0,0}  \left< \tau_{0,0}^2 \tau_{5,0} \right>_{1,1} + 2  a_{98}  \left< \tau_{0,0}^2 \tau_{0,1}^3 \right>_{0,0}  \left< \tau_{0,0}^2 \tau_{5,0} \right>_{1,1}\\
& +& 2  a_{99}  \left< \tau_{0,0}^4 \tau_{0,1}\tau_{5,0} \right>_{0,1}  \left< \tau_{0,1}^2 \right>_{1,0} + 2  a_{70}  \left< \tau_{0,0}^4 \tau_{0,1}\tau_{5,0} \right>_{0,1}  \left< \tau_{0,1} \right>_{1,0}^2
+ 4  a_{103}  \left< \tau_{0,0}^4 \tau_{0,1}^2 \tau_{5,0} \right>_{0,1}  \left< \tau_{0,1} \right>_{1,0}\\
& +& 8  a_{104}  \left< \tau_{0,0}^4 \tau_{0,1}^3 \tau_{5,0} \right>_{0,1} .
 \end{eqnarray*}
Applying the string, dilation, and divisor equations, the previous equation reduces to
\begin{eqnarray*}
Q&=&   -\frac{1}{12} a_{37} \{   \left< \tau_{3,0} \right>_{1,1} +  \left< \tau_{2,1} \right>_{1,1} \}  -\frac{1}{12} a_{44} \{   \left< \tau_{3,0} \right>_{1,1} +  \left< \tau_{2,1} \right>_{1,1} \}-\frac{1}{12}  a_{53} \{    \left< \tau_{3,0} \right>_{1,1} +  \left< \tau_{2,1} \right>_{1,1} \}\\
&-& \frac{1}{12} a_{57} \{    \left< \tau_{3,0} \right>_{1,1} +  \left< \tau_{2,1} \right>_{1,1} \}+  a_{92} \{ \frac{1}{12}+ 3  \left< \tau_{3,0} \right>_{1,1} + 6  \left< \tau_{2,1} \right>_{1,1} \}
+  a_{84} \{  4  \left< \tau_{3,0} \right>_{1,1} \\
&+& 8  \left< \tau_{2,1} \right>_{1,1} \}+  a_{90} \{  4  \left< \tau_{3,0} \right>_{1,1} + 8  \left< \tau_{2,1} \right>_{1,1} \}+  a_{15} \{  2  \left< \tau_{5,0} \right>_{2,1} + 2  \left< \tau_{4,1} \right>_{2,1} \}+  a_{24} \{  2  \left< \tau_{5,0} \right>_{2,1} + 2  \left< \tau_{4,1} \right>_{2,1} \}\\
&+&  a_{2} \{  2  \left< \tau_{5,0} \right>_{2,1} + 2  \left< \tau_{4,1} \right>_{2,1} + 2  \left< \tau_{1,1}\tau_{4,0} \right>_{2,1} - 2  \left< \tau_{3,0} \right>_{2,0} \}
+  a_{16} \{  \left< \tau_{1,1}\tau_{4,0} \right>_{2,1}
+ 5  \left< \tau_{5,0} \right>_{2,1} + 5  \left< \tau_{4,1} \right>_{2,1}\\
&-& \left< \tau_{3,0} \right>_{2,0}\}- \left< \tau_{2,0}\tau_{6,0} \right>_{3,1}- \left< \tau_{5,0} \right>_{3,0}
- \frac{1}{8}  a_{10}  \left< \tau_{5,0} \right>_{2,1} - \frac{1}{24}  a_{9}  \left< \tau_{5,0} \right>_{2,1} - \frac{1}{288}  a_{70}
  + \frac{1}{576}  a_{79}  \left< \tau_{3,0} \right>_{1,1}\\
 &+& \frac{1}{576}  a_{75}  \left< \tau_{3,0} \right>_{1,1}
 + 8  a_{104}  + \frac{1}{6912}  a_{82}.
\end{eqnarray*}
Plugging in the values for the correlators using Gathmann's program based on the Virasoro constraints, c.f. the Appendix of \cite{KL}, we obtain
\begin{eqnarray*}
Q&=&-\frac{757}{1451520}- \frac{1}{288}  a_{70}  - \frac{1}{288}  a_{57}  - \frac{1}{288}  a_{53}  - \frac{1}{288}  a_{44}  - \frac{1}{288}  a_{37}  - \frac{1}{4608}  a_{10}  - \frac{1}{13824}  a_{9}  + 8  a_{104}  + \frac{1}{3}  a_{92}  \\
&+& \frac{1}{3}  a_{90}  + \frac{1}{3}  a_{84}  + \frac{1}{6912}  a_{82}  + \frac{13}{2880}  a_{24}  + \frac{13}{720}  a_{16}  + \frac{13}{2880}  a_{15}  + \frac{13}{720}  a_{2}
\end{eqnarray*}
Plugging in the equations from Lemma \ref{lem:coeff} implies that $Q=0$ for all $a_2$.

\section{A relation in the tautological ring of $\overline{\cal M}_{3,2}$}
\label{sec:dualgraph}

In this appendix, we give the dual graph representation
for the relation in the tautological ring of $\overline{\cal M}_{3,2}$
corresponding to Theorem~\ref{thm:g3T2Tskew}.
 We adopt the conventions of \cite{Ge2} for dual graphs with a slight
modification. We denote vertices of genus $0$  by a hollow circle
$\begin{picture}(10,8)(0,0)\put(5,3){\circle{7}}\end{picture}$, and
vertices of genus $g \geq 1$ by $\begin{picture}(10,8)(0,0)\put(5,3){\circle{7}}
\put(3,2){$\scriptscriptstyle g$}\end{picture}$.
A vertex with an incident arrowhead denotes the $\psi$ class associated to
the marked point (or a node) on the irreducible component associated
to that vertex. When translating relations in the tautological ring of $\overline{\cal M}_{g, n}$
to universal equations for Gromov-Witten invariants, we need to divide the coefficient
of each stratum by the number of
elements in the automorphism group of the corresponding dual graph.

Let
{\allowdisplaybreaks
\begin{eqnarray}
 G_{i, j}
&:=&   \frac{4}{9} \hspace{10pt}
\begin{picture}(80, 30)
\put(8, 16){$\scriptscriptstyle j$}
\put(7, 5){\line(0, 1){10}}
\put(5.5, -1){$\scriptscriptstyle 2$}
\put(7, 1){\circle{7}}
\put(40.5, 2){\vector(-1, 0){30}}
\put(40, 1){ \circle{7}}
\put(44, 16){$\scriptscriptstyle i$}
\put(43, 5){\line(0, 1){10}}
\curve(48, 2, 60, 9, 70, 2, 60, -7, 48, 0)
\end{picture}
+ \frac{5}{12} \hspace{10pt}
\begin{picture}(80, 30)
\put(8, 16){$\scriptscriptstyle j$}
\put(7, 15){\vector(0, -1){10}}
\put(5.5, -1){$\scriptscriptstyle 2$}
\put(7, 1){\circle{7}}
\put(40.5, 2){\line(-1, 0){30}}
\put(40, 1){ \circle{7}}
\put(44, 16){$\scriptscriptstyle i$}
\put(43, 5){\line(0, 1){10}}
\curve(48, 2, 60, 9, 70, 2, 60, -7, 48, 0)
\end{picture}
 \nonumber \\
&& +\frac{16}{3} \hspace{10pt}
\begin{picture}(100, 30)
\put(8, 16){$\scriptscriptstyle j$}
\put(7, 5){\line(0, 1){10}}
\put(5.5, -1){$\scriptscriptstyle 2$}
\put(7, 1){\circle{7}}
\put(40.5, 2){\vector(-1, 0){30}}
\put(44, 16){$\scriptscriptstyle i$}
\put(40, 1){ \circle{7}}
\put(43, 5){\line(0, 1){10}}
\put(48, 1){\line(1, 0){30}}
\put(81, 1){\circle{7}}
\put(79.5, -1){$\scriptscriptstyle 1$}
\end{picture}
+5 \hspace{10pt}
\begin{picture}(100, 30)
\put(9, 16){$\scriptscriptstyle j$}
\put(7, 15){\vector(0, -1){10}}
\put(5.5, -1){$\scriptscriptstyle 2$}
\put(7, 1){\circle{7}}
\put(40.5, 2){\line(-1, 0){30}}
\put(44, 16){$\scriptscriptstyle i$}
\put(40, 1){ \circle{7}}
\put(43, 5){\line(0, 1){10}}
\put(48, 1){\line(1, 0){30}}
\put(81, 1){\circle{7}}
\put(79.5, -1){$\scriptscriptstyle 1$}
\end{picture}
\nonumber \\
&&
+\frac{40}{3} \hspace{10pt}
\begin{picture}(100, 30)
\put(5.5, -1){$\scriptscriptstyle 2$}
\put(7, 1){\circle{7}}
\put(40.5, 2){\vector(-1, 0){30}}
\put(44, 16){$\scriptscriptstyle i$}
\put(40, 1){ \circle{7}}
\put(43, 5){\line(0, 1){10}}
\put(48, 1){\line(1, 0){30}}
\put(83, 16){$\scriptscriptstyle j$}
\put(82, 5){\line(0, 1){10}}
\put(81, 1){\circle{7}}
\put(79.5, -1){$\scriptscriptstyle 1$}
\end{picture}
 +\frac{1}{6} \hspace{10pt}
\begin{picture}(80, 30)
\put(8, 16){$\scriptscriptstyle j$}
\put(7, 5){\line(0, 1){10}}
\put(5.5, -1){$\scriptscriptstyle 2$}
\put(7, 1){\circle{7}}
\put(40.5, 2){\vector(-1, 0){30}}
\put(44, 16){$\scriptscriptstyle i$}
\put(40, 1){ \circle{7}}
\curve(11, -1, 26, -6, 41, -1)
\put(43, 5){\line(0, 1){10}}
\end{picture}
\nonumber \\
&&
+ \hspace{10pt}
\begin{picture}(80, 30)
\put(8, 16){$\scriptscriptstyle j$}
\put(7, 15){\vector(0, -1){10}}
\put(5.5, -1){$\scriptscriptstyle 2$}
\put(7, 1){\circle{7}}
\put(40.5, 2){\line(-1, 0){30}}
\put(44, 16){$\scriptscriptstyle i$}
\put(40, 1){ \circle{7}}
\curve(11, -1, 26, -6, 41, -1)
\put(43, 5){\line(0, 1){10}}
\end{picture}
+ \frac{4}{9} \hspace{10pt}
\begin{picture}(120, 30)
\put(8, 16){$\scriptscriptstyle i$}
\put(7, 5){\line(0, 1){10}}
\put(5.5, -1){$\scriptscriptstyle 2$}
\put(7, 1){\circle{7}}
\put(10.5, 2){\line(1, 0){30}}
\put(44, 16){$\scriptscriptstyle j$}
\put(40, 1){ \circle{7}}
\put(43, 5){\line(0, 1){10}}
\put(48, 1){\line(1, 0){30}}
\put(81, 1){\circle{7}}
\curve(84.5, 2, 97, 9, 107, 2, 97, -7, 84.5, 0)
\end{picture} \nonumber \\
&& + \frac{1}{9} \hspace{10pt}
\begin{picture}(120, 30)
\put(5.5, -1){$\scriptscriptstyle 1$}
\put(8, 16){$\scriptscriptstyle i$}
\put(7, 5){\line(0, 1){10}}
\put(7, 1){\circle{7}}
\put(10.5, 2){\line(1, 0){30}}
\put(40, 1){ \circle{7}}
\put(42.5, -1){$\scriptscriptstyle 1$}
\put(44, 16){$\scriptscriptstyle j$}
\put(43, 5){\line(0, 1){10}}
\put(48, 1){\line(1, 0){30}}
\put(81, 1){\circle{7}}
\curve(85, 2, 97, 9, 107, 2, 97, -7, 85, 0)
\end{picture}
+ \frac{1}{9} \hspace{10pt}
\begin{picture}(120, 30)
\put(5.5, -1){$\scriptscriptstyle 1$}
\put(8, 16){$\scriptscriptstyle i$}
\put(7, 5){\line(0, 1){10}}
\put(7, 1){\circle{7}}
\put(10.5, 2){\line(1, 0){30}}
\put(40, 1){ \circle{7}}
\put(42.5, -1){$\scriptscriptstyle 1$}
\put(48, 1){\line(1, 0){30}}
\put(81, 1){\circle{7}}
\put(84, 16){$\scriptscriptstyle j$}
\put(82, 5){\line(0, 1){10}}
\curve(85, 2, 97, 9, 107, 2, 97, -7, 85, 0)
\end{picture}  \nonumber \\
&& + \frac{1}{15} \hspace{10pt}
\begin{picture}(120, 30)
\put(5.5, -1){$\scriptscriptstyle 1$}
\put(7, 1){\circle{7}}
\put(10.5, 2){\line(1, 0){30}}
\put(40, 1){ \circle{7}}
\put(42.5, -1){$\scriptscriptstyle 1$}
\put(44, 16){$\scriptscriptstyle j$}
\put(43, 5){\line(0, 1){10}}
\put(48, 1){\line(1, 0){30}}
\put(81, 1){\circle{7}}
\put(83, 16){$\scriptscriptstyle i$}
\put(82, 5){\line(0, 1){10}}
\curve(85, 2, 97, 9, 107, 2, 97, -7, 85, 0)
\end{picture}
+ \frac{1}{15} \hspace{10pt}
\begin{picture}(100, 30)
\put(5.5, -1){$\scriptscriptstyle 1$}
\put(7, 1){\circle{7}}
\put(10.5, 2){\line(1, 0){30}}
\put(43.5, -1){$\scriptscriptstyle 1$}
\put(44, 1){\circle{7}}
\put(44, 16){$\scriptscriptstyle j$}
\put(43, 5){\line(0, 1){10}}
\curve(48, 2, 66, 6, 81, 2)
\curve(48, -2, 66, -6, 81, -2)
\put(80, 1){ \circle{7}}
\put(84, 16){$\scriptscriptstyle i$}
\put(83, 5){\line(0, 1){10}}
\end{picture} \nonumber \\
&& +\frac{9}{10} \hspace{10pt}
\begin{picture}(100, 30)
\put(5.5, -1){$\scriptscriptstyle 1$}
\put(7, 1){\circle{7}}
\put(8, 16){$\scriptscriptstyle j$}
\put(7, 5){\line(0, 1){10}}
\put(10.5, 2){\line(1, 0){30}}
\put(40, 1){ \circle{7}}
\put(44, 16){$\scriptscriptstyle i$}
\put(43, 5){\line(0, 1){10}}
\put(48, 1){\line(1, 0){30}}
\put(81, 1){\circle{7}}
\put(79.5, -1){$\scriptscriptstyle 1$}
\curve(7, -2, 22, -8, 44, -10, 66, -8, 81, -2)
\end{picture}
 + \frac{1}{15} \hspace{10pt}
\begin{picture}(120, 30)
\put(5.5, -1){$\scriptscriptstyle 1$}
\put(7, 1){\circle{7}}
\put(10.5, 2){\line(1, 0){30}}
\put(40, 1){ \circle{7}}
\put(44, 16){$\scriptscriptstyle i$}
\put(43, 5){\line(0, 1){10}}
\put(48, 1){\line(1, 0){30}}
\put(81, 1){\circle{7}}
\put(83, 16){$\scriptscriptstyle j$}
\put(82, 5){\line(0, 1){10}}
\put(80.5, -1){$\scriptscriptstyle 1$}
\curve(85, 2, 97, 9, 107, 2, 97, -7, 85, 0)
\end{picture}
\nonumber \\
&&
+ \frac{4}{15} \hspace{10pt}
\begin{picture}(120, 30)
\put(5.5, -1){$\scriptscriptstyle 1$}
\put(8, 16){$\scriptscriptstyle i$}
\put(7, 5){\line(0, 1){10}}
\put(7, 1){\circle{7}}
\put(10.5, 2){\line(1, 0){30}}
\put(40, 1){ \circle{7}}
\put(44, 16){$\scriptscriptstyle j$}
\put(43, 5){\line(0, 1){10}}
\put(48, 1){\line(1, 0){30}}
\put(81, 1){\circle{7}}
\put(80.5, -1){$\scriptscriptstyle 1$}
\curve(85, 2, 97, 9, 107, 2, 97, -7, 85, 0)
\end{picture}
+ \frac{6}{5} \hspace{10pt}
\begin{picture}(100, 30)
\put(7, 16){$\scriptscriptstyle j$}
\put(5.5, 5){\line(0, 1){10}}
\put(5.5, -1){$\scriptscriptstyle 1$}
\put(7, 1){\circle{7}}
\curve(11, 2, 27, 6, 44, 2)
\curve(11, -2, 27, -6, 44, -2)
\put(43, 1){ \circle{7}}
\put(47, 16){$\scriptscriptstyle i$}
\put(46, 5){\line(0, 1){10}}
\put(50.5, 2){\line(1, 0){30}}
\put(83.5, -1){$\scriptscriptstyle 1$}
\put(85, 1){\circle{7}}
\end{picture}
 \nonumber \\
&&
+ \frac{1}{3} \hspace{10pt}
\begin{picture}(100, 30)
\put(7, 16){$\scriptscriptstyle i$}
\put(5.5, 5){\line(0, 1){10}}
\put(5.5, -1){$\scriptscriptstyle 1$}
\put(7, 1){\circle{7}}
\curve(11, 2, 27, 6, 44, 2)
\curve(11, -2, 27, -6, 44, -2)
\put(43, 1){ \circle{7}}
\put(50.5, 2){\line(1, 0){30}}
\put(83.5, -1){$\scriptscriptstyle 1$}
\put(85, 1){\circle{7}}
\put(86, 16){$\scriptscriptstyle j$}
\put(85, 5){\line(0, 1){10}}
\end{picture}
+ \frac{32}{15} \hspace{10pt}
\begin{picture}(100, 30)
\put(5.5, -1){$\scriptscriptstyle 1$}
\put(7, 1){\circle{7}}
\curve(11, 2, 27, 6, 44, 2)
\curve(11, -2, 27, -6, 44, -2)
\put(43, 1){ \circle{7}}
\put(47, 16){$\scriptscriptstyle i$}
\put(46, 5){\line(0, 1){10}}
\put(50.5, 2){\line(1, 0){30}}
\put(83.5, -1){$\scriptscriptstyle 1$}
\put(85, 1){\circle{7}}
\put(87, 16){$\scriptscriptstyle j$}
\put(85.5, 5){\line(0, 1){10}}
\end{picture}
 \nonumber \\
&&
+ \frac{41}{45} \hspace{10pt}
\begin{picture}(100, 30)
\put(8.5, 16){$\scriptscriptstyle j$}
\put(7, 5){\line(0, 1){10}}
\put(5.5, -1){$\scriptscriptstyle 1$}
\put(7, 1){\circle{7}}
\put(10.5, 2){\line(1, 0){30}}
\put(40, 1){ \circle{7}}
\put(45, 16){$\scriptscriptstyle i$}
\put(44, 5){\line(0, 1){10}}
\curve(41, -2, 38, -8, 44, -16, 50, -8, 47, -2)
\put(48, 1){\line(1, 0){30}}
\put(81, 1){\circle{7}}
\put(79.5, -1){$\scriptscriptstyle 1$}
\end{picture}
+ \frac{4}{5} \hspace{10pt}
\begin{picture}(100, 30)
\put(5.5, -1){$\scriptscriptstyle 1$}
\put(7, 1){\circle{7}}
\put(10.5, 2){\line(1, 0){30}}
\put(42.5, -1){$\scriptscriptstyle 1$}
\put(44, 1){\circle{7}}
\put(45, 16){$\scriptscriptstyle j$}
\put(44, 5){\line(0, 1){10}}
\put(47, 2){\line(1, 0){29}}
\put(76, 1){ \circle{7}}
\put(81, 16){$\scriptscriptstyle i$}
\put(80, 5){\line(0, 1){10}}
\put(84, 1){\line(1, 0){30}}
\put(118, 1){\circle{7}}
\put(116.5, -1){$\scriptscriptstyle 1$}
\end{picture}
 \nonumber \\
&&
+ \frac{16}{5} \hspace{10pt}
\begin{picture}(100, 30)
\put(5.5, -1){$\scriptscriptstyle 1$}
\put(7, 1){\circle{7}}
\put(10.5, 2){\line(1, 0){30}}
\put(42.5, -1){$\scriptscriptstyle 1$}
\put(44, 1){\circle{7}}
\put(47, 2){\line(1, 0){29}}
\put(76, 1){ \circle{7}}
\put(81, 16){$\scriptscriptstyle j$}
\put(80, 5){\line(0, 1){10}}
\put(84, 1){\line(1, 0){30}}
\put(118, 1){\circle{7}}
\put(120, 16){$\scriptscriptstyle i$}
\put(119, 5){\line(0, 1){10}}
\put(116.5, -1){$\scriptscriptstyle 1$}
\end{picture}
\nonumber \\
&&
+ \frac{184}{15} \hspace{10pt}
\begin{picture}(100, 30)
\put(8.5, 16){$\scriptscriptstyle j$}
\put(7, 5){\line(0, 1){10}}
\put(5.5, -1){$\scriptscriptstyle 1$}
\put(7, 1){\circle{7}}
\put(10.5, 2){\line(1, 0){30}}
\put(40, 1){ \circle{7}}
\put(45, 16){$\scriptscriptstyle i$}
\put(44, 5){\line(0, 1){10}}
\put(43.5, -2.5){\line(0, -1){20}}
\put(43.5, -26.5){\circle{7}}
\put(42.5, -28.5){$\scriptscriptstyle 1$}
\put(48, 1){\line(1, 0){30}}
\put(81, 1){\circle{7}}
\put(79.5, -1){$\scriptscriptstyle 1$}
\end{picture}
+ \frac{1}{180} \hspace{10pt}
\begin{picture}(110, 30)
\curve(28.5, 2, 15, 8, 2, 0, 15, -8, 30, -2)
\put(32.5, 1){\circle{7}}
\put(31, -1){$\scriptscriptstyle 1$}
\put(36, 2){\line(1, 0){30}}
\put(33.5, 16){$\scriptscriptstyle j$}
\put(32, 5){\line(0, 1){10}}
\put(66, 1){ \circle{7}}
\put(70, 16){$\scriptscriptstyle i$}
\put(69, 5){\line(0, 1){10}}
\curve(73, 2, 87, 9, 98, 2, 87, -7, 73, 0)
\end{picture}
 \nonumber \\
&&
+ \frac{1}{180} \hspace{10pt}
\begin{picture}(100, 50)
\curve(28.5, 2, 15, 8, 2, 0, 15, -8, 30, -2)
\put(31, -1){$\scriptscriptstyle 1$}
\put(32.5, 1){\circle{7}}
\put(33.5, 16){$\scriptscriptstyle j$}
\put(32, 5){\line(0, 1){10}}
\curve(37, 2, 56, 6, 71, 2)
\curve(37, -2, 56, -6, 71, -2)
\put(70, 1){ \circle{7}}
\put(74, 16){$\scriptscriptstyle i$}
\put(73, 5){\line(0, 1){10}}
\end{picture}
+ \frac{2}{45} \hspace{10pt}
\begin{picture}(60, 50)
\put(5.5, -1){$\scriptscriptstyle 1$}
\put(7, 1){\circle{7}}
\put(8.5, 16){$\scriptscriptstyle j$}
\put(7, 5){\line(0, 1){10}}
\curve(11, 2, 27, 6, 44, 2)
\put(10, 0){\line(1, 0){33}}
\curve(11, -2, 27, -6, 44, -2)
\put(43, 1){ \circle{7}}
\put(47, 16){$\scriptscriptstyle i$}
\put(46, 5){\line(0, 1){10}}
\end{picture}
\nonumber \\
&&
+ \frac{1}{10} \hspace{10pt}
\begin{picture}(100, 30)
\curve(28.5, 2, 15, 8, 2, 0, 15, -8, 30, -2)
\put(32.5, 1){\circle{7}}
\put(33.5, 16){$\scriptscriptstyle i$}
\put(32, 5){\line(0, 1){10}}
\curve(36, 2, 56, 6, 71, 2)
\curve(35, -2, 56, -6, 71, -2)
\put(70, 1){ \circle{7}}
\put(72.5, -1){$\scriptscriptstyle 1$}
\put(74.5, 16){$\scriptscriptstyle j$}
\put(73, 5){\line(0, 1){10}}
\end{picture}
+ \frac{1}{15} \hspace{10pt}
\begin{picture}(80, 30)
\put(8.5, 16){$\scriptscriptstyle j$}
\put(7, 5){\line(0, 1){10}}
\put(5.5, -1){$\scriptscriptstyle 1$}
\put(7, 1){\circle{7}}
\put(10.5, 2){\line(1, 0){30}}
\put(40, 1){ \circle{7}}
\put(40.5, 16){$\scriptscriptstyle i$}
\put(43, 5){\line(0, 1){10}}
\curve(44.5, 5, 52, 14, 64, 18, 54, 6, 47, 3)
\curve(43, -3, 52, -14, 64, -18, 54, -6, 46, -2)
\end{picture}
\nonumber
\end{eqnarray}}
where $i$ and $j$ label the two marked points (i.e. tails of dual graphs).
Then Theorem~\ref{thm:g3T2Tskew} corresponds to the following
\begin{thm}
In the tautological ring of $\overline{\cal M}_{3, 2}$, the following relation holds
\[ \psi_1^2 \psi_2 - \psi_1 \psi_2^2 \,\, =\, \frac{1}{7} \, \left\{ G_{1,2} - G_{2,1} \right\}. \]
\end{thm}



\vspace{30pt}
\noindent
Takashi Kimura \\

\noindent
Department of Mathematics and Statistics; \\
111 Cummington Mall, Boston University; \\
Boston, MA 02215, USA.  \\
\noindent
Email: {\it kimura@math.bu.edu}

\vspace{40pt}

\noindent
Xiaobo Liu \\

\noindent
Beijing International Center for Mathematical Research, \\
Beijing University, Beijing, China. \\
Email: {\it xbliu@math.pku.edu.cn}

\noindent
\&

\noindent
Department of Mathematics, \\
University of Notre Dame, \\
Notre Dame, IN 46556, USA \\
\noindent
Email: {\it xliu3@nd.edu}

\end{document}